\documentclass[11pt]{article}
\usepackage[english]{babel}
\usepackage{amsmath,amsthm,amssymb}
\usepackage{mathrsfs}
\usepackage{bbm}
\usepackage{amsfonts}
\usepackage[margin=0.8in]{geometry}
\usepackage{graphicx}
\usepackage{caption}
\usepackage{subcaption}
\usepackage{wasysym}
\usepackage{enumerate}
\usepackage{algorithm2e}
\usepackage{tabularx}
\usepackage{color}
\usepackage{wrapfig}
\usepackage{todonotes}

\newtheorem{thm}{Theorem}[section]
\newtheorem{ass}[thm]{Assumption}
\newtheorem{cor}[thm]{Corollary}
\newtheorem{lem}[thm]{Lemma}
\newtheorem{prop}[thm]{Proposition}
\newtheorem*{hyp*}{Hypothesis}
\theoremstyle{definition}
\newtheorem{defn}[thm]{Definition}

\theoremstyle{remark}
\newtheorem{rem}[thm]{Remark}
\numberwithin{equation}{section}
\newcommand{\R}{\mathbb R}

\newcommand{\bbD}{\mathbb D}
\newcommand{\bbF}{\mathbb F}

\newcommand{\bbQ}{\mathbb Q}

\newcommand{\mcA}{\mathcal{A}}

\newcommand{\mcB}{\mathcal{B}}
\newcommand{\mcD}{\mathcal D}
\newcommand{\mcE}{\mathcal E}

\newcommand{\mcF}{\mathcal F}
\newcommand{\mcI}{\mathcal I}

\newcommand{\mcT}{\mathcal T}
\newcommand{\mcP}{\mathcal P}

\newcommand{\mcO}{\mathcal O}
\newcommand{\mcH}{\mathcal H}
\newcommand{\mcK}{\mathcal K}
\newcommand{\mcU}{\mathcal U}
\newcommand{\mcS}{\mathcal S}

\newcommand{\E}{\mathbb{E}}
\newcommand{\Prob}{\mathbb{P}}

\newcommand{\vect}{\mathbf{t}}
\newcommand{\vecb}{\mathbf{b}}

\newcommand{\vecv}{\mathbf{v}}
\newcommand{\esssup}{\mathop{\rm{ess}\,\sup}}

\newcommand{\ett}{\mathbbm{1}}

\newcommand{\cadlag}{c\`adl\`ag~}

\newcommand{\caglad}{c\`agl\`ad~}

\newcommand{\PrM}{\mathfrak{P}}%{{\rm{Pr}}}

\newcommand{\ie}{\textit{i.e.\ }}
\newcommand{\eg}{\textit{e.g.\ }}
\newcommand{\etal}{\textit{et.~al.\ }}

\begin{document}

%\title{Sequential Systems of Reflected Backward Stochastic Differential Equations with Application to Mixed Continuous and Impulse Control in Infinite Horizon\footnote{This work was supported by the Swedish Energy Agency through grant number 42982-1}}

\title{Finite Horizon Robust Impulse Control in a Non-Markovian Framework and Related Systems of Reflected BSDEs\footnote{This work was supported by the Swedish Energy Agency through grant number 42982-1}}

\author{Magnus Perninge\footnote{M.\ Perninge is with the Department of Physics and Electrical Engineering, Linnaeus University, V\"axj\"o,
Sweden. e-mail: magnus.perninge@lnu.se.}} %
\maketitle
% ----------------------------------------------------------------
\begin{abstract}
We consider a robust impulse control problem in finite horizon where the underlying uncertainty stems from an impulsively and continuously controlled functional stochastic differential equation (FSDE) driven by Brownian motion. We assume that the controller acts upon the system by impulses while the adversary player (nature) acts through continuous controls. We look for a weak solution which leads us to consider a system of sequentially interconnected, obliquely reflected backward stochastic differential equations (RBSDEs) with stochastic Lipschitz coefficients. We show existence of solutions to our system of RBSDEs by applying a Picard iteration approach. Uniqueness then follows by relating the limit to an auxiliary impulse control problem.
\end{abstract}

% ----------------------------------------------------------------
\section{Introduction}
The standard stochastic impulse control problem is an optimal control problem that arises when an operator controls a dynamical system with noise by intervening on the system at a discrete set of stopping times. In impulse control the control-law, thus, takes the form $u=(\tau_1,\ldots,\tau_N;\beta_1,\ldots,\beta_N)$, where $\tau_1\leq\tau_2\leq\cdots\leq\tau_N$ is a sequence of times when the operator intervenes on the system and $\beta_j$ is the impulse that the operator affects the system with at time $\tau_j$. We restrict our attention to the case of a Brownian filtration $\bbF:=\{\mcF_t\}_{t\geq 0}$ and assume that the $\tau_j$ are $\bbF$-stopping times and that $\beta_j$ is $\mcF_{\tau_j}$-measurable and take values in a compact subset $U$ of $\R^d$.

As impulse control problems appear in a vast number of real-world applications (see \eg~\cite{Korn,PSImpulsive} for applications in finance and \cite{BaseiImpulse,CarmLud} for applications in energy) a lot of attention has been given to various types of problems where the control is of impulse type. In the standard Markovian setting the relation to quasi-variational inequalities has frequently been exploited to find optimal controls (see the seminal work in~\cite{BensLionsImpulse} or turn to \cite{OksenSulemBok} for a more recent textbook). In the non-Markovian framework an impulse control problem where the impulses affect the dynamics in an additive way was solved in~\cite{DjehiceImpulse}. This approach was extended to incorporate delivery lag in~\cite{Hdhiri} and, more recently, also to an infinite horizon setting in~\cite{DjehicheInfHorImp}.

A different approach to non-Markovian impulse control was initiated in \cite{SwitchElephant} and then further developed in \cite{JonteSFDE} where interconnected Snell envelopes indexed by controls were used to find solutions to problems with impulsively controlled functional stochastic differential equations (FSDEs) under additional $L^p$-type Lipschitz conditions on the coefficients. We mention also the extension to infinite horizon in~\cite{ImpulsIH}.

In many applications the controller is ambiguous concerning the probabilistic model for the driving noise or unsure of the accuracy of the dynamic model employed. One way of incorporating such ambiguity into the model itself goes through robust formulations of the control problem. We consider %the situation when the ambiguity is incorporated into the model by having an
a robust formulation of the problem where an adversary player (nature) chose a worst case continuous control in response to the operators impulse control, also referred to as Knightian uncertainty. As opposed to \cite{BayraktarRobust}, which considers robust optimal switching in a Markovian framework, we consider more general impulse control problems and a non-Markovian setting and utilize the theory of backward stochastic differential equations (BSDEs). In particular, we extend the results in \cite{JonteSFDE} to a robust framework by considering weak formulations to the problem of maximizing the reward functional
\begin{align}
J(u,\alpha):=\E\Big[\int_0^T\phi(t,X^{u,\alpha}_t,\alpha_t)dt+\psi(X^{u,\alpha}_T) -\sum_{j=1}^N\ell(\tau_j,X^{[u]_{j-1},\alpha}_{\tau_j},\beta_j)\Big]
\end{align}
over impulse controls, $u$, when simultaneously a minimization is performed over continuous controls $\alpha:=(\alpha_s)_{0\leq s\leq T}$, taking values in a compact subset $A$ of $\R^d$. Here, $[u]_j:=(\tau_1,\ldots,\tau_{N\wedge j};\beta_1,\ldots,\beta_{N\wedge j})$ and $X^{u,\alpha}$ solves the impulsively controlled functional SDE
\begin{align}
X^{u,\alpha}_t&=x_0+\int_0^ta(s,(X^{u,\alpha}_r)_{ r\leq s},\alpha_s)ds+\int_0^t\sigma(s,(X^{u,\alpha}_r)_{ r\leq s})dW_s\label{ekv:forward-sde1}
\end{align}
for $t\in [0,\tau_{1})$ and
\begin{align}
X^{u,\alpha}_{t}&=\Gamma(\tau_j,X^{[u]_{j-1},\alpha}_{\tau_j},\beta_j)+\int_{\tau_j}^ta(s,(X^{u,\alpha}_r)_{ r\leq s},\alpha_s)ds+\int_{\tau_j}^t\sigma(s,(X^{u,\alpha}_r)_{ r\leq s})dW_s,\label{ekv:forward-sde2}
\end{align}
when $t\in [\tau_{j},\tau_{j+1})$ for $j=1,\ldots,N$, with $\tau_{N+1}:=\infty$. By considering systems of reflected BSDEs with drivers that satisfy a stochastic Lipschitz condition we are able to relax the common assumption that $|\sigma^{-1}(t,x)a(t,x,\alpha)|$ is bounded and instead assume a linear growth, \ie that $|\sigma^{-1}(t,x)a(t,x,\alpha)|\leq k_L(1+\sup_{s\in[0,t]}|x_s|)$, for some constant $k_L>0$.

El Karoui \etal introduced the notion of reflected backward stochastic differential equations (RBSDEs) and demonstrated a link between RBSDEs and optimal stopping in~\cite{ElKaroui1}. This was later exploited in a series of articles~\cite{HamJean,HuTang,HamZhang} as a means of finding solutions to optimal switching problems. In \cite{HamZhang} existence and uniqueness of solutions to an interconnected systems of reflected BSDEs were shown. Furthermore, it was shown that the solutions are related to optimal switching problems under Knightian uncertainty. Important contributions from the perspective of the present work are also \cite{Bender00,Briand08} that consider BSDEs where the Lipschitz coefficient on the $z$-variable of the driver is a stochastic process and the more recent work presented in~\cite{ElAsri2020} where a RBSDE with stochastic Lipschitz coefficient is solved.

As a means of finding weak solutions to the above mentioned robust impulse control problem we aim to add to the results on RBSDEs by finding solutions to the system of RBSDEs
\begin{align}\label{ekv:seq-bsde-FH}
\begin{cases}
  Y^{v}_t=\xi^v+\int_t^T f^v(s,Y^{v}_s,Z^{v}_s)ds-\int_t^T Z^{v}_sdW_s+ K^{v}_T-K^{v}_t,\quad\forall t\in[0,T], \\
  Y^{v}_t\geq\sup_{b\in U}\{Y^{v\circ(t,b)}_t-c^v(t,b)\},\quad\forall t\in[0,T],\\
  \int_0^T(Y^{v}_t-\sup_{b\in U}\{Y^{v\circ(t,b)}_t-c^v(t,b)\})dK^{v}_t=0
\end{cases}
\end{align}
indexed by impulse controls $v$, under the assumption that the driver $f$ satisfies a stochastic Lipschitz condition on the $z$-variable. We rely on a Piccard iteration approach and the main obstacle we face is showing continuity of the map $(t,b)\mapsto Y^{v\circ(t,b)}_t$.

The remainder of the article is organised as follows. In the next section we set the notation, specify what we mean by a solution to \eqref{ekv:seq-bsde-FH} and recall some well known results for RBSDEs with deterministic Lipschitz coefficients from the original work~\cite{ElKaroui1}. Then in Section~\ref{sec:simple-rbsde} we derive moment and stability estimates for solutions to RBSDEs with stochastic Lipschitz coefficients. In Section~\ref{sec:seq-rbsde-FH} we turn to sequential systems of RBSDEs and show that \eqref{ekv:seq-bsde-FH} admits a unique solution. Finally, in Section~\ref{sec:robust-impulse} we show how to find optimal controls for our robust impulse control problem by relating solutions to \eqref{ekv:seq-bsde-FH} to weak formulations of the control problem at hand.

\section{Preliminaries\label{sec:prelim}}
We let $(\Omega,\mcF,\bbF,\Prob)$ be a complete filtered probability space, where $\bbF:=(\mcF_t)_{0\leq t\leq T}$ is the augmented natural filtration of a $d$-dimensional Brownian motion $W$ and $\mcF:=\mcF_T$, where $T\in(0,\infty)$ is the horizon.\\

\noindent Throughout, we will use the following notation:
\begin{itemize}
  \item We let $\E$ denote expectation with respect to $\Prob$ and for any other probability measure $\bbQ$ on $(\Omega,\mcF)$, we denote by $\E^{\bbQ}$ expectation with respect to $\bbQ$.
  \item $\mcP_{\bbF}$ is the $\sigma$-algebra of $\bbF$-progressively measurable subsets of $[0,T]\times \Omega$.
  \item For $p\geq 1$, we let $\mcS^{p}$ be the set of all $\R$-valued, $\mcP_{\bbF}$-measurable, continuous processes $(Z_t: t\in[0,T])$ such that $\|Z\|_{\mcS^p}:=\E\big[\sup_{t\in[0,T]} |Z_t|^p\big]^{1/p}<\infty$.
  \item For $p\geq 1$, we let $\mcS^{p}_l$ be the set of all $\R$-valued, $\mcP_{\bbF}$-measurable, \caglad processes $(Z_t: t\in[0,T])$ such that $\|Z\|_{\mcS^p}<\infty$.
  \item We let $\mcH^{p}$ denote the set of all $\R^d$-valued $\mcP_{\bbF}$-measurable processes $(Z_t: t\in[0,T])$ such that $\|Z\|_{\mcH^p}:=\E\big[\big(\int_0^T |Z_t|^2 dt\big)^{p/2}\big]^{1/p}<\infty$.
  \item For any probability measure $\bbQ$, we let $\mcS^p_\bbQ$ and $\mcH^p_\bbQ$ be defined as $\mcS^p$ and $\mcH^p$, respectively, with the exception that the norm is defined with expectation taken with respect to $\bbQ$, \ie $\|Z\|_{\mcS^p_\bbQ}:=\E^\bbQ\big[\sup_{t\in[0,T]} |Z_t|^p\big]$ and $\|Z\|_{\mcH^p_\bbQ}:=\E^\bbQ\Big[\big(\int_0^T |Z_t|^2 dt\big)^{p/2}\Big]<\infty$.
  %\item We let $\mcK^p$ ($\mcK^{p,d}$) be the subset of $\mcS^p$ ($\mcS^{p,d}$) of non-negative processes.
  \item We let $\mcT$ be the set of all $\bbF$-stopping times and for each $\eta\in\mcT$ we let $\mcT_\eta$ be the corresponding subsets of stopping times $\tau$ such that $\tau\geq \eta$, $\Prob$-a.s.
  \item For each $\tau\in\mcT$, we let $\mcI(\tau)$ be the set of all $\mcF_\tau$-measurable random variables taking values in $U$, so that $\mcI(\tau)$ is the set of all admissible interventions at time $\tau$.
  \item We let $\mcU$ be the set of all $u=(\tau_1,\ldots,\tau_N;\beta_1,\ldots,\beta_N)$, where $N$ is the (random) number of interventions, $(\tau_j)_{j=1}^N$ is a non-decreasing sequence of $\bbF$-stopping times taking values in $[0,T]$ and $\beta_j\in\mcI(\tau_j)$.
  \item We let $\mcU^f$ denote the subset of $u\in\mcU$ for which $N$ is $\Prob$-a.s.~finite (\ie $\mcU^f:=\{u\in\mcU:\: \Prob\big[\{\omega\in\Omega :$ $ N>k, \:\forall k>0\}\big]=0\}$) and for all $k\geq 0$ we let $\mcU^k:=\{u\in\mcU:\:N\leq k,\,\Prob{\rm - a.s.}\}$.%$N:=\sup\{j:\tau_j\leq T\}$
  \item For $\eta\in\mcT$ we let $\mcU_{\eta}$ (and $\mcU_{\eta}^f$ resp.~$\mcU_{\eta}^k$) be the subset of $\mcU$ (and $\mcU^f$ resp.~$\mcU^k$) with $\tau_1\geq\eta$, $\Prob$-a.s.
  \item We let $\mcA$ be the set of all $\mcP_{\bbF}$-measurable processes $(\alpha_t)_{0\leq t\leq T}$ taking values in $A$ and for each $t\in[0,T]$ we let $\mcA_t$ be the set of all $\mcP_{\bbF}$-measurable processes $(\alpha_s)_{t\leq s\leq T}$ taking values in $A$.
  \item We denote by $\mcD$ the set of all double sequences $(t_1,\ldots;b_1,\ldots)$ where $(t_j)_{j\geq 1}$ is a non-decreasing sequence in $[0,T]$ and $b_j\in U$ for $j\geq 1$.
  \item We let $\mcD^f$ be the subset of $\mcD$ with all finite sequences and for $k\geq 0$ we let $\mcD^k$ be the subset of sequences with precisely $k$ interventions, \ie sequences of the type $(t_1,\ldots,t_k;b_1,\ldots,b_k)$.
  \item Throughout, we let $\vecv=(\vect,\vecb)$, with $\vect:=(t_1,\ldots,t_n)$ and $\vecb:=(b_1,\ldots,b_n)$, where $n$ is possibly infinite, denote a generic element of $\mcD$.
  \item For $\vecv=(\vect,\vecb)\in\mcD^f$ and $\vecv'=(\vect',\vecb')\in\mcD$ we introduce the concatenation, denoted by $\circ$, defined as $\vecv\circ\vecv':=(t_1,\ldots,t_n,t'_1\vee t_n,\ldots,t'_{n'}\vee t_n;b_1,\ldots,b_n,b'_1,\ldots,b'_{n'})$. Furthermore, for $\vecv\in\mcD$, we define the truncation to $k\geq 0$ interventions as $[\vecv]_{k}:=(t_1,\ldots,t_{k\wedge n};b_1,\ldots,b_{k\wedge n})$.
  \item For each $u\in\mcU^f$ we let $u(t)=[u]_{N(t)}$ with $N(t):=\max\{j\geq 0:\tau_j\leq t\}$.
  %\item For $l\geq 0$, we let $\Pi_l:=\{0,1/2^{l},2/2^{l},\ldots\}$ and set $\bold \Pi:=\cup_{l=1}^T \Pi_l$.
  \item We introduce the norm $\|\vecv\|_{\mcD^f}:=\sum_{j=1}^n(|t_j|+|b_j|)$ on $\mcD^f$ and let $\|\vecv-\vecv'\|_{\mcD^f}:=\infty$ whenever $n\neq n'$.
  %\item We extend the above norm to $\mcU^f$ by setting $\|u\|_{\mcU^f,C}:=\E\Big[\|u\|_{\mcD^f,C}\Big]$ for all $u\in\mcU^f$.
  \item We let $*$ denote stochastic integration and set $(X*W)_{t,s}=\int_t^s X_rdW_r$.
  \item We let $\mcE$ denote the Dol\'eans-Dade exponential and use the notation
  \begin{align*}
    \mcE(X*W)_{t,s}=e^{\int_t^s X_rdW_r-\frac{1}{2}\int_t^s |X_r|^2dr}.
  \end{align*}
  Also, we write $\mcE(X*W)_{t}:=\mcE(X*W)_{0,t}$.
  \item For any $\mcP_{\bbF}$-measurable process $\zeta$ such that $\E[\mcE(\zeta*W)_T]=1$, we define $\bbQ^{\zeta}$ to be the probability measure equivalent to $\Prob$, such that $d\bbQ^\zeta=\mcE(\zeta*W)_Td\Prob$.
  \item For any non-negative, $\mcP_{\bbF}$-measurable \cadlag process $L$ we let $\PrM^L$ denote the set of all probability measures $\bbQ$ on $(\Omega,\mcF)$ such that $d\bbQ=\mcE(\zeta*W)_Td\Prob$, for some $\mcP_{\bbF}$-measurable process $\zeta$, with $|\zeta_t|\leq L_t$ for all $t\in[0,T]$ (outside of a $\Prob$-null set).
\end{itemize}
In addition, we will throughout assume that, unless otherwise specified, all inequalities hold in the $\Prob$-a.s.~sense.

Furthermore, we define the following set:
\begin{defn}\label{defn:mcO}
We let $\mcO_\bbF$ be the set of all $\mcP_\bbF\otimes\mcB(U)$-measurable maps\footnote{Throughout, we generally suppress dependence on $\omega$ and refer to $h\in\mcO_\bbF$ as a map $(t,b)\to h(t,b)$.} $h:\Omega\times [0,T]\times U\to\R$ such that for each $\tau\in\mcT^f$ and $\beta\in\mcI(\tau)$ we have $h(\tau,\beta)\in L^p(\Omega,\mcF,\Prob)$ for all $p\geq 0$ and (outside of a $\Prob$-null set) the map $(t,b)\mapsto h(t,b)$ is jointly continuous.
%we have for all $(t,b)\in [0,T]\times U$:
%\begin{enumerate}[i)]
%  \item\label{mcO:lim} The limit $\lim_{(t',b')\to(t,b)}h(t',b')$ exists,
%  \item\label{mcO:cont} $\lim_{t'\to t}\sup_{b'\in U}|h(t',b')-h(t,b')|=0$, and
%  \item\label{mcO:usc} $\lim_{b'\to b}h(t,b')\leq h(t,b)$.
%\end{enumerate}
\end{defn}

\begin{defn}\label{defn:consistency}
We refer to a family of processes $((X^v_t)_{0\leq t\leq T}:v\in \mcU^f)$ as being \emph{consistent} if for each $u\in\mcU^f$, the map $h:[0,T]\times \Omega\times U\to\R$ given by $h(t,b)=X^{u\circ(t,b)}_t$ is $\mcP_\bbF\otimes\mcB(U)$-measurable and for each $\tau\in\mcT$ and each $\beta\in\mcI(\tau)$ we have $X^{u\circ(\tau,\beta)}_\tau=h(\tau,\beta)$, $\Prob$-a.s.
\end{defn}

%\subsection{Problem formulation}

One of the main objectives of the present work is to show that \eqref{ekv:seq-bsde-FH} admits a unique solution. We, therefore, need to define what we mean by a solution to \eqref{ekv:seq-bsde-FH}.

\begin{defn}\label{defn:rbsde-solution}
A solution to \eqref{ekv:seq-bsde-FH} is a family $(Y^v,Z^v,K^v)_{v\in\mcU^f}$, where
\begin{enumerate}[i)]
\item the family $(Y^v)_{v\in\mcU^f}$ is consistent and for each $v\in\mcU^f$, we have $Y^v\in \mcS^2$ with a norm that is uniformly bounded in $v$ (\ie $\sup_{u\in\mcU^f}\|Y^u\|_{\mcS^2})<\infty$) and $(t,b)\mapsto Y^{v\circ(t,b)}_t\in\mcO_\bbF$,
\item $Z^v\in\mcH^2$ for each $v\in\mcU^f$; and
\item $K^v\in\mcS^2$ is non-decreasing with $K^v_0=0$.
\end{enumerate}
\end{defn}

%%%%%%%%%%%%%%%%%%%%%%%%%%%%%%%%%%%%%%%%%%%%%%%%%%%%%%%%%%%%%%%%%%%%%%%%%%%%%%%%%%%%%%%%%%%%%%%%%%%%%%%%%%%%%%%%%%%%%%%%%%%%%%%%%%%%%%%%%%%%%%%%%%

\subsection{Prior results on RBSDEs}

Our approach will rely heavily on the available theory of reflected backward SDEs and we, therefore, recall the following important result:

\begin{thm}\label{thm:ElKaroui-rbsde}\emph{(El Karoui \etal \cite{ElKaroui1})}
  Assume that
  \begin{enumerate}[a)]
  \item $\xi\in L^2(\Omega,\mcF_T,\Prob)$.
  \item The barrier $S$ is real-valued, $\mcP_\bbF$-measurable and continuous with $S^+\in\mcS^2$ and $S_T\leq\xi$.
  \item $f:[0,T]\times \Omega\times \R\times\R^d\to\R$ is such that $f(\cdot,y,z)\in\mcH^2$ for all $(y,z)\in \R\times\R^{d}$ and for some $k_f>0$ and all $(y,y',z,z')\in\R^{2(1+d)}$ we have
  \begin{align*}
     |f(t,y',z')-f(t,y,z)|\leq k_f(|y'-y|+|z'-z|).
  \end{align*}
  \end{enumerate}
  Then, there exists a unique triple $(Y,Z,K):=(Y_t,Z_t,K_t)_{0\leq t\leq T}$ with $Y,K\in\mcS^2$ and $Z\in\mcH^2$, where $K$ is non-decreasing with $K_0=0$, such that
\begin{align}\label{ekv:ElKaroui-rbsde}
  \begin{cases}
    Y_t=\xi+\int_t^T f(s,Y_s,Z_s)ds-\int_t^T Z_s dW_s+K_T-K_t,\\
    Y_t\geq S_t,\, \forall t\in [0,T] \mbox{ and }\int_0^T \left(Y_t-S_t\right)dK_t=0.
  \end{cases}
\end{align}
Furthermore\footnote{Throughout, $C$ will denote a generic positive constant that may change value from line to line.},
\begin{align}\label{ekv:ElKaroui-bound}
\|Y\|_{\mcS^2}^2+\|Z\|_{\mcH^2}^2+\|K\|_{\mcS^2}^2&\leq C\E\Big[|\xi|^{2}+\int_0^T|f(s,0,0)|^2ds+\sup_{t\in[0,T]}|(S_t)^+|^2\Big].
\end{align}
In addition, $Y$ can be interpreted as the Snell envelope in the following way
\begin{equation*}
      Y_t=\esssup_{\tau\in\mcT_t}\E\bigg[\int_t^\tau f(s,Y_s,Z_s)ds+S_\tau\ett_{[\tau<T]}+\xi \ett_{[\tau=T]}\Big|\mcF_t\bigg]
\end{equation*}
and with $D_t:=\inf\{r\geq t: Y_r=S_r\}\wedge T$ we have the representation
\begin{equation*}
      Y_t=\E\bigg[\int_t^{D_t} f(s,Y_s,Z_s)ds+S_{D_t}\ett_{[D_t<T]}+\xi \ett_{[D_t=T]}\Big|\mcF_t\bigg]
\end{equation*}
and $K_{D_t}-K_t=0$, $\Prob$-a.s.

Moreover, if $(\tilde Y,\tilde Z,\tilde K)$ is the solution to the reflected BSDE with parameters $(\tilde\xi,\tilde f,\tilde S)$, then
\begin{align}\nonumber
\|\tilde Y-Y\|_{\mcS^2}^2+\|\tilde Z-Z\|_{\mcH^2}^2+\|\tilde K-K\|_{\mcS^2}^2&\leq C(\|\tilde S-S\|_{\mcS^{2}}\Psi_T^{1/2}+\E\Big[|\tilde \xi-\xi|^{2}
\\
&\quad+\int_0^T |\tilde f(s,Y_s,Z_s)-f(s,Y_s,Z_s)|^2ds\Big]),\label{ekv:ElKaroui-diff}
\end{align}
where
\begin{align*}
\Psi_T:=\E\Big[|\tilde \xi|^{2}+|\xi|^{2}+\int_0^T (|\tilde f(s,0,0)|^2+|f(s,0,0)|^2)ds+\sup_{t\in[0,T]}|(\tilde S_{t})^+ + (S_{t})^+|^{2}\Big].
\end{align*}
\end{thm}

Existence and uniqueness of solutions to \eqref{ekv:ElKaroui-rbsde} under the assumption that the Lipschitz coefficient on the $z$-variable is a stochastic process was recently shown in \cite{ElAsri2020}. In the next section we show that the moment and stability estimates in Theorem~\ref{thm:ElKaroui-rbsde} translates to the case with stochastic Lipschitz coefficients and exponents $p\neq 2$ as well. The latter turns out to play a crucial role when later on showing continuity of the map $(t,b)\mapsto Y^{v\circ(t,b)}$.

%%%%%%%%%%%%%%%%%%%%%%%%%%%%%%%%%%%%%%%%%%%%%%%%%%%%%%%%%%%%%%%%%%%%%%%%%%%%%%%%%%%%%%%%%%%%%%%%%%%%%%%%%%%%%%%%%%%%%%%%%%%%%%%%%%%%%%%%%%%%%%%%%%

\section{Reflected BSDEs with stochastic Lipschitz coefficient\label{sec:simple-rbsde}}
In this section we adapt some known results on BSDEs and in particular reflected BSDEs to our specific setting and give some prior estimates that will be useful later on. In particular, we will consider the non-reflected BSDE with parameters $(f,\xi)$,
\begin{align}\label{ekv:bsde-no-reflection}
  Y_t=\xi+\int_t^T f(s,Y_s,Z_s)ds-\int_t^T Z_sdW_s,\quad \forall\,t\in[0,T]
\end{align}
and the reflected BSDE with parameters $(f,\xi,S)$,
\begin{align}\label{ekv:bsde-reflection}
  \begin{cases}
    Y_t=\xi+\int_t^T f(s,Y_s,Z_s)ds-\int_t^T Z_s dW_s+K_T-K_t,\\
    Y_t\geq S_t,\, \forall t\in [0,T] \mbox{ and }\int_0^T \left(Y_t-S_t\right)dK_t=0.
  \end{cases}
\end{align}

Throughout this section we will assume that for some $q'> 1$ and all $p\geq 1$ the following holds:
\begin{ass}\label{ass:on-rbsde-FH}
\begin{enumerate}[(i)]
  \item There is a $\Prob$-a.s.~non-negative, $\mcP_\bbF$-measurable, continuous process $(L_t:t\in[0,T])$ (our stochastic Lipschitz coefficient) such that for all $\mcP_\bbF$-measurable processes $(\zeta_t:t\in[0,T])$ with $|\zeta_t|\leq L_t$ for all $t\in[0,T]$ (outside of a $\Prob$-null set) we have $\E[|\mcE(\zeta*W)_T|^{q'}]<\infty$ and $\E^{\bbQ^{\zeta}}[|\mcE(-\zeta*W^{\zeta})_T|^{q'}]<\infty$.
  \item\label{ass:xi-FH} The terminal value $\xi\in L^{p}(\Omega,\mcF_T,\Prob)$.
  \item\label{ass:on-f-FH} The driver $(t,\omega,y,z)\mapsto f(t,y,z): [0,T]\times\Omega\times\R\times \R^{d}\to\R$ is $\mcP_{\bbF}\otimes\mcB(\R)\otimes\mcB(\R^{d})$-measurable. Furthermore, we have
  \begin{enumerate}[a)]
  \item The bound
  \begin{align}
  \E\Big[\int_0^T|f(s,0,0)|^{p}ds\Big]<\infty.
  \end{align}
  \item There is a constant $k_f>0$ such that
  \begin{align}
     |f(t,y',z')-f(t,y,z)|\leq k_f|y'-y|+L_t|z'-z|,
  \end{align}
  for all $(t,y,y',z,z')\in[0,T]\times\R^{2(1+d)}$, $\Prob$-a.s.
  \item The barrier $S$ is real-valued, $\mcP_\bbF$-measurable and continuous with $S^+\in\mcS^p$ and $S_T\leq\xi$, $\Prob$-a.s.
  \end{enumerate}
\end{enumerate}
\end{ass}

\bigskip

Under this set of assumptions, letting $q>1$ be such that $\frac{1}{q}+\frac{1}{q'}=1$, we have the following:
\begin{thm}\label{thm:bsde-no-reflection}
The BSDE \eqref{ekv:bsde-no-reflection} has a unique solution $(Y,Z)\in \mcS^p\times \mcH^p$, with
\begin{align*}
\|Y\|_{\mcS^p}^p+\|Z\|_{\mcH^p}^p\leq C\E\Big[|\xi|^{q^2 p}+\big(\int_t^T |f(s,0,0)|ds\big)^{q^2 p}\Big]^{1/q^2}.
\end{align*}
Furthermore, if $(\tilde Y,\tilde Z)$ is a solution to \eqref{ekv:bsde-no-reflection} with parameters $(\tilde f,\tilde \xi)$, then
\begin{align*}
\|\tilde Y-Y\|_{\mcS^p}^p+\|\tilde Z-Z\|_{\mcH^p}^p\leq C\E\Big[|\tilde \xi-\xi|^{q^2 p}+\big(\int_t^T |\tilde f(s,Y_s,Z_s)-f(s,Y_s,Z_s)|ds\big)^{q^2 p}\Big]^{1/q^2}.
\end{align*}
\end{thm}

\noindent\emph{Proof.} This follows by repeating the arguments in the proofs of Lemma 3.1 and Lemma 3.2 in~\cite{Imkeller}.\qed\\

As a step towards obtaining moment estimates for solution to \eqref{ekv:bsde-reflection} we follow the convention in \cite{ElAsri2020} and let $(Y^{m,n},Z^{m,n},K^{m,n})\in\mcS^2\times\mcH^2\times\mcS^2$ be the unique solution to \eqref{ekv:bsde-reflection} with parameters $(\ett_{[L_s\leq m]}f^+-\ett_{[L_s\leq n]}f^-,\xi,S)$, \ie
\begin{align*}
\begin{cases}
    Y^{m,n}_t=\xi+\int_t^Tf^{m,n}(s,Y^{m,n}_s,Z^{m,n}_s)ds+\int_t^TZ^{m,n}_sdW_s+K^{m,n}_T-K^{m,n}_t\\
    Y^{m,n}_t\geq S_t,\, \forall t\in [0,T] \mbox{ and }\int_0^T \left(Y^{m,n}_t-S_t\right)dK^{m,n}_t=0.
  \end{cases}
\end{align*}
where $f^{m,n}:=\ett_{[L_s\leq m]}f^+-\ett_{[L_s\leq n]}f^-$ is a standard Lipschitz driver.

\begin{lem}\label{lem:Z-mn-is-bbQ-martingale}
For all $m,n\geq 0$ we have $Z^{m,n}_s\in\mcH^2_\bbQ$ for all $\bbQ\in\PrM^{L\wedge(m\vee n)}$.
\end{lem}

\noindent\emph{Proof.} For each $\bbQ\in\PrM^{L\wedge (m\vee n)}$, there is a $\mcP_\bbF$-measurable process, $\zeta$, such that $|\zeta_t|\leq L_t\wedge(m\vee n)$, for all $t\in [0,T]$ (outside of a $\Prob$-null set) and $d\bbQ :=\mcE(\zeta*W)_Td\Prob$. By the Girsanov theorem (see \eg Chapter 15 in \cite{CohenElliottBook}), it follows that under the probability measure $\bbQ$, the process $W^\zeta:=W-\int_0^\cdot\zeta_sds$ is a Brownian motion. Moreover, the triple $(Y^{m,n},Z^{m,n},K^{m,n})$ solves
\begin{align*}
  \begin{cases}
    Y_t=\xi+\int_t^T (f^{m,n}(s,Y_s,Z_s)-\zeta_sZ_s)ds-\int_t^T Z_s dW^\zeta_s+K_T-K_t,\\
    Y_t\geq S_t,\, \forall t\in [0,T] \mbox{ and }\int_0^T \left(Y_t-S_t\right)dK_t=0.
  \end{cases}
\end{align*}
Since this is a reflected BSDE with standard Lipschitz driver we can apply \eqref{ekv:ElKaroui-bound} to get that
\begin{align*}
\E^\bbQ\Big[\int_0^T|Z^{m,n}_s|^2ds\Big]&\leq C\E^\bbQ\Big[|\xi|^{2}+\int_0^T|f^{m,n}(s,0,0)|^2ds+\sup_{t\in[0,T]}|(S_t)^+|^2\Big]
\\
&\leq C\E\Big[\mcE(\zeta*W)_T\Big(|\xi|^{2}+\int_0^T|f^{m,n}(s,0,0)|^2ds+\sup_{t\in[0,T]}|(S_t)^+|^2\Big)\Big]
\\
&\leq C\E\big[|\mcE(\zeta*W)_T|^{q'}]^{1/q'}\E\Big[|\xi|^{2q}+\int_0^T|f^{m,n}(s,0,0)|^{2q}ds+\sup_{t\in[0,T]}|(S_t)^+|^{2q}\Big]^{1/q}
\end{align*}
and the result follows by Assumption~\ref{ass:on-rbsde-FH}.\qed\\

\begin{lem}\label{lem:rbsde-trunk-bound}
There is a $C>0$ such that
\begin{align}\label{ekv:unif-bound-FH}
\|Y^{m,n}\|_{\mcS^p}^p + \|Z^{m,n}\|_{\mcH^p}^p + \|K^{m,n}\|_{\mcS^p}^p \leq C\E\Big[|\xi|^{q^2p}+\big(\int_0^T |f^{m,n}(s,0,0)|ds\big)^{q^2p}+\sup_{t\in[0,T]}|(S_{t})^+|^{q^2p}\Big]^{1/q^2}.
\end{align}
for all $m,n\geq 0$.
\end{lem}

\noindent\emph{Proof.} We define
\begin{align*}
\gamma_s:=\frac{f^{m,n}(s,Y^{m,n}_s,Z^{m,n}_s)-f^{m,n}(s,0,Z^{m,n}_s)}{Y^{m,n}_s}\ett_{[Y^{m,n}_s\neq 0]}
\end{align*}
and set $e_t:=e^{\int_0^t\gamma_sds}$. By Assumption~\ref{ass:on-rbsde-FH} we have that $e^{-k_fT}\leq e_t\leq e^{k_fT}$. Applying Ito's formula to $e_tY^{m,n}_t$ gives that for any $\tau\in\mcT_t$ we have
\begin{align*}
e_tY^{m,n}_t&=e_\tau Y^{m,n}_\tau +\int_t^\tau e_s(f^{m,n}(s,Y^{m,n}_s,Z^{m,n}_s)-\gamma_sY^{m,n}_s)ds-\int_t^\tau e_sZ^{m,n}_sdW_s+\int_t^\tau e_sdK^{m,n}_s
\\
&=e_\tau Y^{m,n}_\tau +\int_t^\tau e_s(f^{m,n}(s,0,0)+\zeta_sZ^{m,n}_s)ds-\int_t^\tau e_sZ^{m,n}_sdW_s+\int_t^\tau e_sdK^{m,n}_s
\end{align*}
where
\begin{align*}
\zeta_s:=\frac{f^{m,n}(s,0,Z^{m,n}_s)-f^{m,n}(s,0,0)}{|Z^{m,n}_s|^2}(Z^{m,n}_s)^\top\ett_{[Z^{m,n}_s\neq 0]}.
\end{align*}
By the Girsanov theorem (see \eg Chapter 15 in \cite{CohenElliottBook}), it follows that under the probability measure $\bbQ^{\zeta}$ defined as $d\bbQ^{\zeta}:=\mcE(\zeta*W)_Td\Prob$ the process $W^{\zeta}_t:=W_t-\int_0^t\zeta_sds$ is a Brownian motion. Furthermore, we have
\begin{align*}
e_t Y^{m,n}_t=e_\tau Y^{m,n}_\tau +\int_t^\tau e_sf^{m,n}(s,0,0)ds-\int_t^\tau e_sZ^{m,n}_sdW^{\zeta}_s+\int_t^\tau e_sdK^{m,n}_s.
\end{align*}
Taking the conditional expectation while picking $\tau=D_t(:=\inf\{t\geq 0:Y^{m,n}_t=S_t\}\wedge T)$ and appealing to Lemma~\ref{lem:Z-mn-is-bbQ-martingale} which implies that the stochastic integral is a $\bbQ^{\zeta}$-martingale gives
\begin{align}\label{ekv:Y-mn-alt-desc}
e_tY^{m,n}_t=\E^{\bbQ^{\zeta}}\Big[\ett_{[D_t<T]}e_{D_t}S_{D_t}+\ett_{[D_t=T]}e_{T}\xi+\int_t^{D_t} e_sf^{m,n}(s,0,0)ds\big|\mcF_t\Big].
\end{align}
Doob's maximal inequality together with % the fact that $Y^{m,n}$ is bounded from below by a non-reflected BSDE and
the bounds on $e_t$ gives that
\begin{align*}
\E^{\bbQ^{\zeta}}\Big[\sup_{t\in[0,T]}|Y^{m,n}_t|^p\Big]\leq C\E^{\bbQ^{\zeta}}\Big[\big(\int_0^{T}|f^{m,n}(s,0,0)| ds\big)^p+\sup_{t\in[0,T]}|(S_{t})^+|^p+|\xi|^p\Big].
\end{align*}
Changing back to our original probability measure $\Prob$ on the right hand side gives
\begin{align*}
\E^{\bbQ^{\zeta}}\Big[\sup_{t\in[0,T]}|Y^{m,n}_t|^p\Big]&\leq C\E\Big[\mcE(\zeta*W)_T\Big(\big(\int_0^{T}|f^{m,n}(s,0,0)| ds\big)^p+\sup_{t\in[0,T]}|(S_{t})^+|^p+|\xi|^p\Big)\Big]
\\
&\leq C\E\big[\mcE(\zeta*W)_T^{q'}]^{1/q'}\E\Big[\big(\int_0^{T}|f^{m,n}(s,0,0)|ds\big)^{qp}+\sup_{t\in[0,T]}|(S_{t})^+|^{qp}+|\xi|^{qp}\Big]^{1/q}.
\end{align*}
Moreover, as
\begin{align*}
\E\Big[\sup_{t\in[0,T]}|Y^{m,n}_t|^p\Big]&=\E^{\bbQ^{\zeta}}\Big[\mcE(-\zeta*W^{\zeta})_T\sup_{t\in[0,T]}|Y^{m,n}_t|^p\Big]
\\
&\leq C\E^{\bbQ^{\zeta}}\Big[\sup_{t\in[0,T]}|Y^{m,n}_t|^{qp}\Big]^{1/q}
\end{align*}
the estimate for $Y^{m,n}$ follows.\\

To arrive at the estimate for $Z^{m,n}$ we apply Ito's formula to $(Y^{m,n})^2$ and get that
\begin{align}\nonumber
|Y^{m,n}_t|^2+\int_t^T|Z^{m,n}_s|^2ds &= |\xi|^2+2\int_t^TY^{m,n}_sf^{m,n}(s,Y^{m,n}_s,Z^{m,n}_s)ds
\\
&\quad-2\int_t^TY^{m,n}_sZ^{m,n}_sdW_s + 2\int_t^TY^{m,n}_sdK^{m,n}_s\nonumber
\\
&= |\xi|^2+2\int_t^T(\gamma_s|Y^{m,n}_s|^2+Y^{m,n}_sf^{m,n}(s,0,0))ds\nonumber
\\
&\quad-2\int_t^TY^{m,n}_sZ^{m,n}_sdW^{\zeta}_s + 2\int_t^TY^{m,n}_sdK^{m,n}_s.\label{ekv:ito-on-y2}
\end{align}
Furthermore, applying the relation $ab\leq \tfrac{1}{2}\kappa a^2+\tfrac{1}{2\kappa}b^2$ with $\kappa>0$ and the Skorokhod condition $\int_0^T(Y^{m,n}_s-S_s)dK^{m,n}_s=0$ yields
\begin{align*}
\int_t^TY^{m,n}_sdK^{m,n}_s\leq\int_t^T(S_s)^+dK^{m,n}_s\leq \frac{\kappa}{2}\sup_{r\in[t,T]}|(S_r)^+|^2+\frac{1}{2\kappa}|K^{m,n}_T|^2.
\end{align*}

Now, Proposition~2.2 in \cite{ElKaroui1} gives that
\begin{align*}
K^{m,n}_T-K^{m,n}_t&=\sup_{r\in[t,T]}\Big(\xi+\int_r^Tf^{m,n}(s,Y^{m,n}_s,Z^{m,n}_s)ds-\int_r^TZ^{m,n}_sdW_s-S_r\Big)^-
\\
&=\sup_{r\in[t,T]}\Big(\xi+\int_r^T(\gamma_sY^{m,n}_s+f^{m,n}(s,0,0))ds-\int_r^TZ^{m,n}_sdW^{\zeta}_s-S_r\Big)^-
\end{align*}
and, in particular, we have that
\begin{align}\label{ekv:K_T-bound}
K^{m,n}_T&\leq |\xi|+\int_0^T(k_f|Y^{m,n}_s|+|f^{m,n}(s,0,0)|)ds+2\sup_{r\in[0,T]}\Big|\int_0^r Z^{m,n}_sdW^{\zeta}_s\Big|+\sup_{r\in[0,T]}|(S_r)^+|
\end{align}
Inserted into \eqref{ekv:ito-on-y2} this gives
\begin{align*}
\int_0^T|Z^{m,n}_s|^2ds &\leq \frac{C}{\kappa}\sup_{r\in[0,T]}\Big|\int_0^TZ^{m,n}_sdW^{\zeta}_s\Big|^2 -2\int_0^TY^{m,n}_sZ^{m,n}_sdW^{\zeta}_s
\\
&\quad +C(1+\kappa)\big(|\xi|^2+\sup_{r\in[0,T]}|Y^{m,n}_r|^2+\sup_{r\in[0,T]}|(S_r)^+|^2+\big(\int_0^T|f^{m,n}(s,0,0)|ds\big)^2\big)
\end{align*}
or
\begin{align*}
&\E^{\bbQ^{\zeta}}\Big[\big(\int_0^T|Z^{m,n}_s|^2ds)^{p/2}-\frac{C}{\kappa^{p/2}}\sup_{r\in[0,T]}\Big|\int_0^TZ^{m,n}_sdW^{\zeta}_s\Big|^p -C\Big|\int_0^TY^{m,n}_sZ^{m,n}_sdW^{\zeta}_s\Big|^{p/2}\Big]
\\
& \leq C(1+\kappa^{p/2})\E^{\bbQ^{\zeta}}\Big[|\xi|^p+\sup_{r\in[0,T]}|Y^{m,n}_r|^p+\sup_{r\in[0,T]}|(S_r)^+|^p+\big(\int_0^T|f^{m,n}(s,0,0)|ds\big)^p\Big].
\end{align*}
On the other hand, applying the Burkholder-Davis-Gundy inequality (BDG for short) gives that
\begin{align*}
&\E^{\bbQ^{\zeta}}\Big[\big(\int_0^T|Z^{m,n}_s|^2ds)^{p/2}-\frac{C}{\kappa^{p/2}}\sup_{r\in[0,T]}\Big|\int_0^TZ^{m,n}_sdW^{\zeta}_s\Big|^p -C\Big|\int_0^TY^{m,n}_sZ^{m,n}_sdW^{\zeta}_s\Big|^{p/2}\Big]
\\
&\geq  (1-\frac{C}{\kappa^{p/2}})\E^{\bbQ^{\zeta}}\Big[\big(\int_0^T|Z^{m,n}_s|^2ds\big)^{p/2}\Big] - C\E^{\bbQ^{\zeta}}\Big[\big(\int_0^T|Y^{m,n}_sZ^{m,n}_s|^2ds\big)^{p/4}\Big]
\end{align*}
By again using that $ab\leq \tfrac{1}{2}\kappa a^2+\tfrac{1}{2\kappa}b^2$ we get
\begin{align*}
\E^{\bbQ^{\zeta}}\Big[\big(\int_0^T|Y^{m,n}_sZ^{m,n}_s|^2ds\big)^{p/4}\Big]\leq C\E^{\bbQ^{\zeta}}\Big[\kappa\sup_{s\in [0,T]}|Y^{m,n}_s|^p+\frac{1}{\kappa}\big(\int_0^T|Z^{m,n}_s|^2ds\big)^{p/2}\Big]
\end{align*}
and the estimate follows by choosing $\kappa>0$ sufficiently large and repeating the steps above to change back to the original measure $\Prob$. Finally, the bound for $K^{m,n}$ is immediate from \eqref{ekv:K_T-bound} the above.\qed\\

\begin{lem}\label{lem:rbsde-trunk-diff}
If $(Y^{m,n},Z^{m,n},K^{m,n})$ and $(\tilde Y^{\tilde m,\tilde n},\tilde Z^{\tilde m,\tilde n},\tilde K^{\tilde m,\tilde n})$ are solutions to \eqref{ekv:bsde-reflection} with parameters $(f^{m,n},\xi,S)$ and $(\tilde f^{\tilde m,\tilde n},\tilde\xi,\tilde S)$, respectively, then
\begin{align*}
&\|\tilde Y^{\tilde m,\tilde n}-Y^{m,n}\|_{\mcS^p}^p+\|\tilde Z^{\tilde m,\tilde n}-Z^{m,n}\|_{\mcH^p}^p + \|\tilde K^{\tilde m,\tilde n}-K^{m,n}\|_{\mcS^p}^p
\\
&\leq C(\|\tilde S-S\|_{\mcS^{q^2p}}^{p/2}\Psi_T^{1/2}+\E\Big[|\tilde \xi-\xi|^{q^2p}+\big(\int_0^T |\tilde f^{\tilde m,\tilde n}(s,Y_s,Z_s)-f^{m,n}(s,Y_s,Z_s)|ds\big)^{q^2p}\Big]^{1/q^2}),
\end{align*}
where
\begin{align*}
\Psi_T:=\E\Big[|\tilde \xi|^{q^2p}+|\xi|^{q^2p}+\big(\int_0^T (|\tilde f^{\tilde m,\tilde n}(s,0,0)|+|f^{m,n}(s,0,0)|)ds\big)^{q^2p} +\sup_{t\in[0,T]}|(\tilde S_{t})^+ + (S_{t})^+|^{q^2p}\Big]^{1/q^2}.
\end{align*}
\end{lem}

\noindent\emph{Proof.} For ease of notation we omit the superscripts $m,n$ and $\tilde m,\tilde n$ and have
\begin{align*}
\tilde Y_{t}-Y_{t}= \tilde \xi-\xi+\int_t^T (\tilde f(s,\tilde Y_s,\tilde Z_s)-f(s,Y_s,Z_s))ds-\int_t^T (\tilde Z_s-Z_s)dW_s+\tilde K_T-\tilde K_t-(K_T-K_t)
\end{align*}
Now, let
\begin{align*}
\gamma_s:=\frac{\tilde f(s,\tilde Y_s,\tilde Z_s)-\tilde f(s,Y_s,\tilde Z_s)}{\tilde Y_s-Y_s}\ett_{[\tilde Y_s\neq Y_s]}
\end{align*}
and
\begin{align*}
\zeta_s:=\frac{\tilde f(s,Y_s,\tilde Z_s)-\tilde f(s,Y_s,Z_s)}{|\tilde Z_s-Z_s|^2}(\tilde Z_s-Z_s)^\top\ett_{[\tilde Z_s\neq Z_s]}
\end{align*}
and set $e_t:=e^{\int_0^t\gamma_sds}$. Letting $\delta Y:=\tilde Y-Y$, $\delta\xi:=\tilde\xi -\xi$, $\delta Z:=\tilde Z - Z$ and $\delta f:=\tilde f-f$ we get that for any $\tau\in\mcT_t$ we have
\begin{align*}
e_t\delta Y_t&= e_\tau\delta Y_\tau+\int_t^\tau e_s(\delta f(s,Y_s,Z_s)+\zeta_s\delta Z_s)ds-\int_t^\tau e_s\delta Z_sdW_s+\int_t^\tau e_sd(\delta K)_s
\\
&= e_\tau\delta Y_\tau+\int_t^\tau e_s\delta f(s,Y_s,Z_s)ds-\int_t^\tau e_s\delta Z_sdW^{\zeta}_s+\int_t^\tau e_sd(\delta K)_s
\end{align*}
where $W^{\zeta}_t=W_t-\int_0^t\zeta_sds$. Now set $\tau=D_t\,(=\inf\{r\geq t:Y_r=S_r\}\wedge T)$ and we find that
\begin{align*}
e_t\delta Y_t&=  e_\tau\delta Y_{D_t}+\int_t^{D_t} e_s\delta f(s,Y_s,Z_s)ds-\int_t^{D_t} e_s\delta Z_sdW^{\zeta}_s+\int_t^{D_t}e_s(d\tilde K_s-dK_s)
\\
&\geq \ett_{[D_t<T]}e_{D_t}\delta S_{D_t}+\ett_{[D_t=T]}e_T\delta\xi+\int_t^{D_t} e_s \delta f(s,Y_s,Z_s)ds-\int_t^{D_t} e_s\delta Z_sdW^{\zeta}_s.
\end{align*}
On the other hand, by picking $\tau=\tilde D_t:=\inf\{r\geq t:\tilde Y_r=\tilde S_r\}\wedge T$ we get
\begin{align*}
e_t\delta Y_t&\leq \ett_{[\tilde D_t<T]}e_{\tilde D_t}\delta S_{\tilde D_t}+e_T\ett_{[\tilde D_t=T]}\delta\xi+\int_t^{\tilde D_t} e_s\delta f(s,Y_s,Z_s)ds-\int_t^{\tilde D_t} e_s\delta Z_sdW^{\zeta}_s.
\end{align*}
Taking the conditional expectation with respect to the measure $\bbQ^{\zeta}$, with $d\bbQ^{\zeta}=\mcE(\zeta*W)_Td\Prob$, we find that
\begin{align*}
e_t|\delta Y_t|&\leq \E^{\bbQ^{\zeta}}\Big[\sup_{r\in[0,T]}e_r|\delta S_{r}|+e_T|\delta\xi|+\int_0^{T} e_s|\delta f(s,Y_s,Z_s)|ds\big|\mcF_t\Big]
\end{align*}
and the bound for $\|\tilde Y^{\tilde m,\tilde n}-Y^{m,n}\|_{\mcS^p}$ follows by appealing to the argument of Lemma~\ref{lem:rbsde-trunk-bound} while noting that $|\gamma_s|\leq k_f$ and $|\zeta_s|\leq L_s$.\\

Applying Ito's formula to $\delta Y^2$ gives (similarly to \eqref{ekv:ito-on-y2} that)
\begin{align}\nonumber
|\delta Y_t|^2+\int_t^T|\delta Z_s|^2ds &= |\xi|^2+2\int_t^T(\gamma_s|\delta Y_s|^2+\delta Y_s\delta f(s,Y_s,Z_s))ds\nonumber
\\
&\quad-2\int_t^T\delta Y_s\delta Z_sdW^{\zeta}_s + 2\int_t^T\delta Y_sd(\delta K)_s\nonumber
\\
&\leq |\xi|^2+2\int_t^T(k_f|\delta Y_s|^2+\delta Y_s\delta f(s,Y_s,Z_s))ds\nonumber
\\
&\quad-2\int_t^T\delta Y_s\delta Z_sdW^{\zeta}_s + 2\int_t^T\delta S_sd(\delta K)_s,\label{ekv:ito-on-delta_y2}
\end{align}
where the last inequality follows from the fact that $\int_0^T\delta Y_sd(\delta K_s)\leq \int_0^T\delta S_sd(\delta K_s)$. Now,
\begin{align*}
\int_0^T\delta S_s d(\delta K)_s\leq \sup_{r\in[0,T]}|\delta S_r|(\tilde K_T+ K_T)
\end{align*}
which gives that
\begin{align*}
\E^{\bbQ^{\zeta}}\Big[\big(\int_0^T|\delta Z_s|^2ds\big)^{p/2}\Big] &\leq C\E^{\bbQ^{\zeta}}\Big[ |\xi|^p+\sup_{r\in[0,T]}|\delta Y_s|^p+\big(\int_0^T|\delta f(s,Y_s,Z_s)| ds\big)^p
\\
&\quad +\Big|\int_0^T\delta Y_s\delta Z^{m,n}_sdW^{\zeta}_s\Big|^{p/2} + \sup_{r\in[0,T]}|\delta S_r|^{p/2}((\tilde K_T)^{p/2}+ (K_T)^{p/2})\Big]
\end{align*}
and repeating the last steps in the proof of Lemma~\ref{lem:rbsde-trunk-bound} we have that
\begin{align*}
\E\Big[\big(\int_0^T|\delta Z_s|^2ds\big)^{p/2}\Big] &\leq C\E\Big[ |\xi|^{q^2p}+\sup_{r\in[0,T]}|\delta Y_s|^{q^2p}+\big(\int_0^T|\delta f(s,Y_s,Z_s)|ds\big)^{q^2p}
\\
&\quad \sup_{r\in[0,T]}|\delta S_r|^{q^2p/2}((\tilde K_T)^{q^2p/2}+ (K_T)^{q^2p/2})\Big]^{1/q^2}.
\end{align*}
Finally,
\begin{align*}
\E\Big[\sup_{r\in[0,T]}|\delta S_r|^{q^2p/2}((\tilde K_T)^{q^2p/2}+ (K_T)^{q^2p/2})\Big]\leq 2\E\Big[\sup_{r\in[0,T]}|\delta S_r|^{q^2p}\Big]^{1/2}\E\Big[(\tilde K_T)^{q^2p}+ (K_T)^{q^2p}\Big]^{1/2}
\end{align*}
%and the bound for $\|\tilde Z^{\tilde m,\tilde n}-Z^{m,n}\|_{\mcH^p}$ follows by noting that writing
%while noting that
%\begin{align*}
%|\delta K_T|&\leq |\delta \xi|+\int_r^T(k_f|\delta Y_s|+|\delta f(s,Y_s,Z_s)|)ds+\sup_{r\in[0,T]}\Big|\int_0^r \delta Z_sdW^{\zeta}_s\Big|+\sup_{r\in[0,T]}|\delta S_r|
%\end{align*}
and the result follows by applying Lemma~\ref{lem:rbsde-trunk-bound} to the last term.\qed\\

\begin{prop}\label{prop:rbsde-solu-FH}
The reflected BSDE \eqref{ekv:bsde-reflection} has a unique solution $(Y,Z,K)\in \mcS^2\times \mcH^2\times \mcS^2$, with $K$ non-decreasing and $K_0=0$. Furthermore, we have
\begin{align}\label{ekv:YZK-bound-FH}
\|Y\|_{\mcS^p}^p+\|Z\|_{\mcH^p}^p+\|K\|_{\mcS^p}^p\leq C(\|S^+\|_{\mcS^{q^2p}}^p+\E\Big[|\xi|^{q^2p}+\big(\int_0^T |f(s,0,0)|ds\big)^{q^2p}\Big]^{1/q^2})
\end{align}
and if $(\tilde Y,\tilde Z,\tilde K)$ is a solution to \eqref{ekv:bsde-reflection} with parameters $(\tilde f,\tilde \xi,\tilde S)$ then
\begin{align}\nonumber
&\|\tilde Y-Y\|_{\mcS^p}^p+\|\tilde Z-Z\|_{\mcH^p}^p+\|\tilde K-K\|_{\mcS^p}^p
\\
&\leq C(\|\tilde S-S\|_{\mcS^{q^2p}}^{p/2}\Psi_T^{1/2}+\E\Big[|\tilde \xi-\xi|^{q^2p}+\big(\int_0^T |\tilde f(s,Y_s,Z_s)-f(s,Y_s,Z_s)|ds\big)^{q^2p}\Big]^{1/q^2}),\label{ekv:YZK-diff-FH}
\end{align}
where
\begin{align*}
\Psi_T:=\E\Big[|\tilde \xi|^{q^2p}+|\xi|^{q^2p}+\big(\int_0^T (|\tilde f(s,0,0)|+|f(s,0,0)|)ds\big)^{q^2p}+\sup_{t\in[0,T]}|(\tilde S_{t})^+ + (S_{t})^+|^{q^2p}\Big]^{1/q^2}.
\end{align*}
\end{prop}

\noindent\emph{Proof.} The first part of the proof largely follows that of Lemma 2.3 in \cite{ElAsri2020} and is included for the sake of completeness. Taking the limit on the right-hand side of \eqref{ekv:unif-bound-FH} and appealing to dominated convergence we find that
\begin{align*}
\|Y^{m,n}\|_{\mcS^p}^p+\|Z^{m,n}\|_{\mcH^p}^p+\|K^{m,n}\|_{\mcS^p}^p&\leq C(\|S^+\|_{\mcS^{q^2p}}^p+\E\Big[|\xi|^{q^2p}+\big(\int_0^T |f(s,0,0)|ds\big)^{q^2p}\Big]^{1/q^2})
\\
&\leq C,
\end{align*}
by Assumption~\ref{ass:on-rbsde-FH} for all $m,n\geq 0$ and $p\geq 1$. In particular, since for fixed $m$ comparison implies that the sequence $(Y^{m,n})_{n\geq 0}$ of continuous processes is non-increasing, this implies that $Y^{m}:=\lim_{n\to\infty}Y^{m,n}$ exists, $\Prob$-a.s., and by Fatou's lemma it satisfies $\|Y^{m}\|_{\mcS^p}^p\leq C$.  Furthermore, we have by Lemma~\ref{lem:rbsde-trunk-diff} that
\begin{align*}
\|Y^{m,n}-Y^{m,n'}\|_{\mcS^2}\to 0,
\end{align*}
as $n,n'\to\infty$ implying that $Y^{m}$ is a continuous process. Moreover, as $(Z^{m,n})_{m,n\geq 0}$ is a uniformly bounded double-sequence in $\mcH^2$ and by again appealing to Lemma~\ref{lem:rbsde-trunk-diff} we have that
\begin{align*}
\|Z^{m,n}-Z^{m,n'}\|_{\mcH^2}\to 0,
\end{align*}
as $n,n'\to\infty$ and we conclude that there is a $Z^m$ and consequently also a $K^m$ such that $Z^{m,n}\to Z^m$ in $\mcH^2$ and $K^{m,n}\to K^m$, in $\mcS^2$ as $n\to\infty$. Now, let $\eta_l:=\inf\{s>0:L_s\geq l\}\wedge T$ and note that by Theorem~\ref{thm:ElKaroui-rbsde} there is a unique triple $(\hat Y^m_t,\hat Z^m_t,\hat K^m_t)_{0\leq t\leq \eta_l}$ such that
\begin{align*}
  \begin{cases}
  \hat Y^m_t=Y^m_{\eta_l}+\int_t^{\eta_l} f^m(s,\hat Y^m_s,\hat Z^m_s)ds-\int_t^{\eta_l} \hat Z^m_sdW_s+\hat K^m_{\eta_l}-\hat K^m_t,\quad \forall t\in [0, {\eta_l}],\\
  \hat Y^m_t\geq S^m_t,\, \forall t\in [0, {\eta_l}] \mbox{ and }\int_0^{\eta_l}(\hat Y^m_t-S^m_t)d\hat K^m_t=0,
  \end{cases}
\end{align*}
with $f^{m}:=\ett_{[L_s\leq m]}f^+-f^-$. A trivial modification of the stability result in \eqref{ekv:ElKaroui-diff} to random terminal times gives that
\begin{align*}
&\|(Y^{m,n}-\hat Y^{m})\ett_{[0,\eta_l]}\|_{\mcS^2}^2+\|(Z^{m,n}-\hat Z^{m})\ett_{[0,\eta_l]}\|_{\mcH^2}^2 + \|(K^{m,n}-\hat K^{m})\ett_{[0,\eta_l]}\|_{\mcS^2}^2 \to 0,
\end{align*}
as $n\to\infty$ and by uniqueness of limits it follows that $(Y^m,Z^m,K^m)=(\hat Y^m,\hat Z^m,\hat K^m)$ in $\mcS^2\times\mcH^2\times\mcS^2$. Since $\eta_l=T$ outside of a $\Prob$-null set for $l\,(=l(\omega))$ sufficiently large, we conclude that $(Y^m,Z^m,K^m)$ is the unique solution to
\begin{align*}
  \begin{cases}
  Y^m_t=\xi+\int_t^{T} f^m(s,Y^m_s,Z^m_s)ds-\int_t^{T} Z^m_sdW_s+K^m_{T}-K^m_t,\quad \forall t\in [0, T],\\
  Y^m_t\geq S^m_t,\, \forall t\in [0, T] \mbox{ and }
  \int_0^{T}(Y^m_t-S^m_t)dK^m_t=0.
  \end{cases}
\end{align*}
From the above, the bounds \eqref{ekv:YZK-bound-FH} and \eqref{ekv:YZK-diff-FH} with $f=f^m$ are obtained through Lemma~\ref{lem:rbsde-trunk-bound} and Lemma~\ref{lem:rbsde-trunk-diff}, respectively, by Fatou's lemma and dominated convergence.\\

Letting $m\to\infty$ and repeating the above argument the result follows.\qed\\

\begin{cor}\label{cor:rbsde-char-FH}
If $(Y,Z,K)$ solves \eqref{ekv:bsde-reflection}, then $Y$ can be interpreted as the Snell envelope in the following way
\begin{equation}\label{ekv:Y-is-Snell-FH}
      Y_t=\esssup_{\tau\in\mcT_t}\E\bigg[\int_t^\tau f(s,Y_s,Z_s)ds+S_\tau\ett_{[\tau<T]}+\xi \ett_{[\tau=T]}\Big|\mcF_t\bigg].
\end{equation}
In particular, with $D_t:=\inf\{r\geq t: Y_r=S_r\}\wedge T$ we have the representation
\begin{equation}\label{ekv:Snell-attained}
      Y_t=\E\bigg[\int_t^{D_t} f(s,Y_s,Z_s)ds+S_{D_t}\ett_{[D_t<T]}+\xi \ett_{[D_t=T]}\Big|\mcF_t\bigg]
\end{equation}
and $K_{D_t}-K_t=0$, $\Prob$-a.s.
\end{cor}

\noindent\emph{Proof.} Let $\eta_l:=\inf\{s>0:L_s\geq l\}\wedge T$, then Theorem~\ref{thm:ElKaroui-rbsde} gives that for all $l\geq 0$ we have
\begin{align*}
      Y_t=\E\bigg[\int_t^{D_t\wedge\eta_l} f(s,Y_s,Z_s)ds+S_\tau\ett_{[D_t<\eta_l]}+\ett_{[D_t\geq \eta_l]}(\ett_{[\eta_l<T]}Y_{\eta_l}+\ett_{[\eta_l=T]}\xi)\Big|\mcF_t\bigg]
\end{align*}
and \eqref{ekv:Snell-attained} follows by letting $l\to\infty$ and using dominated convergence. Optimality of $D_t$ follows similarly, implying that \eqref{ekv:Y-is-Snell-FH} holds.\qed\\

%Let $\hat Y$ denote the right-hand side of \eqref{ekv:Y-is-Snell-FH}. Then $\hat Y$ is the Snell envelope of a \cadlag process with positive jumps and we have the characterisation
%\begin{equation*}
%      \hat Y_t=\E\bigg[\int_t^{\hat D_t} f(s,Y_s,Z_s)ds+S_{\hat D_t}\ett_{[\hat D_t<T]}+\xi \ett_{[\hat D_t=T]}\Big|\mcF_t\bigg],
%\end{equation*}
%where $\hat D_t:=\inf\{s\geq t: \hat Y_s=S_s\}\wedge T$. and uniqueness of the Snell envelope that $Y=\hat Y$.\qed\\

%%%%%%%%%%%%%%%%%%%%%%%%%%%%%%%%%%%%%%%%%%%%%%%%%%%%%%%%%%%%%%%%%%%%%%%%%%%%%%%%%%%%%%%%%%%%%%%%%%%%%%%%%%%%%%%%%%%%%%%%%%%%%%%%%%%%%%%%%%%%%%%%%%
%%%%%%%%%%%%%%%%%%%%%%%%%%%%%%%%%%%%%%%%%%%%%%%%%%%%%%%%%%%%%%%%%%%%%%%%%%%%%%%%%%%%%%%%%%%%%%%%%%%%%%%%%%%%%%%%%%%%%%%%%%%%%%%%%%%%%%%%%%%%%%%%%%
%%%%%%%%%%%%%%%%%%%%%%%%%%%%%%%%%%%%%%%%%%%%%%%%%%%%%%%%%%%%%%%%%%%%%%%%%%%%%%%%%%%%%%%%%%%%%%%%%%%%%%%%%%%%%%%%%%%%%%%%%%%%%%%%%%%%%%%%%%%%%%%%%%

\section{Sequential systems of reflected BSDEs in finite horizon\label{sec:seq-rbsde-FH}}
In this section we move on to the sequential system of reflected BSDEs in \eqref{ekv:seq-bsde-FH}. To be able to use our results for impulse control we must allow the stochastic Lipschitz coefficient in \eqref{ekv:seq-bsde-FH} to depend on the control parameter $v\in\mcU^f$, rendering us a family of stochastic Lipschitz coefficients $(L^v:v\in\mcU^f)$ with $L^v\in\mcS^p$ for all $v\in\mcU^f$ and $p\geq 0$.

\subsection{Motivation}
The motivation for our study of the system of reflected BSDEs in \eqref{ekv:seq-bsde-FH} is applications to non-Markovian, robust impulse control. In this context \eqref{ekv:seq-bsde-FH} takes the form of a system of reflected forward-backward SDEs of the form
\begin{align}\label{ekv:seq-fbsde}
\begin{cases}
  Y^{v}_t=\psi(X^v_T)+\int_t^T h(s,X^{v}_s,Z^{v}_s)ds-\int_t^T Z^{v}_sdW_s+ K^{v}_T-K^{v}_t,\quad\forall t\in[0,T], \\
  Y^{v}_t\geq\sup_{b\in U}\{Y^{v\circ(t,b)}_t-\ell(t,X^{v}_t,b))\},\quad\forall t\in[0,T],\\
  \int_0^T(Y^{v}_t-\sup_{b\in U}\{Y^{v\circ(t,b)}_t-\ell(t,X^{v}_t,b))\})dK^{v}_t=0.
\end{cases}
\end{align}
where $X^v$ solves an impulsively controlled forward-SDE (without any continuous control) and the driver, $h$, is stochastic Lipschitz in $z$ with coefficient $L^v$ bounded by $k_L(1+\sup_{s\in [0,\cdot]}|X^v_s|)$ for some $k_L>0$. %see \eqref{ekv:driftless-sde1}-\eqref{ekv:driftless-sde2})

In this section we will consider the case when $\psi$, $h$ and $\ell$ are Lipschitz in $x$ so that, for example, we have $|\psi(X^{v'}_T)-\psi(X^v_T)|\leq C|X^{v'}_T-X^{v}_T|$. When we return to this problem in Section~\ref{sec:robust-impulse}, where a thorough treatment is given, we will extend results obtained in the present section to allow $\psi$, $h$ and $\ell$ to be only locally Lipschitz and of polynomial growth in the $x$-variable.

\subsection{Assumptions}
We introduce the following sets of probability measures on $(\Omega,\mcF)$.
\begin{defn}
We let $\PrM^v:=\PrM^{L^v}$ and define $\mcK^v$ to be the set of all $\mcP_\bbF$-measurable processes $\zeta$ with $|\zeta_t|\leq L^{v}_t$ for all $t\in[0,T]$ (outside of a $\Prob$-null set). Moreover, for all $t\in[0,T]$, we let $\PrM^v_t:=\cup_{u\in\mcU^f_t}\PrM^{v(t)\circ u}$ and $\mcK^v_t:=\cup_{u\in\mcU^f_t}\mcK^{v(t)\circ u}$. We also use the shorthands $\PrM_0:=\PrM^{\emptyset}_0$ and $\mcK_0:=\mcK^{\emptyset}_0$.
\end{defn}

To streamline presentation we will formulate our assumptions on the coefficients in terms of the existence of a family of bounding processes:
\begin{defn}\label{defn:bounding-fam}
We say that a family of processes $(L^v,\Lambda^{\vecv,\vecv',v},\bar K^{v,p},\bar K^{\vecv,\vecv',k,p}: (\vecv,\vecv')\in \cup_{\kappa\geq 1}\mcD^\kappa\times\mcD^\kappa,\,v\in\mcU^f,\,k\geq 0, p\geq 1)$ is a \emph{bounding family} if for each $k\geq 0$ and $p,\kappa\geq 1$, there is a $C>0$ and a $p'\geq 1$ such that for all $v\in\mcU^f$ and $\vecv,\vecv'\in \mcD^\kappa$ and some $q'>1$, we have:
\begin{enumerate}[i)]
  \item\label{b-fam:Lv} $L^v\in\mcS^p$ is non-decreasing, the map $\vecv\mapsto L^{\vecv}$ is continuous from $\mcD^f$ to $\mcH^2$. Moreover, for all $\zeta\in\mcK^v$ and $\bbQ\in\PrM^v$ we have $\E^{\bbQ}[|\mcE(\zeta*W^{\bbQ})_T|^{q'}]\leq C$ (where $W^\bbQ$ is a Brownian motion under $\bbQ$). Moreover, $\E\big[e^{q'\int_0^T|L^v_t|^2dt}\big]\leq C$.% and $\E^{\bbQ}[\mcE(r'\zeta*W^{\bbQ})_T]=1$ for all $r'\in[1,q']$. \todo{Probably not needed.}  and $\vecv\mapsto L^\vecv_T$ is $\Prob$-a.s.~continuous
  \item\label{b-fam:phi} $\Lambda^{\vecv,\vecv',v}\in\mcH^{p}$ is a \cadlag process and the map $(\vecv,\vecv')\mapsto \int_0^T|\Lambda^{\vecv,\vecv',v}_s|^2ds$ is $\Prob$-a.s.~continuous with $\int_0^T|\Lambda^{\vecv,\vecv,v}_s|^2ds=0$.
  \item\label{b-fam:bar-K-v} $\bar K^{v,p}\in\mcS^2$ with $\|\bar K^{v,p}\|_{\mcS^2}^2\leq C$ and $(\bar K^{v,p})^r\leq (\bar K^{v,rp})$ for $r\geq 1$.% and $\bar K^{v,p}\geq 1$, $\Prob$-a.s.
  \item\label{b-fam:bar-K-vp-v} $\bar K^{\vecv,\vecv',k,p}\in\mcS^1$ with $\|\bar K^{\vecv,\vecv',k,p'}\|_{\mcS^1}\leq C\|\vecv'-\vecv\|^p_{\mcD^f}$.
\end{enumerate}
Moreover, for each $r>1$, there is a $C>0$ such that
\begin{align}
  \esssup_{u\in\mcU^f_t}\E\big[|L^{v\circ u}_T|^p\big|\mcF_t\big]&\leq \bar K^{v,p}_t,\label{ekv:L-bound}
  \\
  \esssup_{u\in\mcU^k_t}\E\Big[\big(\int_0^T|\Lambda^{\vecv,\vecv', u}_s|^{2}ds\big)^{p/2}\Big|\mcF_t\Big]&\leq \bar K^{\vecv,\vecv',k,p}_t\label{ekv:phi-bound}
  \\
  \esssup_{u\in\mcU^f_t}\E\Big[\sup_{s\in [t,T]}|\bar K^{v\circ (u(s)),p}_s|^{r}\Big|\mcF_t\Big]&\leq C\bar K^{v,pr}_t,\label{ekv:barK-Doob-equivalent}
\end{align}
for all $\bbQ\in\PrM^{v}_t$.
\end{defn}

\begin{rem}
Going back to the motivation and, in particular, to \eqref{ekv:seq-fbsde}, $L^v$ and $\Lambda^{\vecv,\vecv',v}$ will be processes that satisfy
\begin{align*}
L^v_t&\leq C(1+\sup_{s\in [0,t]}|X^v_s|),
%\\
%\int_{0}^T|L^{v'}_t-L^{v}_t|^2 dt&\leq C\int_0^T|X^{v'}_t-X^{v}_t|^2dt
\\
\Lambda^{\vecv,\vecv',v}_t&=C(|L^{\vecv'\circ v}_t-L^{\vecv\circ v}_t|+|X^{\vecv'\circ v}_t-X^{\vecv\circ v}_t|)
\end{align*}
that will act as Lipschitz coefficients on the driver. The bounding processes on the other hand will take the form
\begin{align*}
\bar K^{v,p}_t&=C(1+\esssup_{u\in\mcU^f_t}\E\Big[\sup_{s\in [0,T]}|X^{v\circ u}_s|^p\Big|\mcF_t\Big]),
\\
\bar K^{\vecv,\vecv',k,p}_t&=C\esssup_{u\in\mcU^k_t}\E\Big[|X^{\vecv'\circ u}_T-X^{\vecv\circ u}_T|^p+\big(\int_0^T|X^{\vecv'\circ u}_s-X^{\vecv\circ u}_s|^2ds\big)^{p/2}
\\
&\quad+|X^{\vecv'\circ [u]_{N-1}}_{\tau_N}-X^{\vecv\circ [u]_{N-1}}_{\tau_N}|^{p}\Big|\mcF_t\Big]
\end{align*}
motivating relations \eqref{ekv:L-bound} and \eqref{ekv:phi-bound} while \eqref{ekv:barK-Doob-equivalent} will follow by Gr\"onwall's inequality combined with Doob's maximal inequality.
\end{rem}

We will make the following assumptions on the involved coefficients:
\begin{ass}\label{ass:on-coeff-FH}
There is a bounding family $(L^v,\Lambda^{\vecv,\vecv',v},\bar K^{v,p},\bar K^{\vecv,\vecv',k,p}: (\vecv,\vecv')\in \cup_{\kappa\geq 1}\mcD^\kappa\times\mcD^\kappa,\,v\in\mcU^f,\,k\geq 0,\,p\geq 1)$ such that for each $k\geq 0$, $p,\kappa\geq 1$, $v\in\mcU^f$ and $\vecv,\vecv'\in\mcD^\kappa$ we have:
\begin{enumerate}[(i)]
  \item\label{ass:varphi} The map $(\omega,\vecv)\mapsto\xi^\vecv:\Omega\times\mcD^f\to\R$ is $\mcF_T\otimes \mcB(\mcD^f)$-measurable, $\Prob$-a.s.~continuous in $\vecv$ and satisfies
      \begin{align}\label{ekv:xi-bound}
        \esssup_{u\in\mcU^f_t}\E\big[|\xi^{v\circ u}|^p\big|\mcF_t\big]\leq \bar K^{v,p}_t
      \end{align}
      and
      \begin{align}
        \esssup_{u\in\mcU^k_t}\E\big[|\xi^{\vecv'\circ u}-\xi^{\vecv\circ u}|^{p}\big|\mcF_t\big]\leq \bar K^{\vecv,\vecv',k,p}_t.
      \end{align}
  \item The intervention cost $c^v$ is such that $(t,b)\mapsto-c^v(t,b)\in\mcO_\bbF$ and satisfies
      \begin{align}
        \inf_{(t,b)\in [0,T]\times U}c^v(t,b)\geq \delta,
      \end{align}
      for some $\delta>0$. Moreover,
      \begin{align}
        \esssup_{u\in\mcU^k_t}\E\big[|c^{\vecv'\circ [u]_{N-1}}(\tau_N,\beta_N)-c^{\vecv\circ [u]_{N-1}}(\tau_N,\beta_N)|^{p}\big|\mcF_t\big]\leq \bar K^{\vecv,\vecv',k,p}_t.
      \end{align}
  \item We have $\xi^{v}\geq \sup_{b\in U}\{\xi^{v\circ(T,b)}-c^v(T,b)\}$, $\Prob$-a.s.
  \item\label{ass:on-f} The map $(\vecv,t,\omega,y,z)\mapsto f^{\vecv}(t,y,z):\mcD^f\times[0,T]\times\Omega\times\R\times \R^{d}\to\R$ is $\mcB(\mcD^f)\otimes\mcP_{\bbF}\otimes\mcB(\R)\otimes\mcB(\R^{d})$-measurable and for each $(y,z)\in\R^{1+d}$ the map $\vecv\mapsto f^\vecv(\cdot,y,z)$ is a continuous map from $\mcD^f$ to $\mcH^2$. Furthermore, we have
  \begin{enumerate}[a)]
  \item the bound
  \begin{align}
  \esssup_{u\in\mcU^f_t}\E\Big[\int_t^{T}|f^{v\circ u}(s,0,0)|^p ds\big|\mcF_t\Big]\leq \bar K^{v,p}_t,
  \end{align}
  \item the Lipschitz condition
      \begin{align}
        |f^{\vecv'\circ u}(t,y',z')-f^{\vecv\circ u}(t,y,z)|\leq k_f|y'-y|+(L^{\vecv\circ u}_t\vee L^{\vecv'\circ u}_t)|z'-z|+(1+|z|+|z'|)\Lambda_t^{\vecv',\vecv,u},
      \end{align}
      for all $(t,\vecv,\vecv',y,y',z,z')\in[0,T]\times \cup_{\kappa\geq 1}(\mcD^\kappa\times\mcD^\kappa)\times\R^{2(1+d)}$, $\Prob$-a.s., for all $u\in\mcU^f$; and
  \item for $u:=(\tau_1,\ldots,\tau_N;\beta_1,\ldots,\beta_N)\in\mcU^f$ we have the causality property
  \begin{align*}
    f^u(t,y,z)=\sum_{j=0}^N\ett_{[\tau_j,\tau_{j+1})}(t)f^{[u]_j}(t,y,z),
  \end{align*}
  where $\tau_0:=0$ and $\tau_{N+1}=\infty$.
%  \item We have
%  \begin{align}
%  \E\Big[\big(\int_0^T |f^{u}(s,0,0)-f^{v}(s,0,0)|ds\big)^2\Big]\leq C\E\big[\|u-v\|^2_{\mcD^f}\big],
%  \end{align}
%  for all $u,v\in\mcU_f$.
  \end{enumerate}
\end{enumerate}
\end{ass}

%\begin{rem}
%Note that in light of the possibility to change between measures, the above assumptions could all have been stated with conditional expectations taken with respect to $\Prob$. However, as will be apparent later on, the above setting is more natural when working with the robust impulse control problem and also more convenient to work with.
%\end{rem}

Before moving on to show existence and uniqueness of solutions to \eqref{ekv:seq-bsde-FH} under Assumption~\ref{ass:on-coeff-FH} we give the following auxiliary result:
\begin{lem}\label{lem:can-Girsanov}
For each $p\geq 1$ there is a $r'>1$ such that $R^v$ defined as
\begin{align*}
  R^v_t:=\esssup_{\zeta\in\mcK^v_t}\E\big[|\mcE(\zeta*W)_{t,T}|^{r'}\big|\mcF_t\big],
\end{align*}
and $\tilde R^v$ defined as
\begin{align*}
  \tilde R^v_t:=\esssup_{\zeta\in\mcK^v_t}\E^{\bbQ^\zeta}\big[|\mcE(-\zeta*W^{\zeta})_{t,T}|^{r'}\big|\mcF_t\big],
\end{align*}
where $W^{\zeta}$ is a $\bbQ^\zeta$-Brownian motion, are both $\mcP_\bbF$-measurable, \cadlag processes such that $\|R^v\|_{\mcS^p}$ and $\|\tilde R^v\|_{\mcS^p}$ are uniformly bounded in $v\in\mcU^f$.
\end{lem}

\noindent\emph{Proof.} By continuity of the map $\vecv\mapsto L^\vecv:\mcD^f\to\mcH^2$ it follows that $R^u$ is continuous on $[\tau_j,\tau_{j+1})$ for $j=0,\ldots,N+1$. For $x\geq 0$, we let $\tau^x:=\inf\{s\geq 0: R^v_s\geq x\}\wedge T$ and note that for $\theta> 1$, we have by right-continuity that
\begin{align*}
\Prob[\sup_{t\in[0,T]}|R^v_t|^{p\theta}\geq x^\theta]&=\Prob[|R^v_{\tau^x}|^{p\theta}\geq x^\theta]\leq \frac{\E\big[|R^v_{\tau^x}|^{p\theta}\big]}{x^\theta}.
\end{align*}
However, for each $\epsilon>0$, there is a $u^\epsilon \in\mcU^f_{\tau^x}$ and a $\zeta^\epsilon$, with $|\zeta^\epsilon_t|\leq \ett_{[\tau^x,T]}(t){L^{v(\tau^x)\circ u^\epsilon}_t}$ such that $R^v_{\tau^x}<\E\big[|\mcE(\zeta^{\epsilon}*W)_{\tau^x,T}|^{r'}\big|\mcF_{\tau^x}\big]+\epsilon$. In particular, it follows that
\begin{align*}
|R^v_{\tau^x}|^{p\theta}&\leq 2^{{p\theta}-1}(\E\Big[\mcE(r'\zeta^{\epsilon}*W)_{\tau^x,T}e^{\frac{(r')^2 -r'}{2}\int_{\tau^x}^T|\zeta^{\epsilon}_s|^2ds}\Big|\mcF_{\tau^x}\Big]^{p\theta} +\epsilon^{p\theta})
\\
&\leq 2^{{p\theta}-1}(\E^{\bbQ^{r'\zeta^{\epsilon}}}\Big[e^{p\theta\frac{(r')^2 -r'}{2}\int_{\tau^x}^T|\zeta^{\epsilon}_s|^2ds}\Big|\mcF_{\tau^x}\Big] +\epsilon^{p\theta})
\end{align*}
which implies that
\begin{align*}
\E\big[|R^v_{\tau^x}|^{p\theta}\big]&\leq 2^{{p\theta}-1}(\E\Big[\mcE(r'\zeta^{\epsilon}*W)_{\tau^x,T}e^{p\theta\frac{(r')^2 -r'}{2}\int_{\tau^x}^T|\zeta^{\epsilon}_s|^2ds}\Big] +\epsilon^{p\theta})
\\
&=2^{{p\theta}-1}(\E^{\bbQ^{r'\zeta^{\epsilon}}}\Big[e^{{p\theta}\frac{(r')^2 -r'}{2}\int_{\tau^x}^T|\zeta^{\epsilon}_s|^2ds}\Big] +\epsilon^{p\theta})
\end{align*}
Since $\epsilon>0$ was arbitrary we find that
\begin{align*}
\Prob\big[\sup_{t\in[0,T]}|R^v_t|^{p}\geq x\big]&=\Prob\big[\sup_{t\in[0,T]}|R^v_t|^{p\theta}\geq x^{\theta}\big]
\\
&\leq 2^{p\theta-1}\sup_{\zeta\in\mcK_0}\frac{\E^{\bbQ^{r'\zeta}}\big[e^{{p\theta}\frac{(r')^2 -r'}{2}\int_{0}^T|\zeta_s|^2ds}\big] }{x^\theta}.
\end{align*}
Now, with $\tilde q'=q'/r'$ and $\tilde q=\tilde q'/(\tilde q'-1)=q'/(q'-r')$ we have
\begin{align*}
\E^{\bbQ^{r'\zeta}}\big[e^{{p\theta}\frac{(r')^2 -r'}{2}\int_{0}^T|\zeta_s|^2ds}\big]\leq \E\big[|\mcE(r'\zeta*W)|^{\tilde q'}\big]^{1/\tilde q'}\E\big[e^{{p\theta\tilde q}\frac{(r')^2 -r'}{2}\int_{0}^T|\zeta_s|^2ds}\big]^{1/\tilde q}.
\end{align*}
For any $\theta>1$ and arbitrary $p\geq 1$ the coefficient ${p\theta\tilde q(r')}\frac{(r')^2 -r'}{2}$ can be made arbitrarily small by choosing $r'>1$ sufficiently small and since there is a $C\geq 0$ such that
\begin{align*}
\E\big[e^{q'\int_{0}^T|\zeta_s|^2ds}\big]\leq \sup_{u\in\mcU^f}\E\big[e^{q'\int_{0}^T|L^u_s|^2ds}\big]\leq C
\end{align*}
by Definition~\ref{defn:bounding-fam}.i) we conclude that
\begin{align*}
\Prob\big[\sup_{t\in[0,T]}|R^v_t|^{p}\geq x\big]&\leq \frac{C}{x^\theta}.
\end{align*}
In particular, using integration by parts, we find that
\begin{align*}
\E\big[\sup_{t\in[0,T]}|R^v_t|^p\big]=\int_0^\infty (\frac{C}{x^\theta}\wedge 1)dx<\infty,
%\int_0^\infty \Prob\big[\sup_{t\in[0,T]}R^v_t\geq x\big]dx<\infty
\end{align*}
showing that $R^v\in\mcS^p$ with norm uniformly bounded in $v$.\\

Concerning $\tilde R^v$ we again have that for each $\epsilon>0$, there is a $u^\epsilon \in\mcU^f_{\tau^x}$ and a $\zeta^\epsilon$, with $|\zeta^\epsilon_t|\leq \ett_{[\tau^x,T]}(t){L^{v(\tau^x)\circ u^\epsilon}_t}$ such that $\tilde R^v_{\tau^x}<\E^{\bbQ^{\zeta^\epsilon}}\big[|\mcE(-\zeta^{\epsilon}*W^{\zeta^\epsilon})_{\tau^x,T}|^{r'}\big|\mcF_{\tau^x}\big]+\epsilon$. Hence,
\begin{align*}
|\tilde R^v_{\tau^x}|^{p\theta}&\leq 2^{{p\theta}-1}(\E^{\bbQ^{\zeta^\epsilon}}\Big[\mcE(-r'\zeta^{\epsilon}*W^{\zeta^\epsilon})_{\tau^x,T}e^{\frac{(r')^2 -r'}{2}\int_{\tau^x}^T|\zeta^{\epsilon}_s|^2ds}\Big|\mcF_{\tau^x}\Big]^{p\theta} +\epsilon^{p\theta})
\\
&\leq 2^{{p\theta}-1}(\E^{\bbQ^{(1-r')\zeta^{\epsilon}}}\Big[e^{p\theta\frac{(r')^2 -r'}{2}\int_{\tau^x}^T|\zeta^{\epsilon}_s|^2ds}\Big|\mcF_{\tau^x}\Big] +\epsilon^{p\theta})
\end{align*}
implying that
\begin{align*}
\E\big[|R^v_{\tau^x}|^{p\theta}\big]&\leq 2^{{p\theta}-1}(\E^{\bbQ^{(1-r')\zeta^{\epsilon}}}\Big[e^{{p\theta}\frac{(r')^2 -r'}{2}\int_{\tau^x}^T|\zeta^{\epsilon}_s|^2ds}\Big] +\epsilon^{p\theta}).
\end{align*}
Consequently,
\begin{align*}
\Prob\big[\sup_{t\in[0,T]}|R^v_t|^{p}\geq x\big]&\leq 2^{p\theta-1}\sup_{\zeta\in\mcK_0}\frac{\E^{\bbQ^{(1-r')\zeta}}\big[e^{{p\theta}\frac{(r')^2 -r'}{2}\int_{0}^T|\zeta_s|^2ds}\big] }{x^\theta}.
\end{align*}
and the assertion follows by repeating the above argument.\qed\\

\begin{lem}\label{lem:can-Girsanov-twice}
For each $p\geq 1$ there is a $r'>1$ such that $\bar R^v$ defined as
\begin{align*}
  \bar R^v_t:=\esssup_{u\in\mcU^f_t}\E\Big[\sup_{s\in[t,T]}\esssup_{\zeta\in\mcK^{v(t)\circ u}_s} \E\big[|\mcE(\zeta*W)_{s,T}|^{r'}\big|\mcF_s\big]\Big|\mcF_t\Big]
\end{align*}
is a $\mcP_\bbF$-measurable, \cadlag process such that $\|\bar R^v\|_{\mcS^p}$ is uniformly bounded in $v\in\mcU^f$.
\end{lem}

\noindent\emph{Proof.} Let $r'>1$ be such that $\|R^v\|_{\mcS^{p\theta}}$ is uniformly bounded in $v\in \mcU^f$ for some $\theta>1$. For $x\geq 0$, we let $\tau^x:=\inf\{s\geq 0: \bar R^v_s\geq x\}\wedge T$ and note that for each $\epsilon>0$, there is a $u^\epsilon \in\mcU^f_{\tau^x}$ such that $\bar R^v_{\tau^x}<\E\Big[\sup_{s\in[\tau^x,T]}R^{v\circ u_\epsilon}_s\big|\mcF_{\tau^x}\Big]+\epsilon$. Jensen's inequality now gives
\begin{align*}
\E\big[|\bar R^v_{\tau^x}|^{p\theta}\big]\leq 2^{p\theta-1}(\E\Big[\sup_{s\in[\tau^x,T]}|R^{v(\tau^x)\circ u_\epsilon}_s|^{p\theta}\Big]+\epsilon^{p\theta}).
\end{align*}
Since $\epsilon>0$ was arbitrary and the first term is bounded by $2^{p\theta-1}\sup_{u\in\mcU^f}\|R^u\|_{\mcS^{p\theta}}$ the result follows by repeating the last steps in the proof of Lemma~\ref{lem:can-Girsanov}.\qed\\

Throughout this section, we assume that $r'>1$ is small enough that $\|R^v\|_{\mcS^3}$, $\|\tilde R^v\|_{\mcS^3}$ and $\|\bar R^v\|_{\mcS^3}$ are bounded uniformly in $v\in\mcU^f$ and let $r$ be such that $\frac{1}{r'}+\frac{1}{r}=1$.

\subsection{An approximating sequence}
In this section we outline a Piccard type approximation scheme, that will ultimately lead us to the conclusion that \eqref{ekv:seq-bsde-FH} has a solution under Assumption~\ref{ass:on-coeff-FH}. We note that for all $v\in\mcU^f$,
\begin{align}\label{ekv:bsde_0}
  Y^{v,0}_t=\xi^{v}+\int_t^T f^v(s,Y^{v,0}_s,Z^{v,0}_s)ds-\int_t^T Z^{v,0}_sdW_s
\end{align}
admits a unique solution by Theorem~\ref{thm:bsde-no-reflection}. We, thus, consider the following sequence of families of reflected BSDEs
\begin{align}\label{ekv:rbsde_k}
\begin{cases}
  Y^{v,k}_t=\xi^{v}+\int_t^T f^v(s,Y^{v,k}_s,Z^{v,k}_s)ds-\int_t^T Z^{v,k}_sdW_s+ K^{v,k}_T-K^{v,k}_t \\
  Y^{v,k}_t\geq\sup_{b\in U}\{Y^{v\circ(t,b),k-1}_t-c^v(t,b)\}\\
  \int_0^T(Y^{v,k}_t-\sup_{b\in U}\{Y^{v\circ(t,b),k-1}_t-c^v(t,b)\})dK^{v,k}_t=0.
\end{cases}
\end{align}
for $k\geq 1$.

We will make use of the following induction hypothesis:

\begin{hyp*}[{\bf RBSDE.$l$}] There is a family of pairs $((Y^{v,0},Z^{v,0}): v\in \mcU^f)$ and a sequence of families of triples $((Y^{v,k},Z^{v,k},K^{v,k}): v\in \mcU^f)_{1\leq k\leq l}$ such that:
\begin{enumerate}[i)]
\item\label{Yk:rec} For each $v\in \mcU^f$, the pair $(Y^{v,0},Z^{v,0})\in\mcS^2\times\mcH^2$ solves \eqref{ekv:bsde_0} and the triple $(Y^{v,k},Z^{v,k},K^{v,k})\in\mcS^2\times\mcH^2\times\mcS^2$ solves \eqref{ekv:rbsde_k} for $k=1,\ldots,l$.
%\item There is a constant $C>0$ (that does not dependent of $l$) such that
%\begin{align*}
%  \sup_{v\in\mcU^f}(\|Y^{v,k}\|_{\mcS^2}^2+\|Z^{v,k}\|_{\mcH^2}^2+\|K^{v,k}\|_{\mcS^2}^2)\leq C,
%\end{align*}
%for $k=0,\ldots,l$.
\item\label{Yk:consist} For all $k\in\{0,\ldots,l\}$ and each $j\geq 0$, the map $h^{k,j}:[0,T]\times \mcD^j\to L^2(\Omega,\mcF,\Prob)\,((t,\vecv)\mapsto Y^{\vecv,k}_t)$ is jointly continuous (outside of a $\Prob$-null set) and, moreover, for each $v\in\mcU^f$ we have $Y^{v,k}_t=\sum_{j=0}^\infty\ett_{[N=j]}h^{k,j}(t,v)$ for all $t\in[0,T]$ outside of a $\Prob$-null set.
\end{enumerate}
\end{hyp*}
Recall Definition~\ref{defn:rbsde-solution} specifying what we mean by a solution to \eqref{ekv:bsde_0} and \eqref{ekv:rbsde_k}. Through this definition the first condition in Hypothesis {\bf RBSDE.$l$} implies that $\sup_{u\in\mcU^f}\|Y^{v,k}\|<\infty$ and also dictates the regularity of the map $(t,b)\mapsto Y^{v\circ(t,b)}_t$. The second statement, on the other hand, is a stronger version of the consistency property for families of processes introduced in Definition~\ref{defn:consistency}. To simplify presentation we will refer to the second property as \emph{strong consistency}.

\begin{prop}\label{prop:seq-bsde_0-FH}
Hypothesis {\bf RBSDE.$0$} holds.
\end{prop}

\noindent\emph{Proof.} Existence of solutions to \eqref{ekv:bsde_0} with $Y^v\in\mcS^2$ and $Z^v\in\mcH^2$ and uniform boundedness (specified in \emph{(i)} of Definition~\ref{defn:rbsde-solution}) follows from Theorem~\ref{thm:bsde-no-reflection}.

For $\kappa\geq 1$ and $\vecv,\vecv'\in\mcD^\kappa$ we have by Theorem~\ref{thm:bsde-no-reflection}, Assumption~\ref{ass:on-coeff-FH} and Definition~\ref{defn:bounding-fam}.(\ref{b-fam:bar-K-vp-v}) that for each $p\geq 1$, there is a $p'\geq 1$ such that
\begin{align*}
&\|Y^{\vecv',0}-Y^{\vecv,0}\|_{\mcS^{p'}}^{p'}+\|Z^{\vecv',0}-Z^{\vecv,0}\|_{\mcH^{p'}}^{p'}
\\
&\leq C\E\Big[|\xi^{\vecv'}-\xi^{\vecv}|^{q^2p'}+\big(\int_0^T |f^{\vecv'}(s,Y^{\vecv,0}_s,Z^{\vecv,0}_s)-f^{\vecv}(s,Y^{\vecv,0}_s,Z^{\vecv,0}_s)|ds\big)^{q^2p'}\Big]^{1/q^2}
\\
&\leq C\E\Big[|\xi^{\vecv'}-\xi^{\vecv}|^{q^2p'}+\big(\int_0^T (1+|Z^{\vecv,0}_s|)\Lambda^{\vecv,\vecv',\emptyset}_sds\big)^{q^2p'}\Big]^{1/q^2}
\\
&\leq C\E\Big[|\xi^{\vecv'}-\xi^{\vecv}|^{q^2p'}+\big(\int_0^T (1+|Z^{\vecv,0}_s|^2)ds\int_0^T |\Lambda^{\vecv,\vecv',\emptyset}_s|^2ds\big)^{q^2p'/2}\Big]^{1/q^2}
\\
&\leq C(\E\Big[|\xi^{\vecv'}-\xi^{\vecv}|^{q^2p'}\Big]^{1/q^2}+(1+\E\Big[\big(\int_0^T |Z^{\vecv,0}_s|^{2}ds\big)^{q^2p'}\Big]^{1/2q^2})\E\Big[\big(\int_0^T|\Lambda^{\vecv,\vecv',\emptyset}_s|^2ds\big)^{q^2p'}ds\Big]^{1/2q^2}
\\
&\leq C\|\vecv'-\vecv\|_{\mcD^f}^{p},
\end{align*}
where $C>0$ does not depend on $\vecv,\vecv'$. By picking $p=(m+1)\kappa+1$, Kolmogorov's continuity theorem (see \eg Theorem 72 in Chapter IV of~\cite{Protter}) guarantees the existence of a family of processes $(\hat Y^{\vecv},\hat Z^{\vecv}:\vecv\in\mcD^\kappa)$, where for each $\vecv\in\mcD^\kappa$, $\hat Y^{\vecv}\in\mcS^2$ and $\hat Z^{\vecv}\in\mcH^2$, such that, outside of a $\Prob$-null set, $(s,\vecv)\mapsto\hat Y^{\vecv}_s$ is continuous in $\vecv$ uniformly in $s$ and $\vecv\mapsto\hat Z^{\vecv}_\cdot$ is $L^2([0,T])$-continuous (and in particular that $(s,\vecv)\mapsto\int_s^TZ^{\vecv}_rdW_r$ is continuous in $\vecv$ uniformly in $s$) and moreover $(\hat Y^{\vecv}_s,\hat Z^{\vecv}_s)=(Y^{\vecv,0}_s,Z^{\vecv,0}_s)$, $\Prob$-a.s.

Now, since this holds for all $\kappa\geq 0$ it is clear that for each $\vecv\in\mcD^f$, the pair $(\hat Y^{\vecv}_s,\hat Z^{\vecv}_s:s\in[0,T])$ solves the corresponding BSDE \eqref{ekv:bsde_0} and that the map $\vecv\mapsto (\hat Y^{\vecv},\hat Z^{\vecv}):\mcD^f\to\mcS^2\times\mcH^2$ is continuous.

To establish the strong consistency in Hypothesis {\bf RBSDE.$0$}-\emph{(\ref{Yk:consist})}, we need to show that for any $v\in\mcU^f$, the pair $(\hat Y^{v},\hat Z^{v})\in\mcS^2\times\mcH^2$ solves the BSDE corresponding to the control $v$. We let $(v_j)_{j\geq 0}$ be an approximating sequence in $\mcU^f$ (\ie $v_j\in \mcU^f$ and $v_j\to v$, $\Prob$-a.s.~as $j\to\infty$) taking values in a countable subset of $\mcD^f$. %so that and $\tau_j$ are stopping times such that $\tau_j\searrow\tau$, $\Prob$-a.s.
By continuity we have
\begin{align*}
\sup_{s\in [0,T]}|\hat Y^{v_j,0}_s-\hat Y^{v,0}_s|\to 0,
\end{align*}
$\Prob$-a.s., as $j\to\infty$. Moreover, due to continuity of the map $\vecv\mapsto \xi^{\vecv}$ we have
\begin{align*}
  \xi^{v_j}\to \xi^{v}
\end{align*}
and by Assumption~\ref{ass:on-coeff-FH} we get
\begin{align*}
\int_0^T |f^{v_j}(s,\hat Y^{v_j,0}_s,\hat Z^{v_j,0}_s)-f^{v}(s,\hat Y^{v,0}_s,\hat Z^{v,0}_s)|ds&\leq \int_0^T |f^{v_j}(s,\hat Y^{v_j,0}_s,\hat Z^{v_j,0}_s)-f^{v}(s,\hat Y^{v_j,0}_s,\hat Z^{v_j,0}_s)|ds+
\\
&\quad +\int_0^T |f^{v}(s,\hat Y^{v_j,0}_s,\hat Z^{v_j,0}_s)-f^{v}(s,\hat Y^{v,0}_s,\hat Z^{v,0}_s)|ds
\\
& \leq \int_0^T (1+2|\hat Z^{v_j,0}_s|)\Lambda_s^{v_j,v,\emptyset}ds
\\
&\quad + \int_0^T (k_f|\hat Y^{v_j,0}_s-\hat Y^{v,0}_s|+L^{v}_s|\hat Z^{v_j,0}_s-\hat Z^{v,0}_s|)ds
\end{align*}
%\begin{align*}
%\int_0^T |f^{v_j}(s,\hat Y^{v_j,0}_s,\hat Z^{v_j,0}_s)-f^{v}(s,\hat Y^{v,0}_s,\hat Z^{v,0}_s)|ds&\leq \int_0^T |f^{v_j}(s,\hat Y^{v_j,0}_s,\hat Z^{v_j,0}_s)-f^{v_j}(s,\hat Y^{v,0}_s,\hat Z^{v,0}_s)|ds+
%\\
%&\quad +\int_0^T |f^{v}(s,\hat Y^{v,0}_s,\hat Z^{v,0}_s)-f^{v_j}(s,\hat Y^{v,0}_s,\hat Z^{v,0}_s)|ds
%\\
%& \leq \int_0^T (k_f|\hat Y^{v_j,0}_s-\hat Y^{v,0}_s|+L^{v_j}_s|\hat Z^{v_j,0}_s-\hat Z^{v,0}_s|)ds
%\\
%&\quad +\int_0^T |(1+|\hat Z^{v,0}_s|)\Lambda_s^{v_j,v,\emptyset}ds
%\end{align*}
which tends to 0, $\Prob$-a.s., as $j\to\infty$. Moreover, for each $j\geq 0$, the pair $(\hat Y^{v_j},\hat Z^{v_j})$ solves \eqref{ekv:bsde_0} with control $v_j$ and we conclude that
\begin{align*}
\hat Y^{v,0}_\eta&=\lim_{j\to\infty}\Big\{\xi^{v_j}+\int_\eta^T f^{v_j}(s,\hat Y^{v_j,0}_s,\hat Z^{v_j,0}_s)ds-\int_\eta^T \hat Z^{v_j,0}_sdW_s\Big\}
\\
&=\xi^{v}+\int_\eta^T f^{v}(s,\hat Y^{v,0}_s,\hat Z^{v,0}_s)ds-\int_\eta^T \hat Z^{v,0}_sdW_s,
\end{align*}
for each $\eta\in\mcT$ and strong consistency follows.\qed\\

We now turn to the reflected BSDEs \eqref{ekv:rbsde_k}. To obtain estimates for the triple $(Y^{v,k},Z^{v,k},K^{v,k})$ we rely on Corollary~\ref{cor:rbsde-char-FH} to reduce the system of reflected BSDEs to a single non-reflected BSDE with jumps. We, thus, introduce the following BSDE:

\begin{defn}\label{defn:U-V}
For $v,u\in\mcU^f$, let the pair $(U^{v,u}, V^{v,u})\in\mcS^2_l\times\mcH^2$ (recall that $\mcS^2_l$ is the set of $\mcP_\bbF$-measurable \caglad processes with finite $\mcS^2$-norm) be the unique solution to the BSDE
\begin{align}\label{ekv:U-def}
 U_t^{v,u}&=\xi^{v\circ u}+\int_t^{T}f^{v\circ u}(s, U^{v,u}_s, V^{v,u}_s)ds
-\int_t^{T} V^{v,u}_sdW_s-\sum_{j=1}^N \ett_{[\tau_j\geq t]}c^{v\circ [u]_{j-1}}(\tau_j,\beta_j),
\end{align}
whenever a unique solution exists and let $U^{v,u}\equiv -\infty$, otherwise.
\end{defn}

\bigskip

%Expressing the process $Y^{v,k}$ in terms of the pair $(U^{v,u}, V^{v,u})\in\mcS^2\times\mcH^2$ for suitable $u$ will be a basis for obtaining estimates throughout this section and

\begin{prop}\label{prop:UV-has-solution}
For each $k\geq 0$, $v\in\mcU^f$ and $u\in\mcU^k$ the BSDE \eqref{ekv:U-def} admits a unique solution and $V^{v,u}\in\mcH^2_\bbQ$ for all $\bbQ\in\PrM^v$.
\end{prop}

\noindent\emph{Proof.} Existence of a unique solution to \eqref{ekv:U-def} follows from repeated use of Theorem~\ref{thm:bsde-no-reflection} since the intervention costs belong to $L^p(\Omega,\mcF,\Prob)$ for all $p\geq 1$. Moreover, a similar argument gives that
\begin{align*}
\|V^{v,u}\|_{\mcH^p}^p&\leq C\E\Big[|\xi^{v\circ u}|^{q^2 p}+\int_0^T |f^{v\circ u}(s,0,0)|^{q^2 p}ds+\sum_{j=1}^N |c^{v\circ [u]_{j-1}}(\tau_j,\beta_j)|^{q^2 p}\Big]^{1/q^2}
\\
&\leq C.
\end{align*}
Now,
\begin{align*}
\|V^{u,v}\|^2_{\mcH^2_\bbQ} &= \E\Big[\mcE(\zeta*W)_T\int_0^T |V^{u,v}_s|^2ds\Big]
\\
&\leq \E\big[|\mcE(\zeta*W)_T|^{q'}\big]^{1/q'}\|V^{v,u}\|_{\mcH^q}^2
\end{align*}
and the assertion follows.\qed\\

In addition, we introduce the following notation:

\begin{defn}\label{defn:bbQ-gamma}
For $v,u\in\mcU^f$ such that \eqref{ekv:U-def} admits a unique solution, we define
\begin{align*}
\gamma^{v,u}_s:=\frac{f^{v\circ u}(s,U^{v,u}_s,V^{v,u})-f^{v\circ u}(s,0,V^{v,u})}{U^{v,u}_s}\ett_{[U^{v,u}_s\neq 0]}
\end{align*}
and for $0\leq s\leq t\leq T$, we set $e^{v,u}_{s,t}:=e^{\int_s^t\gamma^{v,u}_rdr}$ and $e^{v,u}_t:=e^{v,u}_{0,t}$. Moreover, we define
\begin{align*}
\zeta^{v,u}:=\frac{f^{v\circ u}(s,0,V^{v,u})-f^{v\circ u}(s,0,0)}{|V^{v,u}_s|^2}(V^{v,u}_s)^\top\ett_{[V^{v,u}_s\neq 0]}
\end{align*}
and let $\bbQ^{v,u}:=\bbQ^{\zeta^{v,u}}$, the probability measure, equivalent to $\Prob$, under which $W^{v,u}:=W-\int_0^\cdot \zeta^{v,u}_sds$ is a Brownian motion.
\end{defn}

\bigskip

%%%%%%%%%%%%%%%%%%%%%%%%%%%%%%%%%%%%%%%%%%%%%%%%%%%%%%%%%%%%%%%%%%%%%%%%%%%%%%%%%%%%%%%%%%%%%%%%%%%%%%%%%%%%%%%%%%%%%%%%%%%%%%%%%%%%%%%%%%%%%%%%%%

Before we move on to show that Hypothesis {\bf RBSDE.$l$} holds for all $l\geq 0$ we give three helpful lemmas.

\begin{lem}\label{lem:seq-diag-bound-FH}
Assume that Hypothesis {\bf RBSDE.l} holds for some $l\geq 0$, then for each $p\geq 1$, there is a $C>0$ (that does not depend on $l$) such that
\begin{align}\label{ekv:seq-diag-bound-FH}
\|(\sup_{b\in U}|Y^{v\circ (t,b),l}_{t}|: t\in[0,T])\|_{\mcS^p}\leq C.
\end{align}
\end{lem}

\noindent\emph{Proof.} We let $\tilde v:=v\circ(t,b)$ and set $\tau^*_1:=\inf \{s\geq t:Y^{\tilde v,l}_s=\sup_{b\in U}\{Y^{\tilde v\circ(s,b),l-1}_s-c^{\tilde v}(s,b)\}\}\wedge T$ and have by Corollary~\ref{cor:rbsde-char-FH} and consistency that
\begin{align*}
Y^{\tilde v,l}_t&=\ett_{[\tau^*_1<T]}\sup_{b'\in U}\{Y^{\tilde v\circ(\tau^*_1,b'),l-1}_{\tau^*_1}-c^{\tilde v}(\tau^*_1,b')\}+\ett_{[\tau^*_1=T]}\xi^{\tilde v}
\\
&\quad+\int_t^{\tau^*_1}f^{\tilde v}(s,Y^{\tilde v,l}_s,Z^{\tilde v,l}_s)ds-\int_t^{\tau^*_1}Z^{\tilde v,l}_sdW_s
\\
&=\ett_{[\tau^*_1<T]}\{Y^{\tilde v\circ(\tau^*_1,\beta^*_1),l-1}_{\tau^*_1}-c^{\tilde v}(\tau^*_1,\beta^*_1)\}+\ett_{[\tau^*_1=T]}\xi^{\tilde v}
\\
&\quad+\int_t^{\tau^*_1}f^{\tilde v}(s,Y^{\tilde v,l}_s,Z^{\tilde v,l}_s)ds-\int_t^{\tau^*_1}Z^{\tilde v,l}_sdW_s,
\end{align*}
where $\beta^*_1$ can be chosen to be $\mcF_{\tau^*_1}$-measurable by continuity of the map $b'\mapsto Y^{\tilde v\circ(\tau^*_1,b'),l-1}_{\tau^*_1}-c^{\tilde v}(\tau^*_1,b')$ and the measurable selection theorem (see \eg Chapter 7 in \cite{BertsekasShreve} or \cite{ElKarouiTan}).

Now, we can continue and inductively define $\tau^*_{j}:=\inf \{s\geq \tau_{j-1}^*:Y^{\tilde v\circ(\tau^*_1,\ldots,\tau^*_{j-1},\beta^*_1,\ldots,\beta^*_{j-1}),l+1-j}_s=\sup_{b\in U}\{Y^{\tilde v\circ(\tau^*_1,\ldots,\tau^*_{j-1},\beta^*_1,\ldots,\beta^*_{j-1})\circ(s,b),l-j}_s-c^{\tilde v\circ(\tau^*_1,\ldots,\tau^*_{j-1},\beta^*_1,\ldots,\beta^*_{j-1})}(s,b)\}\}\wedge T$ for $j=1,\ldots, l$, and take $\beta^*_j$ to be the corresponding $\mcF_{\tau_j^*}$-measurable maximizer. By induction we get that
\begin{align}\nonumber
Y^{\tilde v,l}_t&=\xi^{\tilde v\circ u^*}+\int_t^{T}\sum_{j=0}^{N^*}\ett_{[\tau^*_j,\tau^*_{j+1})}(s)f^{\tilde v\circ [u^*]_j}(s,Y^{\tilde v\circ [u^*]_j,l-j}_s,Z^{\tilde v\circ [u^*]_j,l-j}_s)ds
\\
&\quad-\int_t^{T}\sum_{j=0}^{N^*}\ett_{[\tau^*_j,\tau^*_{j+1})}(s)Z^{\tilde v\circ [u^*]_j,l-j}_sdW_s-\sum_{j=1}^{N^*}c^{\tilde v\circ [u^*]_{j-1}}(\tau^*_j,\beta^*_j),\label{ekv:Y_k-alter-def}
\end{align}
where $u^*:=(\tau_1^*,\ldots,\tau^*_{N^*};\beta^*_1,\ldots,\beta^*_{N^*})$ with $N^*:=\max\{j\in\{0,\ldots,l\}:\tau^*_j<T\}$ and using the convention that $\tau^*_0=0$ and $\tau^*_{N^*+1}=T$.

In particular, \eqref{ekv:Y_k-alter-def} implies by comparison and positivity of the intervention cost that $U^{\tilde v,\emptyset}_t\leq Y^{\tilde v,l}_t\leq \esssup_{u\in\mcU^{l}_t}U^{\tilde v\circ u,\emptyset}_t$, $\Prob$-a.s., and we find that $|Y^{\tilde v,l}_t|\leq \esssup_{u\in\mcU^{l}_t}|U^{\tilde v\circ u,\emptyset}_t|$. Furthermore, by continuity of the map $b\mapsto Y^{v\circ(t,b)}_t$ we can find a $\beta^\diamond\in\mcI(t)$ such that $\sup_{b\in U}|Y^{v\circ(t,b)}_t|=|Y^{v\circ(t,\beta^\diamond)}_t|$, $\Prob$-a.s., and we conclude that $\sup_{b\in U}|Y^{v\circ(t,b)}_t|\leq \esssup_{u\in\mcU^{l+1}_t}|U^{v\circ u,\emptyset}_t|$. Now, for arbitrary $u\in\mcU^f$ we have that
\begin{align*}
U^{v\circ u,\emptyset}_t%&=\xi^{v\circ u}+\int_t^{T}f^{v\circ u}(s,Y^{v\circ u,0}_s,Z^{v\circ u,0}_s)ds-\int_t^{T}Z^{v\circ u,0}_sdW_s
%\\
&=e^{v\circ u,\emptyset}_{t,T}\xi^{v\circ u}+\int_t^{T}e^{v\circ u,\emptyset}_{t,s}f^{v\circ u}(s,0,0)ds-\int_t^{T}e^{v\circ u,\emptyset}_{t,s}V^{v\circ u,\emptyset}_sdW^{v\circ u,\emptyset}_s.
\end{align*}
This gives, since $V^{v\circ u,\emptyset}\in\mcH^2_{\bbQ^{v\circ u,\emptyset}}$ (see Proposition~\ref{prop:UV-has-solution}), that
\begin{align}\nonumber
|U^{v\circ u,\emptyset}_t|^p&\leq C\E^{\bbQ^{v\circ u,\emptyset}}\Big[|\xi^{v\circ u}|^p+\int_t^{T}|f^{v\circ u}(s,0,0)|^pds\big|\mcF_t\Big]
\\
&\leq C\E\big[|\mcE(\zeta^{v\circ u,\emptyset}*W)_{t,T}|^{r'}\big|\mcF_t\big]^{1/r'}\E\Big[\big(|\xi^{v\circ u}|^p+\int_t^{T}|f^{v\circ u}(s,0,0)|^pds\big)^{r}\big|\mcF_t\Big]^{1/r}\nonumber
\\
&\leq C (R^v_t)^{1/r'}(\bar K^{v,pr}_t)^{1/r},\label{ekv:U-vu-emptyset-bound}
\end{align}
where the last inequality follows by Assumption~\ref{ass:on-coeff-FH} since $|\zeta^{v\circ u,\emptyset}|\leq L^{v\circ u}$. In particular, continuity of $(\sup_{b\in U}Y^{v\circ (t,b),l}_{t}:t\in[0,T])$ implies that, outside of a $\Prob$-null set, we have
\begin{align*}
\sup_{b\in U}|Y^{v\circ (t,b),l}_{t}|^p\leq C (R^v_t)^{1/r'}(\bar K^{v,pr}_t)^{1/r}
\end{align*}
for all $t\in[0,T]$. We can thus apply H\"older's inequality to find that
\begin{align*}
\|(\sup_{b\in U}Y^{v\circ (t,b),l}_{t}: t\in[0,T])\|_{\mcS^p}^p\leq C\|(R^v)^{1/r'}(\bar K^{v,pr})^{1/r}\|_{\mcS^1}\leq C\|R^v\|_{\mcS^1}^{1/r'}\|\bar K^{v,pr}\|_{\mcS^1}^{1/r}\leq C
\end{align*}
by Lemma~\ref{lem:can-Girsanov} and since $\bar K^{v,p}_t\in\mcS^1$ for all $p\geq 1$.\qed\\

\begin{lem}\label{lem:seq-unif-bound-FH}
Assume that Hypothesis {\bf RBSDE.$l$} holds for some $l\geq 0$, then for each $p\geq 1$ there is a $C>0$, that does not depend on $l$, such that
\begin{align}\label{ekv:seq-unif-bound-FH}
\|Y^{v,l+1}\|_{\mcS^p}+\|Z^{v,l+1}\|_{\mcH^p}+\|K^{v,l+1}\|_{\mcS^p}\leq C,
\end{align}
for all $v\in\mcU^f$.
\end{lem}

\noindent\emph{Proof.} This is immediate from Proposition~\ref{prop:rbsde-solu-FH} and Lemma~\ref{lem:seq-diag-bound-FH}.\qed\\

%%%%%%%%%%%%%%%%%%%%%%%%%%%%%%%%%%%%%%%%%%%%%%%%%%%%%%%%%%%%%%%%%%%%%%%%%%%%%%%%%%%%%%%%%%%%%%%%%%%%%%%%%%%%%%%%%%%%%%%%%%%%%%%%%%%%%%%%%%%%%%%%%%

\begin{lem}\label{lem:seq-diff-FH}
Assume that Hypothesis {\bf RBSDE.l} holds for some $l\geq 0$, then for each $\kappa\geq 0$ and $p\geq 1$, there is a $C>0$ and a $p'\geq 1$ such that for any $\vecv,\vecv'\in\mcD^\kappa$, we have
\begin{align}\label{ekv:seq-diff-FH}
\|(\sup_{b\in U}|Y^{\vecv'\circ(t,b),l}_t-Y^{\vecv\circ(t,b),l}_t|:t\in[0,T])\|_{\mcS^{p'}}^{p'}\leq C\|\vecv'-\vecv\|_{\mcD^f}^{p}.
\end{align}
\end{lem}

\noindent\emph{Proof.} We let $u^*$ and $u':=(\tau'_1,\ldots,\tau'_{N'};\beta'_1,\ldots,\beta'_{N'})$ be the controls obtained by repeating the construction in the proof of Lemma~\ref{lem:seq-unif-bound-FH}, starting from $Y^{\vecv\circ(t,b),l}_t$ and $Y^{\vecv'\circ(t,b),l}_t$, respectively, instead of $Y^{v\circ(t,b),l}_t$. By Proposition~\ref{prop:UV-has-solution} it follows that for each $v\in\mcU^f$ and $u\in\mcU^{l+1}$, there is a unique pair $(U^{v,u},V^{v,u})\in\mcS^2\times\mcH^2$ that solves \eqref{ekv:U-def}, \ie
\begin{align*}
U_t^{v,u}&=\xi^{v\circ u}+\int_t^{T}f^{v\circ u}(s,U^{v,u}_s,V^{v,u}_s)ds
-\int_t^{T}V^{v,u}_sdW_s-\sum_{j=1}^N \ett_{[\tau_j\geq t]}c^{v\circ [u]_{j-1}}(\tau_j,\beta_j).
\end{align*}
By a trivial argument we find that $Y^{\vecv\circ(t,b),l}_t=U_t^{\vecv\circ(t,b),u^*}=\esssup_{u\in\mcU^l}U_t^{\vecv\circ(t,b),u}$ and conclude that
\begin{align*}
|Y^{\vecv\circ(t,b),l}_t-Y^{\vecv'\circ(t,b),l}_t|&\leq |U_t^{\vecv\circ(t,b),u^*}-U_t^{\vecv'\circ(t,b),u^*}|+|U_t^{\vecv\circ(t,b),u'} - U_t^{\vecv'\circ(t,b),u'}|
\\
&\leq 2\esssup_{u\in\mcU^{l+1}} |U_t^{\vecv,u}-U_t^{\vecv',u}|
\end{align*}
Letting $\delta U^u:=U^{\vecv',u}-U^{\vecv,u}$ and $\delta V^u:=V^{\vecv',u}-V^{\vecv,u}$ and writing $\delta \Box^u:=\Box^{\vecv'\circ u}-\Box^{\vecv\circ u}$ (for $\Box=\xi,f$ and $c$) gives
\begin{align}
e_t\delta U^u_t&=e_T\delta\xi^{u}+\int_t^{T}e_s\delta f^{u}(s,U^{\vecv,u}_s,V^{\vecv,u}_s)ds
-\int_t^{T}e_s\delta V^{u}dW^\zeta_s-\sum_{j=1}^N e_{\tau_j}\delta c^{[u]_{j-1}}(\tau_j,\beta_j),\label{ekv:U-V-delta}
\end{align}
where $e_t:=e^{\int_0^t \gamma_sds}$ with
\begin{align*}
\gamma_s:=\frac{f^{\vecv'\circ u}(s,U^{\vecv',u}_s,V^{\vecv',u}_s)-f^{\vecv'\circ u}(s,U^{\vecv,u}_s,V^{\vecv',u}_s)}{U^{\vecv',u}_s-U^{\vecv,u}_s}\ett_{[U^{\vecv',u}_s\neq U^{\vecv,u}_s]}
\end{align*}
and $W^\zeta:=W-\int_0^\cdot \zeta_sds$, with
\begin{align*}
\zeta_s:=\frac{f^{\vecv'\circ u}(s,U^{\vecv,u}_s,V^{\vecv',u}_s)-f^{\vecv'\circ u}(s,U^{\vecv,u}_s,V^{\vecv,u}_s)}{|V^{\vecv',u}_s-V^{\vecv,u}_s|^2}(V^{\vecv',u}_s-V^{\vecv,u}_s)^\top\ett_{[V^{\vecv',u}_s\neq V^{\vecv,u}_s]},
\end{align*}
is a Brownian motion under the measure $\bbQ^\zeta\in\PrM^{\vecv'}_t$ given by $d\bbQ^\zeta=\mcE(\zeta*W)_T d\Prob$.

For $u\in\mcU^{l+1}_t$, by again appealing to Proposition~\ref{prop:UV-has-solution}, we have that $\delta V^{u}\in \mcH^2_{\bbQ^\zeta}$ and since $e^{-k_fT}\leq e_t\leq e^{k_fT}$, taking conditional expectation in \eqref{ekv:U-V-delta} gives
\begin{align*}
|\delta U^u_t|&\leq C\E^{\bbQ^\zeta}\Big[|\delta\xi^u|+\int_0^T|\delta f^u(s,U^{\vecv,u}_s,V^{\vecv,u}_s)| ds+\sum_{j=1}^N |\delta c^{[u]_{j-1}}(\tau_j,\beta_j)|\Big|\mcF_t\Big]
\\
&\leq C\E^{\bbQ^\zeta}\Big[|\delta\xi^u|+(\int_t^T|\Lambda_s^{\vecv',\vecv,u}|^2ds)^{1/2}(1+\int_t^T|V^{\vecv,u}_s|^2ds)^{1/2}+\sum_{j=1}^N |\delta c^{[u]_{j-1}}(\tau_j,\beta_j)|\Big|\mcF_t\Big],
%\\
%&\leq C(\bar K^{\vecv,\vecv'}_t+\E^{\bbQ^\zeta}\Big[\Lambda_T^{\vecv',\vecv,u}(1+\int_t^T|V^{\vecv,u}_s|ds)\big|\mcF_t\Big])
%\\
%&\leq C(\bar K^{\vecv,\vecv'}_t+\E^{\bbQ^\zeta}\big[|\Lambda_T^{\vecv',\vecv,u}|^2\big|\mcF_t\big]^{1/2}(1 + \E^{\bbQ^{\zeta}}\Big[\int_t^T|V^{\vecv,u}_s|^2ds\big|\mcF_t\Big]^{1/2})).
\end{align*}
where we have used the Lipschitz condition on $f$ to arrive at the last inequality. In particular, since $N\leq l+1$, $\Prob$-a.s., this gives that
\begin{align*}
|\delta U^u_t|^p&\leq C\E^{\bbQ^\zeta}\Big[|\delta\xi^u|^p+(\int_t^T|\Lambda_s^{\vecv',\vecv,u}|^2ds)^{p/2}(1+\big(\int_t^T|V^{\vecv,u}_s|^2ds\big)^{p/2})+\sum_{j=1}^N |\delta c^{[u]_{j-1}}(\tau_j,\beta_j)|^p\Big|\mcF_t\Big]
\\
&\leq C(R^{\vecv'}_t)^{1/r'}((\bar K^{\vecv,\vecv',l+1,pr}_t)^{1/r}+(\bar K^{\vecv,\vecv',l+1,2pr}_t)^{1/2r}(1+\E\Big[\big(\int_t^T|V^{\vecv,u}_s|^2ds\big)^{pr}\Big|\mcF_t\Big]^{1/2r})),
\end{align*}
as $|\zeta_t|\leq L^{\vecv'\circ u}_t$. By arguing as in the proof of Lemma~\ref{lem:rbsde-trunk-diff} under a conditional setting it follows that
\begin{align*}
\E^{\bbQ^{\vecv,u}}\Big[\big(\int_t^T|V^{\vecv,u}_s|^2 ds\big)^{p}\big|\mcF_t\Big]&\leq C\E^{\bbQ^{\vecv,u}}\Big[|\xi^{\vecv\circ u}|^{2p}+\big(\int_t^T|f^{\vecv\circ u}(s,0,0)|ds\big)^{2p}
\\
&\quad+\sum_{j=1}^N |c^{\vecv\circ[u]_{j-1}}(\tau_j,\beta_j)|^{2p}\Big|\mcF_t\Big].
\end{align*}
%where $\bbQ^{\vecv,u}$ is the measure defined by $d\bbQ^{\vecv,u}=\mcE(\zeta^{\vecv,u}*W)d\Prob$ with
%\begin{align*}
%\zeta^{\vecv,u}_s:=\frac{f^{\vecv\circ u}(s,0,V^{\vecv,u}_s)-f^{\vecv\circ u}(s,0,0)}{|V^{\vecv,u}_s|^2}(V^{\vecv,u}_s)^\top\ett_{[V^{\vecv,u}_s\neq 0]}.
%\end{align*}
Changing back to $\Prob$-expectation we get by the Girsanov theorem that
\begin{align*}
\E\Big[\big(\int_t^T|V^{\vecv,u}_s|^2\big)^{pr} ds\big|\mcF_t\Big]&=\E^{\bbQ^{\vecv,u}}\Big[\mcE(-\gamma^{\vecv,u}*W^{\vecv,u})_{t,T}\big(\int_t^T|V^{\vecv,u}_s|^2 ds\big)^{pr}\Big|\mcF_t\Big]
\\
&\leq C(\tilde R^{\vecv}_t)^{1/r'}\E^{\bbQ^{\vecv,u}}\Big[|\xi^{\vecv\circ u}|^{2pr^2}
\\
&\quad+\big(\int_t^T|f^{\vecv\circ u}(s,0,0)|ds\big)^{2pr^2}+\sum_{j=1}^N |c^{\vecv\circ[u]_{j-1}}(\tau_j,\beta_j)|^{2pr^2}\Big|\mcF_t\Big]^{1/r}
\\
&\leq C(\tilde R^{\vecv}_t)^{1/r'}(R^{\vecv}_t)^{1/r'}(\bar K^{\vecv,2pr^3}_t)^{1/r^2}.
\end{align*}
%where
%\begin{align*}
%R_t:=\E^{\bbQ^{\vecv,u}}\big[|\mcE((\gamma-\gamma^{\vecv,u})*W^{\vecv,u})_{t,T}|^{r'}\big|\mcF_t\big].
%\end{align*}
Combined, this gives that
\begin{align*}
|\delta U^u_t|^p&\leq C(R^{\vecv'}_t)^{1/r'}((\bar K^{\vecv,\vecv',l+1,pr}_t)^{1/r}+(\bar K^{\vecv,\vecv',l+1,2pr}_t)^{1/2r}(1+(\tilde R^{\vecv}_t)^{1/2rr'}(R^{\vecv}_t)^{1/2rr'}(\bar K^{\vecv,2pr^3}_t)^{1/2r^3})).
\end{align*}
Hence, repeated application of H\"older's inequality gives, with $\delta Y_t:=\sup_{b\in U}|Y^{\vecv'\circ(t,b),l}_t-Y^{\vecv\circ(t,b),l}_t|$, that
\begin{align*}
\|\delta Y\|_{\mcS^{p}}^{p}&\leq C(\|R^{\vecv'}\|_{\mcS^1}^{1/r'}(\|\bar K^{\vecv,\vecv',l+1,pr}\|_{\mcS^1}^{1/r}+\|(\bar K^{\vecv,\vecv',l+1,2pr})^{1/2}(1+(\tilde R^{\vecv})^{1/2r'}(R^{\vecv})^{1/2r'}(\bar K^{\vecv,2pr^3})^{1/2r^2})\|_{\mcS^1}^{1/r})
\\
&\leq C(\|R^{\vecv'}\|_{\mcS^1}^{1/r'}(\|\bar K^{\vecv,\vecv',l+1,pr}\|_{\mcS^1}^{1/r}+\|\bar K^{\vecv,\vecv',l+1,2pr}\|_{\mcS^1}^{1/2}(1+\|(\tilde R^{\vecv})^{1/r'}(R^{\vecv})^{1/r'}(\bar K^{\vecv,2pr^3})^{1/r^2})\|_{\mcS^1}^{1/2r})
\\
&\leq C(\|R^{\vecv'}\|_{\mcS^1}^{1/r'}(\|\bar K^{\vecv,\vecv',l+1,pr}\|_{\mcS^1}^{1/r}+\|\bar K^{\vecv,\vecv',l+1,2pr}\|_{\mcS^1}^{1/2}(1+\|\tilde R^{\vecv}R^{\vecv}\|_{\mcS^1}^{1/2rr'}\|\bar K^{\vecv,2pr^3}\|_{\mcS^1}^{1/2r^3})
\end{align*}
Hence, as $\|\tilde R^{\vecv}R^{\vecv}\|_{\mcS^1}\leq \|\tilde R^{\vecv}\|_{\mcS^2}\|R^{\vecv}\|_{\mcS^2}\leq C$ it follows that
\begin{align*}
\|\delta Y\|_{\mcS^{p'}}^{p'}&\leq C(\|\bar K^{\vecv,\vecv',l+1,p'r}\|_{\mcS^1}^{1/r}+\|\bar K^{\vecv,\vecv',l+1,2p'r}\|_{\mcS^1}^{1/2r})
\end{align*}
%We thus find that
%\begin{align*}
%\|(\sup_{b\in U}|Y^{\vecv'\circ(t,b),l}_t-Y^{\vecv\circ(t,b),l}_t|:t\in[0,T])\|^p_{\mcS^p}&\leq C (\|K^{\vecv,\vecv'}\|^p_{\mcS^p}+\|K^{\vecv,\vecv'}K^{\vecv}\|^{p/2}_{\mcS^p})
%\\
%&\leq C (\|K^{\vecv,\vecv'}\|^p_{\mcS^p}+\|K^{\vecv,\vecv'}\|^{p/2}_{\mcS^{2p}}\|K^{\vecv}\|^{p/2}_{\mcS^{2p}})
%\end{align*}
%and the result follows by Assumption~\ref{ass:on-coeff-FH}.
and the result follows by Definition~\ref{defn:bounding-fam}.iv).\qed\\

%%%%%%%%%%%%%%%%%%%%%%%%%%%%%%%%%%%%%%%%%%%%%%%%%%%%%%%%%%%%%%%%%%%%%%%%%%%%%%%%%%%%%%%%%%%%%%%%%%%%%%%%%%%%%%%%%%%%%%%%%%%%%%%%%%%%%%%%%%%%%%%%%%

\begin{prop}\label{prop:seq-rbsde_k-FH}
Hypothesis {\bf RBSDE.$l$} holds for all $l\geq 0$.
\end{prop}

\noindent\emph{Proof.} We note that the triple $(Y^{v,l+1},Z^{v,l+1},K^{v,l+1})$ solves a reflected BSDE with barrier $(\sup_{b\in U}\{-c^{v}(t,b)+Y^{v\circ(t,b),l}_t\}:t\in [0,T])$. By Lemma~\ref{lem:seq-diag-bound-FH} we find that $((\sup_{b\in U}\{-c^{v}(t,b)+Y^{v\circ(t,b),l}_t\})^+:t\in [0,T])\in\mcS^p$ for all $p\geq 1$ uniformly in $v$. Whenever the statement in Hypothesis {\bf RBSDE.$l$} holds for some $l\geq 0$, then Proposition~\ref{prop:rbsde-solu-FH} guarantees the existence of a unique triple $(Y^{v,l+1},Z^{v,l+1},K^{v,l+1})$ solving \eqref{ekv:rbsde_k} with $k=l+1$. Moreover, we have
\begin{align*}
&\E\Big[\sup_{t\in[0,T]}|\sup_{b\in U}\{-c^{v}(t,b)+Y^{v\circ(t,b),l}_{t}\}-\sup_{b\in U}\{-c^{v'}(t,b)+Y^{v'\circ(t,b),l}_t\}|^p\Big]
\\
&\leq C\E\Big[\sup_{t\in[0,T]}\sup_{b\in U}|c^{v}(t,b)-c^{v'}(t,b)|^p+\sup_{t\in[0,T]}\sup_{b\in U}|Y^{v\circ(t,b),l}_{t}-Y^{v'\circ(t,b),l}_t|^p\Big].
%\\
%&\leq \|v-\tilde v\|_{\mcU^f,C}^{p/2}.
\end{align*}
By \eqref{ekv:YZK-diff-FH} of Proposition~\ref{prop:rbsde-solu-FH}, Assumption~\ref{ass:on-coeff-FH} and Lemma~\ref{lem:seq-diff-FH} it, thus, follows by repeating the argument in the proof of Proposition~\ref{prop:seq-bsde_0-FH} that for each $p\geq 1$, there is a $p'\geq 1$ such that
\begin{align*}
&\|Y^{\vecv',l+1}-Y^{\vecv,l+1}\|_{\mcS^{p'}}^{p'}+\|Z^{\vecv',l+1}-Z^{\vecv,l+1}\|_{\mcH^{p'}}^{p'} + \|K^{\vecv',l+1}-K^{\vecv,l+1}\|_{\mcS^{p'}}^{p'}\leq C\|\vecv'-\vecv\|_{\mcD^f}^p.
\end{align*}
For $\vecv',\vecv\in\mcD^{\kappa}$ we let $p=(m+1){\kappa}+1$ and Kolmogorov's continuity theorem implies the existence of a family of processes $(\hat Y^{\vecv},\hat Z^{\vecv},\hat K^{\vecv}:\vecv\in\mcD^{\kappa})$, with $(\hat Y^{\vecv},\hat Z^{\vecv},\hat K^{\vecv})\in\mcS^2\times\mcH^2\times\mcS^2$, such that (outside of a $\Prob$-null set) $\vecv\mapsto\hat Y^{\vecv}_\cdot$, $\vecv\mapsto\hat Z^{\vecv}_\cdot$ and $\vecv\mapsto\hat K^{\vecv}_\cdot$ are uniformly continuous, $L^2([0,T])$ continuous (and that $\vecv\mapsto\int_s^TZ^{\vecv}_rdW_r$ is continuous in $\vecv$ uniformly in $s$) and uniformly continuous, respectively, and moreover $(\hat Y^{\vecv}_s,\hat Z^{\vecv}_s,\hat K^{\vecv}_s)=(Y^{\vecv,l+1}_s,Z^{\vecv,l+1}_s,K^{\vecv,l+1}_s)$, $\Prob$-a.s.~for all $\vecv\in\mcD^{\kappa}$ and $s\in[0,T]$.

Now, taking countable unions this extends to $\mcD^f$, and for all $\vecv\in\mcD^f$, the triple $(\hat Y^{\vecv},\hat Z^{\vecv},\hat K^{\vecv})$ solves the BSDE \eqref{ekv:rbsde_k} for $Y^{\vecv,l+1}$.% and that the map $\vecv\mapsto (\hat Y^{\vecv},\hat Z^{\vecv},\hat K^{\vecv}):\mcD^f\to\mcS^2\times\mcH^2\times\mcS^2$ is continuous.

Furthermore, for any $v\in\mcU^f$ and any approximating sequence $(v_j)_{j\geq 0}$ taking values in a countable dense subset of $\mcD^f$ with $v_j\to v$, $\Prob$-a.s., and $v_j\in \mcU^f$, we have
\begin{align*}
\sup_{s\in [0,T]}|\hat Y^{v_j}_s-\hat Y^{v}_s|\to 0,
\end{align*}
$\Prob$-a.s., as $j\to\infty$. Moreover, due to continuity of the map $\vecv\mapsto \xi^{\vecv}$ we have
\begin{align*}
  \xi^{v_j}\to \xi^{v}
\end{align*}
and by repeating the argument in the proof of Proposition~\ref{prop:seq-bsde_0-FH} we get
\begin{align*}
\int_0^T |f^{v_j}(s,\hat Y^{v_j}_s,\hat Z^{v_j}_s)-f^{v}(s,\hat Y^{v}_s,\hat Z^{v}_s)|ds\to 0,
\end{align*}
$\Prob$-a.s., as $j\to\infty$. We conclude that
\begin{align*}
\hat Y^{v}_\eta&=\lim_{j\to\infty}\Big\{\xi^{v_j}+\int_\eta^T f^{v_j}(s,\hat Y^{v_j}_s,\hat Z^{v_j}_s)ds-\int_\eta^T \hat Z^{v_j}_sdW_s+\hat K^{v_j}_T-\hat K^{v_j}_\eta\Big\}
\\
&=\xi^{v}+\int_\eta^T f^{v}(s,\hat Y^{v}_s,\hat Z^{v}_s)ds-\int_\eta^T \hat Z^{v}_sdW_s+\hat K^{v}_T- \hat K^{v}_\eta,
\end{align*}
for all $\eta\in\mcT$. Finally, by Helly's convergence theorem (see \eg \cite{KolmogorovFomin}, p. 370) we have
\begin{align*}
&\lim_{j\to\infty}\int_0^T(\hat Y^{v}_s-\sup_{b\in U}\{-c^{v}(s,b)+\hat Y^{v\circ(s,b),l}_s\})d\hat K^{v_j}_s
\\
&=\int_0^T(\hat Y^{v}_s-\sup_{b\in U}\{-c^{v}(s,b)+\hat Y^{v\circ(s,b),l}_s\})d\hat K^{v}_s
\end{align*}
and since
\begin{align*}
&\int_0^T|(\hat Y^{v}_s-\sup_{b\in U}\{-c^{v}(s,b)+\hat Y^{v\circ(s,b),l}_s\})-(\hat Y^{v_j}_s-\sup_{b\in U}\{-c^{v_j}(s,b)+\hat Y^{v_j\circ(s,b),l}_s\})|d\hat K^{v_j}_s
\\
&\leq \hat K^{v_j}_T(\sup_{s\in[0,T]\times U}|\hat Y^{v}_s-\hat Y^{v_j}_s|+\sup_{(s,b)\in[0,T]\times U}|c^{v}(s,b)-c^{v_j}(s,b)|
\\
&\quad+\sup_{(s,b)\in[0,T]\times U}|\hat Y^{v\circ(s,b),l}_s-\hat Y^{v_j\circ(s,b),l}_s|)
\end{align*}
which tends to zero, $\Prob$-a.s., as $j\to\infty$ we get that
\begin{align*}
\int_0^T(\hat Y^{v}_s-\sup_{b\in U}\{-c^{v}(s,b)+\hat Y^{v\circ(s,b),l}_s\})d\hat K^{v}_s
&=\lim_{j\to\infty}\int_0^T(\hat Y^{v_j}_s-\sup_{b\in U}\{-c^{v_j}(s,b)+\hat Y^{v_j\circ(s,b),l}_s\})d\hat K^{v_j}_s
\\
&=0
\end{align*}
and we conclude that Hypothesis {\bf RBSDE.$l+1$} holds as well. The statement of the proposition now follows by an induction argument.\qed\\

%%%%%%%%%%%%%%%%%%%%%%%%%%%%%%%%%%%%%%%%%%%%%%%%%%%%%%%%%%%%%%%%%%%%%%%%%%%%%%%%%%%%%%%%%%%%%%%%%%%%%%%%%%%%%%%%%%%%%%%%%%%%%%%%%%%%%%%%%%%%%%%%%%

\subsection{Convergence of the scheme}

We now show that there exists a limit family of triples $(\bar Y^v,\bar Z^v,\bar K^v:v\in\mcU^f):=\lim_{k\to\infty}(Y^{v,k},Z^{v,k},K^{v,k}:v\in\mcU^f)$ that solves the sequential system of reflected BSDEs~\eqref{ekv:seq-bsde-FH}. This result relies heavily upon the following two lemmas and their corollaries.

\begin{lem}\label{lem:seq-U-V-bound}
For $v\in\mcU^f$ and $k\geq 0$, assume that $u^*\in\mcU^k_t$ is such that $U_t^{v,u^*}=\esssup_{u\in\mcU^k_t}U_t^{v,u}$. Then, for each $p\geq 1$, there is a $C>0$, that does not depend on $v$ or $k$, such that
\begin{align*}
\E\Big[\sup_{s\in[t,T]}|U_s^{v,u^*}|^p+\big(\int_t^T|V_s^{v,u^*}|^2ds\big)^{p/2}\Big|\mcF_t\Big]\leq C(\tilde R^v_t R^v_t\bar R^v_t)^{1/r'}(1+\bar K^{v,2pr^3}_t)^{1/r^3}.
\end{align*}
\end{lem}

\noindent\emph{Proof.} For the bound on $U^{v,u^*}$ we note that by \eqref{ekv:U-vu-emptyset-bound} we have
\begin{align*}
|U^{v,u^*}_s|^p &= |Y^{v\circ(u^*(s-)),k-N^*(s-)}_s|^p
\\
&\leq C (R^{v\circ u^*(s-)}_s)^{1/r'}(\bar K^{v\circ u^*(s-),pr}_s)^{1/r}
\end{align*}
and since by definition we have $R^{v\circ (u^*(s))}_s=R^{v\circ u^*}_s$, we find that
\begin{align}\nonumber
\E\Big[\sup_{s\in[t,T]}|U^{v,u^*}_s|^{p}\Big|\mcF_t\Big]&\leq\E\Big[\sup_{s\in[t,T]}(R^{v\circ u^*}_s)^{1/r'}(\bar K^{v\circ (u^*(s)),pr}_s)^{1/r}\Big|\mcF_t\Big]
\\
&\leq \esssup_{u\in\mcU^f_t}\E\Big[\sup_{s\in[t,T]}R^{v\circ u}_s\Big|\mcF_t\Big]^{1/r'}\E\Big[\sup_{s\in[t,T]}\bar K^{v\circ (u^*(s)),pr}_s\Big|\mcF_t\Big]^{1/r}\nonumber
\\
&\leq C(\bar R^v_t)^{1/r'}(1+\bar K^{v,2pr}_t)^{1/2r}\label{ekv:U-v-u-k-bound}
\end{align}
where we have used Jensen's inequality and \eqref{ekv:barK-Doob-equivalent} to reach the last inequality. The first bound then follows by Jensen's inequality since $\tilde R^v_t R^v_t\geq 1$, $\Prob$-a.s.

We apply Ito's formula to $(U^{v,u^*})^2$ and get that
\begin{align*}
|U^{v,u^*}_t|^2+\int_t^T|V^{v,u^*}_s|^2ds&=|\xi^{v\circ u^*}|^2+2\int_t^{T}U^{v,u^*}_sf^{v\circ u^*}(s, U^{v,u^*}_s, V^{v,u^*}_s)ds
-2\int_t^{T} U^{v,u^*}_sV^{v,u^*}_sdW_s
\\
&\quad+\sum_{j=1}^{N^*} (-2U^{v,[u^*]_{j-1}}_{\tau^*_j}c^{v\circ [u^*]_{j-1}}(\tau^*_j,\beta^*_j)+|c^{v\circ [u^*]_{j-1}}(\tau^*_j,\beta^*_j)|^2)
\\
&\leq |\xi^{v\circ u^*}|^2+2\int_t^{T}(\gamma^{v,u^*}_s|U^{v,u^*}_s|^2+U^{v,u^*}_s f^{v\circ u^*}(s,0, 0)ds
\\
&\quad-2\int_t^{T} U^{v,u^*}_sV^{v,u^*}_sdW^{v,u^*}_s+4\sup_{s\in[t,T]}|U^{v,u^*}_s|\sum_{j=1}^{N^*}|c^{v\circ [u^*]_{j-1}}(\tau^*_j,\beta^*_j)|,
\end{align*}
where the last term appears after applying the relation $|c^{v\circ [u^*]_{j-1}}(\tau^*_j,\beta^*_j)|\leq 2\sup_{s\in[t,T]}|U^{v,u^*}_s|$. Using the relation $ab\leq \frac{1}{2}(\kappa a^2+\frac{1}{\kappa}b^2)$ for $\kappa>0$ we get
\begin{align*}
|U^{v,u^*}_t|^2+\int_t^T|V^{v,u^*}_s|^2ds&\leq |\xi^{v\circ u^*}|^2+(C+2\kappa)\sup_{s\in[t,T]}|U^{v,u^*}_s|^2+\int_t^{T}|f^{v\circ u^*}(s,0, 0)|^2ds
\\
&\quad-2\int_t^{T} U^{v,u^*}_sV^{v,u^*}_sdW^{v,u^*}_s+\frac{2}{\kappa}\Big(\sum_{j=1}^{N^*} c^{v\circ [u^*]_{j-1}}(\tau^*_j,\beta^*_j)\Big)^2.
\end{align*}
On the other hand, applying the usual manipulations to \eqref{ekv:U-def} we get
\begin{align*}
 U_t^{v,u^*}&=e^{v,u^*}_{t,T}\xi^{v\circ u^*}+\int_t^{T}e^{v,u^*}_{t,s}f^{v\circ u^*}(s, 0, 0)ds
-\int_t^{T} e^{v,u^*}_{t,s}V^{v,u^*}_sdW^{v,u^*}_s-\sum_{j=1}^{N^*} e^{v,u^*}_{t,\tau_j}c^{v\circ [u^*]_{j-1}}(\tau^*_j,\beta^*_j).
\end{align*}
Rearranging terms now gives us (with $e_{t,\cdot}=e_{t,\cdot}^{v,u^*}$)
\begin{align}%\nonumber
\sum_{j=1}^N e_{t,\tau_j^*}c^{v\circ [u^*]_{j-1}}(\tau_j^*,\beta_j^*) &= e_{t,T}\xi^{v\circ u^*}+\int_t^{T}e_{t,s}f^{v\circ u^*}(s, 0, 0)ds-\int_t^{T} e_{t,s}V^{v,u^*}_sdW^{v,u^*}_s-U_t^{v,u^*}
%\\
%&\leq e_{t,T}\xi^{v\circ u^*}+\int_t^{T}e_{t,s}f^{v\circ u^*}(s, 0, 0)ds-\int_t^{T} e_{t,s}V^{v,u^*}_sdW^{v,u^*}_s-Y_t^{v,0}.
\label{ekv:interv-cost-bound}
\end{align}
From \eqref{ekv:interv-cost-bound} we have that
\begin{align*}
\Big(\sum_{j=1}^{N^*} c^{v\circ [u^*]_{j-1}}(\tau^*_j,\beta^*_j)\Big)^2&\leq C(|\xi^{v\circ u^*}|^2+\int_t^{T}|f^{v\circ u^*}(s, 0, 0)|^2ds
\\
&\quad+\Big|\int_t^{T} e_{t,s}V^{v,u^*}_sdW^{v,u^*}_s\Big|^2+|U_t^{v,u^*}|^2).
\end{align*}
Put together this gives
\begin{align*}
|U^{v,u^*}_t|^2+\int_t^T|V^{v,u^*}_s|^2ds&\leq C(1+\kappa+\frac{1}{\kappa})(|\xi^{v\circ u^*}|^2+\int_t^{T}|f^{v\circ u^*}(s,0, 0)|^2ds+\sup_{s\in[t,T]}|U^{v,u^*}_s|^2)
\\
&\quad-2\int_t^{T} U^{v,u^*}_sV^{v,u^*}_sdW^{v,u^*}_s+\frac{C}{\kappa}\Big|\int_t^{T} e_{t,s}V^{v,u^*}_sdW^{v,u^*}_s\Big|^2.
\end{align*}
Raising both sides to $p/2$ and taking the conditional expectation we find that
\begin{align*}
&\E^{\bbQ^{v,u^*}}\Big[\big(\int_t^T|V^{v,u^*}_s|^2ds\big)^{p/2}\Big|\mcF_t\Big]\leq C\E^{\bbQ^{v,u^*}}\Big[(1+\kappa^{p/2}+\frac{1}{\kappa^{p/2}})(|\xi^{v\circ u^*}|^p+\int_t^{T}|f^{v\circ u^*}(s,0, 0)|^pds
\\
&\quad+\sup_{s\in[t,T]}|U^{v,u^*}_s|^p)+\big(\int_t^{T} |U^{v,u^*}_sV^{v,u^*}_s|^2ds\big)^{p/4}+\frac{1}{\kappa^{p/2}}\big(\int_t^{T} |V^{v,u^*}_s|^2ds\big)^{p/2}\Big|\mcF_t\Big]
\end{align*}
and since
\begin{align*}
\E^{\bbQ^{v,u^*}}\Big[\big(\int_t^{T} |U^{v,u^*}_sV^{v,u^*}_s|^2ds\big)^{p/4}\Big|\mcF_t\Big]\leq \frac{1}{2}\E^{\bbQ^{v,u^*}}\Big[\kappa\sup_{s\in[t,T]}|U^{v,u^*}_s|^p+\frac{1}{\kappa}\big(\int_t^{T} |V^{v,u^*}_s|^2ds\big)^{p/2}\Big|\mcF_t\Big]
\end{align*}
we arrive at the inequality
\begin{align}
\E^{\bbQ^{v,u^*}}\Big[\big(\int_t^T|V^{v,u^*}_s|^2ds\big)^{p/2}\Big|\mcF_t\Big]&\leq C\E^{\bbQ^{v,u^*}}\Big[|\xi^{v\circ u^*}|^p+\int_t^{T}|f^{v\circ u^*}(s,0, 0)|^pds+\sup_{s\in[t,T]}|U^{v,u^*}_s|^p\Big|\mcF_t\Big]\label{ekv:U-V_t-bound-bbQ}
\end{align}
by choosing $\kappa>0$ sufficiently large. Under $\Prob$ this rewrites as
\begin{align}\nonumber
\E\Big[\big(\int_t^T|V^{v,u^*}_s|^2ds\big)^{p/2}\Big|\mcF_t\Big]&\leq C(\tilde R^v_t R^v_t)^{1/r'}\E\Big[|\xi^{v\circ u^*}|^{pr^2}
\\
&\quad+\int_t^{T}|f^{v\circ u^*}(s,0, 0)|^{pr^2}ds+\sup_{s\in[t,T]}|U^{v,u^*}_s|^{pr^2}\Big|\mcF_t\Big]^{1/r^2}\label{ekv:U-V_t-bound}
\end{align}
The desired result now follows by setting $p\leftarrow pr^2$ in \eqref{ekv:U-v-u-k-bound} and using Jensen's inequality while noting that $\bar R^v_t\geq 1$, $\Prob$-a.s.\qed\\

\begin{cor}\label{cor:seq-U-V-bound}
For $v\in\mcU^f$, $t\in[0,T]$, $\beta\in \mcI(t)$ and $k\geq 0$, assume that $u^*\in\mcU^k_t$ is such that $U_t^{v\circ(t,\beta),u^*}=\esssup_{u\in\mcU^k_t}U_t^{v\circ(t,\beta),u}$. Then, for each $p\geq 1$, there is a $C>0$ that does not depend on $v$ or $k$, such that
\begin{align*}
\E\Big[\sup_{s\in[t,T]}|U_s^{v\circ(t,\beta),u^*}|^p+\big(\int_t^T|V_s^{v\circ(t,\beta),u^*}|^2ds\big)^{p/2}\Big|\mcF_t\Big]\leq C(\tilde R^v_t R^v_t\bar R^v_t)^{1/r'}(1+\bar K^{v,2pr^3}_t)^{1/r^3}.
\end{align*}
\end{cor}

\noindent\emph{Proof.} This follows immediately by making suitable manipulations, \ie setting $v\leftarrow v\circ(t,\beta)$, in the proof of Lemma~\ref{lem:seq-U-V-bound}.\qed\\

\begin{lem}\label{lem:N-bound}
For $v\in\mcU^f$ and $k\geq 0$, assume that $u^*\in\mcU^k_t$ is such that $U_t^{v,u^*}=\esssup_{u\in\mcU^k_t}U_t^{v,u}$. Then, for each $p\geq 1$, there is a $C>0$ that does not depend on $v$ or $k$, such that
\begin{align*}
\E\big[(N^*)^p\big|\mcF_t\big]\leq C\E\Big[\Big(\sum_{j=1}^N c^{v\circ [u^*]_{j-1}}(\tau_j^*,\beta_j^*)\Big)^p\Big|\mcF_t\Big]\leq C(\tilde R^v_t R^v_t\bar R^v_t)^{1/r'}(1+\bar K^{v,2pr^3}_t)^{1/r^3}.
\end{align*}
\end{lem}

\noindent\emph{Proof.} Since the intervention costs are bounded from below by $\delta>0$, we have
\begin{align*}
\delta N^*&\leq\sum_{j=1}^N c^{v\circ [u^*]_{j-1}}(\tau_j^*,\beta_j^*)
\end{align*}
from which the first inequality follows. Now, from \eqref{ekv:interv-cost-bound} we have
\begin{align*}
&\E^{\bbQ^{v, u^*}}\Big[\Big(\sum_{j=1}^N c^{v\circ [u^*]_{j-1}}(\tau_j^*,\beta_j^*)\Big)^p\Big|\mcF_t\Big]
\\
&\leq C\E^{\bbQ^{v,u^*}}\Big[|\xi^{v\circ u^*}|^p+\int_t^{T}|f^{v\circ u^*}(s,0, 0)|^pds+\big(\int_t^{T}|V^{v, u^*}_s|^2ds\big)^{p/2}+|U^{v,u^*}_t|^p\Big|\mcF_t\Big]
\\
&\leq C\E^{\bbQ^{v,u^*}}\Big[|\xi^{v\circ u^*}|^p+\int_t^{T}|f^{v\circ u^*}(s,0, 0)|^pds+\sup_{s\in[t,T]}|U^{v,u^*}_s|^p\Big|\mcF_t\Big]
\end{align*}
where we have used \eqref{ekv:U-V_t-bound-bbQ} to get the last inequality. Now the result is immediate from the last part in the proof of Lemma~\ref{lem:seq-U-V-bound}.\qed\\

\begin{cor}\label{cor:N-bound}
For $v\in\mcU^f$, $t\in[0,T]$, $\beta\in \mcI(t)$ and $k\geq 0$, assume that $u^*\in\mcU^k_t$ is such that $U_t^{v\circ(t,\beta),u^*}=\esssup_{u\in\mcU^k_t}U_t^{v\circ(t,\beta),u}$. Then, for each $p\geq 1$, there is a $C>0$ that does not depend on $v$ or $k$, such that
\begin{align*}
\E\big[(N^*)^p\big|\mcF_t\big]\leq C\E\Big[\Big(\sum_{j=1}^N c^{v\circ(t,\beta)\circ [u^*]_{j-1}}(\tau_j^*,\beta_j^*)\Big)^p\Big|\mcF_t\Big]\leq C(\tilde R^v_t R^v_t\bar R^v_t)^{1/r'}(1+\bar K^{v,2pr^3}_t)^{1/r^3}.
\end{align*}
\end{cor}

\noindent\emph{Proof.} This follows by repeating the argument in the proof of Lemma~\ref{lem:N-bound} after making the swap $v\leftarrow v\circ(t,\beta)$.\qed\\

%%%%%%%%%%%%%%%%%%%%%%%%%%%%%%%%%%%%%%%%%%%%%%%%%%%%%%%%%%%%%%%%%%%%%%%%%%%%%%%%%%%%%%%%%%%%%%%%%%%%%%%%%%%%%%%%%%%%%%%%%%%%%%%%%%%%%%%%%%%%%%%%%%

We are now ready to tackle the convergence of the sequence $Y^{v,k}$, this is done in the following proposition, where i) and iii) are the important properties and ii) is included only because it is notationally simpler to verify than iii).

\begin{prop}\label{prop:seq-rbsde_k-conv}
There exists a limit family $(\bar Y^v:v\in\mcU^f)$ such that for all $v\in\mcU^f$ (outside of a $\Prob$-null set) and all $p\geq 2$, we have
\begin{enumerate}[i)]
  \item $Y^{v,k}\nearrow\bar Y^v$ pointwisely,
  \item $\|\bar Y^v-Y^{v,k}\|_{\mcS^2}\to 0$ as $k\to\infty$, and
  \item $\|\sup_{b\in U}|\bar Y^{v\circ(\cdot,b)}_\cdot-Y^{v\circ(\cdot,b),k}_\cdot|\|_{\mcS^p}\to 0$ as $k\to\infty$.
\end{enumerate}
\end{prop}

\noindent\emph{Proof.} The sequence $(Y^{v,k})_{k\geq 0}$ is non-decreasing and $\Prob$-a.s.~bounded by Lemma~\ref{lem:seq-unif-bound-FH}. Thus it converges pointwisely, $\Prob$-a.s., and \emph{i)} follows.

We now turn our focus to the second claim and note by Lemma~\ref{lem:N-bound} that if $u^*=(\tau_1^*,\ldots,\tau^*_{N^*};\beta_1^*,\ldots,\beta^*_{N^*})\in\mcU^k_t$ is such that $U^{v,u^*}_t=Y^{v,k}_t$, then $\E\big[N^*\big|\mcF_t\big]\leq C(\tilde R^v_t R^v_t\bar R^v_t)^{1/r'}(1+\bar K^{v,2r^3}_t)^{1/r^3}$ and, in particular, we find that $\E\big[\ett_{[N^*> k']}\big|\mcF_t\big]\leq C(\tilde R^v_t R^v_t\bar R^v_t)^{1/r'}(1+\bar K^{v,2r^3}_t)^{1/r^3}/k'$ for all $k'\geq 1$.

For any $k'$ with $0\leq k'\leq k$, the truncation $[u^*]_{k'}$ belongs to $\mcU^{k'}_t$. We, thus, have
\begin{align*}
U^{v,[u^*]_{k'}}_t\leq Y^{v,k'}_t\leq Y^{v,k}_t.
\end{align*}
Since $Y^{v,k}_t=U^{v,u^*}_t$, this gives
\begin{align*}
|Y^{v,k}_t-Y^{v,k'}_t|&\leq U^{v,u^*}_t-U^{v,[u^*]_{k'}}_t
\end{align*}
Moreover, since the intervention costs are positive, we have that% since the truncation only has affect when $N^*>k'$ and
\begin{align*}
U^{v,u^*}_t-U^{v,[u^*]_{k'}}_t&\leq \xi^{v\circ u^*}-\xi^{v\circ [u^*]_{k'}}+\int_t^T(f^{v\circ u^*}(s,U^{v,u^*}_s,V^{v,u^*}_s)
\\
&\quad-f^{v\circ [u^*]_{k'}}(s,U^{v,[u^*]_{k'}}_s,V^{v,[u^*]_{k'}}_s))ds-\int_t^T(V^{v,u^*}_s-V^{v,[u^*]_{k'}}_s)dW_s.
\end{align*}
Setting $e_{t,s}:=e^{\int_t^s\gamma_s}$ with
\begin{align*}
\gamma_s:=\frac{f^{v\circ [u^*]_{k'}}(s,U^{v,u^*}_s,V^{v,[u^*]_{k'}}_s)-f^{v\circ [u^*]_{k'}}(s,U^{v,[u^*]_{k'}}_s,V^{v,[u^*]_{k'}}_s)}{U^{v,u^*}_s-U^{v,[u^*]_{k'}}_s}\ett_{[U^{v,u^*}_s\neq U^{v,[u^*]_{k'}}_s]}
\end{align*}
and arguing as in the proof of Lemma~\ref{lem:rbsde-trunk-diff} gives
\begin{align*}
U^{v,u^*}_t-U^{v,[u^*]_{k'}}_t&\leq e_{t,T}(\xi^{v\circ u^*}-\xi^{v\circ [u^*]_{k'}})+\int_t^Te_{t,s}(f^{v\circ u^*}(s,U^{v,u^*}_s,V^{v,u^*}_s)
\\
&\quad-f^{v\circ [u^*]_{k'}}(s,U^{v,u^*}_s,V^{v,[u^*]_{k'}}_s))ds-\int_t^Te_{t,s}(V^{v,u^*}_s-V^{v,[u^*]_{k'}}_s)dW_s
\\
&=e_{t,T}(\xi^{v\circ u^*}-\xi^{v\circ [u^*]_{k'}})+\int_t^Te_{t,s}(f^{v\circ u^*}-f^{v\circ [u^*]_{k'}})(s,U^{v,u^*}_s,V^{v,u^*}_s)ds
\\
&\quad-\int_t^Te_{t,s}(V^{v,u^*}_s-V^{v,[u^*]_{k'}}_s)dW^{\zeta}_s,
\end{align*}
where $W^{\zeta}_t:=W_t-\int_0^t\zeta_sds$, with
\begin{align*}
\zeta_s:=\frac{f^{v\circ [u^*]_{k'}}(s,U^{v,u^*}_s,V^{v,u^*}_s)-f^{v\circ [u^*]_{k'}}(s,U^{v,u^*}_s,V^{v,[u^*]_{k'}}_s)}{|V^{v,u^*}_s-V^{v,[u^*]_{k'}}_s|^2}(V^{v,u^*}_s-V^{v,[u^*]_{k'}}_s)^\top \ett_{[V^{v,u^*}_s\neq V^{v,[u^*]_{k'}}_s]}.
\end{align*}
Taking the conditional expectation, using that $V^{v,u^*},V^{v,[u^*]_{k'}}\in\mcH^2_{\bbQ^{\zeta}}$ by Proposition~\ref{prop:UV-has-solution} and noting that the right-hand side is non-zero only when $N^*>k'$ gives
\begin{align*}
U^{v,u^*}_t-U^{v,[u^*]_{k'}}_t&\leq C\E^{\bbQ^{\zeta}}\Big[\ett_{[N^*>k']}\big(|\xi^{v\circ u^*}-\xi^{v\circ [u^*]_{k'}}|+\int_t^T|(f^{v\circ u^*}-f^{v\circ [u^*]_{k'}})(s,U^{v,u^*}_s,V^{v,u^*}_s)|ds\big)\Big|\mcF_t\Big].
\end{align*}
By H\"older's inequality we find that
\begin{align*}
|Y^{v,k}_t-Y^{v,k'}_t|^p&\leq C \E^{\bbQ^{\zeta}}\big[\ett_{[N^*>k']}\big|\mcF_t\big]^{1/2}\E^{\bbQ^{\zeta}}\Big[|\xi^{v\circ u^*}-\xi^{v\circ [u^*]_{k'}}|^{2p}+\big(\int_t^T|(f^{v\circ u^*}-f^{v\circ [u^*]_{k'}})(s,0,0)|ds\big)^{2p}
\\
&\quad+\sup_{s\in[t,T]}|U^{v,u^*}_s|^{2p}+|L^{v\circ u^*}_T|^{4p}+|L^{v\circ [u^*]_{k'}}_T|^{4p}+\big(\int_t^T|V^{v,u^*}_s|^{2}ds\big)^{2p}\Big|\mcF_t\Big]^{1/2}
\\
&\leq C (R^v_t)^{1/r'}\E\big[\ett_{[N^*>k']}\big|\mcF_t\big]^{1/2r}\E\Big[|\xi^{v\circ u^*}-\xi^{v\circ [u^*]_{k'}}|^{2pr}
\\
&\quad+\big(\int_t^T|(f^{v\circ u^*}-f^{v\circ [u^*]_{k'}})(s,0,0)|ds\big)^{2pr}
\\
&\quad+\sup_{s\in[t,T]}|U^{v,u^*}_s|^{2pr}+|L^{v\circ u^*}_T|^{4pr}+|L^{v\circ [u^*]_{k'}}_T|^{4pr}+\big(\int_t^T|V^{v,u^*}_s|^2 ds\big)^{2pr}\Big|\mcF_t\Big]^{1/2r}
\\
&\leq C(R^v_t)^{1/r'}(\tilde R^v_t R^v_t\bar R^v_t)^{1/rr'}(1+\bar K^{v,2r^3}_t)^{1/2r^4}(1+\bar K^{v,8pr^4}_t)^{1/2r^4}/(k')^{1/2r},
\\
&\leq C(R^v_t)^{1/r'}(\tilde R^v_t R^v_t\bar R^v_t)^{1/rr'}(1+\bar K^{v,8pr^4}_t)^{1/r^4}/(k')^{1/2r},
\end{align*}
where we have used \eqref{ekv:L-bound} and Lemma~\ref{lem:seq-U-V-bound} to arrive at the last inequality.
%Now, since $\bbQ^{\zeta}\in\PrM^v_t$ we have that $\E^{\bbQ^{\zeta}}\big[\ett_{[N^*>k']}\big|\mcF_t\big]\leq C(R^v_t)^{1/r'}(\bar K^{v,r}_t)^{1/r}/k'$ which by Assumption~\ref{ass:on-coeff-FH}.(\emph{i}) gives
%\begin{align*}
%|Y^{v,k}_t-Y^{v,k'}_t|^2&\leq C\frac{(R^v_t)^{1/r'}(\bar K^{v,r}_t)^{1/r}}{k'}(K^{v,2}_t+K^{v,4}_t).
%\end{align*}
Since both sides are \cadlag processes this extends to all $t\in[0,T]$ (outside of a $\Prob$-null set) and we can take the $\mcS^1$-norm followed by H\"older's inequality to get that
\begin{align*}
\|Y^{v,k}-Y^{v,k'}\|_{\mcS^p}^p&\leq C\|(R^v)^{1/r'}(\tilde R^v R^v\bar R^v)^{1/rr'}(1+\bar K^{v,8pr^4})^{1/r^4}\|_{\mcS^1}/(k')^{1/2r}
\\
&\leq C\|R^v\|_{\mcS^1}^{1/r'}\|(\tilde R^v R^v\bar R^v)^{1/r'}(1+\bar K^{v,8pr^4})^{1/r^3}\|_{\mcS^1}^{1/r}/(k')^{1/2r}
\\
&\leq C\|R^v\|_{\mcS^1}^{1/r'}\|\tilde R^v R^v\bar R^v\|_{\mcS^1}^{1/rr'}\|(1+\bar K^{v,8pr^4})^{1/r^2}\|_{\mcS^1}^{1/r^2}/(k')^{1/2r}
\\
&\leq C/(k')^{1/2r},
\end{align*}
where $C>0$ is independent of $k,k'$. The last inequality holds since there is a $C>0$ such that
\begin{align*}
\|\tilde R^v R^v\bar R^v\|_{\mcS^1}\leq \|\tilde R^v\|_{\mcS^3}\| R^v\|_{\mcS^3}\|\bar R^v\|_{\mcS^3}\leq C,
\end{align*}
for all $v\in\mcU^f$. Finally, taking the limit as $k\to\infty$, (\emph{i}) and Fatou's lemma gives that $\|\bar Y^{v} - Y^{v,k'}\|^2_{\mcS^2} = C/(k')^{1/2r}$.\\

For the third claim we note that appealing to the corollaries of lemmas~\ref{lem:N-bound} and \ref{lem:seq-U-V-bound} rather than to the lemmas themselves gives analogously that for each $\beta\in\mcI(t)$ we have
\begin{align*}
|Y^{v\circ(t,\beta),k}_t-Y^{v\circ(t,\beta),k'}_t|^p&\leq   C(R^v_t)^{1/r'}(\tilde R^v_t R^v_t\bar R^v_t)^{1/rr'}(1+\bar K^{v,8pr^4}_t)^{1/r^4}/(k')^{1/2r}.
\end{align*}
Now, continuity and measurable selection implies that there is a $\beta\in \mcI(t)$ such that
\begin{align*}
\sup_{b\in U}|Y^{v\circ(t,b),k}_t-Y^{v\circ(t,b),k'}_t|= |Y^{v\circ(t,\beta),k}_t-Y^{v\circ(t,\beta),k'}_t|
\end{align*}
which leads us to conclude that
\begin{align*}
\sup_{b\in U}|Y^{v\circ(t,b),k}_t-Y^{v\circ(t,b),k'}_t|^p\leq  C(R^v_t)^{1/r'}(\tilde R^v_t R^v_t\bar R^v_t)^{1/rr'}(1+\bar K^{v,8pr^4}_t)^{1/r^4}/(k')^{1/2r}
\end{align*}
and the result follows similarly to the above since, under Hypothesis {\bf RBSDE.}$k$, the left hand side is continuous.\qed\\

\begin{prop}\label{prop:seq-existence}
There is a family $(\bar Z^{v},\bar K^v:v\in\mcU^f)$ such that $(\bar Y^v,\bar Z^{v},\bar K^v:v\in\mcU^f)$ is a solution to \eqref{ekv:seq-bsde-FH}.
\end{prop}

\noindent\emph{Proof.} Having established that $\|\sup_{b\in U}|Y^{v\circ(\cdot,b),k}_\cdot-Y^{v\circ(\cdot,b),k'}_\cdot|\|_{\mcS^p}\to 0$ as $k,k'\to\infty$ in the previous proposition it follows by Proposition~\ref{prop:rbsde-solu-FH} that $\|Z^{v,k}-Z^{v,k'}\|_{\mcH^2}\to 0$ as $k,k'\to\infty$. In particular, $(Z^{v,k})_{k\geq 0}$ is a Cauchy sequence in the Hilbert space $\mcH^2$ and we conclude that there is a $\bar Z^v\in \mcH^2$ such that $Z^{v,k}\to \bar Z^v$ in $\mcH^2$.

Now, letting $\bar K^v$ be defined by $\bar K^v=0$ and
\begin{align*}
\bar K^{v}_T-\bar K^{v}_t&=\sup_{r\in[t,T]}\Big(\xi^v+\int_r^Tf^{v}(s,\bar Y^{v}_s,\bar Z^{v}_s)ds-\int_r^T\bar Z^{v}_sdW_s-\sup_{b\in U}\{\bar Y^{v\circ(r,b)}_r-c^{v}(r,b)\}\Big)^-
\end{align*}
we note that $\|(\bar K^v-K^{v,k})\ett_{[0,\eta_l]}\|_{\mcS^2}\to 0$ as $k\to\infty$ where $\eta_l:=\inf\{s\geq 0:L_s^v\geq l\}\wedge T$ and by Lemma~\ref{lem:seq-unif-bound-FH} we have that $\bar K^v\in\mcS^2$. Since $L^v$ is continuous, and thus has $\Prob$-a.s.~bounded trajectories, we find that
\begin{align*}
\begin{cases}
  \bar Y^{v}_t=\xi^v+\int_t^T f^v(s,\bar Y^{v}_s,\bar Z^{v}_s)ds-\int_t^T \bar Z^{v}_sdW_s+ \bar K^{v}_T-\bar K^{v}_t,\quad\forall t\in[0,T], \\
  \bar Y^{v}_t\geq\sup_{b\in U}\{\bar Y^{v\circ(t,b)}_t-c^v(t,b)\},\quad\forall t\in[0,T],\\
  \int_0^T(\bar Y^{v}_t-\sup_{b\in U}\{\bar Y^{v\circ(t,b)}_t-c^v(t,b))\})dK^{v}_t=0.
\end{cases}
\end{align*}
Finally, the map $(t,b)\mapsto \bar Y^{v\circ(t,b)}_t\in\mcO_\bbF$ by uniform convergence.\qed

% (\emph{iii}) of Proposition~\ref{prop:seq-rbsde_k-FH}

\subsection{Uniqueness by a verification argument}

\begin{thm}\label{thm-verification-FH}
The finite horizon sequential system of reflected BSDEs \eqref{ekv:seq-bsde-FH} admits a unique solution $(Y^v,Z^v,K^v:v\in\mcU^f)$ and $Y^v_t=\esssup_{u\in\mcU^f_t}U^{v,u}_t=U^{v,u^*}_t$ with $u^*=(\tau^*_1,\ldots,\tau^*_{N^*};\beta^*_1,\ldots,\beta^*_{N^*})\in\mcU^f_t$ defined as:
\begin{itemize}
  \item $\tau^*_{j}:=\inf \Big\{s \geq \tau^*_{j-1}:\:Y_s^{v\circ[u^*]_{j-1}}=\sup_{b\in U} \{Y^{v\circ[u^*]_{j-1}\circ(s,b)}_s-c^{v\circ[u^*]_{j-1}}(s,b)\}\Big\}\wedge T$,
  \item $\beta^*_j\in\mathop{\arg\max}_{b\in U}\{Y^{v\circ [u^*]_{j-1}\circ(\tau^*_j,b)}_{\tau^*_j}-c^{v\circ [u^*]_{j-1}}(\tau^*_j,b)\}$
\end{itemize}
and $N^*=\sup\{j:\tau^*_j<T\}$, with $\tau_0^*:=t$.
\end{thm}

\noindent\emph{Proof.} Assume that $(Y^v,Z^v,K^v:v\in\mcU^f)$ is a solution to \eqref{ekv:seq-bsde-FH} (\ie $(Y^v:v\in\mcU)$ is a consistent family such that $(t,b)\mapsto Y^{v\circ(t,b)}_t$ is continuous and it satisfies equation \eqref{ekv:seq-bsde-FH}). Using Proposition~\ref{prop:rbsde-solu-FH} together with consistency $k$ times, gives that
\begin{align}\nonumber
Y^{v}_t&=Y^{v\circ [u^*]_{N^*\wedge k}}_{\tau^*_k}+\int_t^{\tau^*_k} \sum_{j=0}^{k-1}\ett_{[\tau^*_j,\tau^*_{j+1})}(s)f^{v\circ u^*}(s,Y^{v\circ [u^*]_j}_s,Z^{v\circ [u^*]_j}_s)ds
\\
&\quad - \int_t^{\tau^*_k}\sum_{j=0}^{k-1}\ett_{[\tau^*_j,\tau^*_{j+1})}(s)Z^{v\circ [u^*]_j}_s dW_s-\sum_{j=1}^{N^*\wedge k}c^{v\circ [u^*]_{j-1}}(\tau^*_j,\beta^*_j).\label{ekv:Yv-trunk}
\end{align}
Now, by the definition of a solution to~\eqref{ekv:seq-bsde-FH} the sequence $Y^{v\circ [u^*]_{N^*\wedge k}}_{\tau^*_k}$ is uniformly bounded in $L^2(\Omega,\bbQ^{v,u^*})$ and repeating the argument in the proof of Lemma~\ref{lem:N-bound} implies that $\Prob[N^*<\infty]=1$ and thus that $u^*\in\mcU^f$. Letting $k\to\infty$ in \eqref{ekv:Yv-trunk} then gives
\begin{align*}
Y^{v}_t=\xi^{v\circ u^*}+\int_t^{T} f^{v\circ u^*}(s,U^{v, u^*}_s,V^{v, u^*}_s)ds - \int_t^{T}V^{v, u^*}_sdW_s-\sum_{j=1}^{N^*}c^{v\circ [u^*]_{j-1}}(\tau^*_j,\beta^*_j)
\end{align*}
and uniqueness follows.\\

Concerning optimality let $\tilde u=(\tilde \tau_1,\ldots,\tilde \tau_{\tilde N};\tilde \beta_1,\ldots,\tilde \beta_{\tilde N})\in\mcU^f$ and note that if $\tilde N\geq 1$, then
\begin{align*}
Y^{v}_t\geq U^{v\circ (\tilde \tau_1,\tilde \beta_1),[\tilde u]_{2:}}_{\tilde \tau_1}+\int_t^{\tilde \tau_1}f^{v\circ \tilde u}(s,U^{v,\tilde u}_s,V^{v,\tilde u}_s)ds - \int_t^{\tilde \tau_1} V^{v,\tilde u}_sdW_s-c^{v}(\tilde \tau_1,\tilde \beta_1),
\end{align*}
where $[u]_{2:}:=(\tilde \tau_2,\ldots,\tilde \tau_{\tilde N};\tilde \beta_2,\ldots,\tilde \beta_{\tilde N})$. Successively repeating this process while considering the fact that $\tilde u\in\mcU^f$ eventually leads us to the conclusion that $Y^v_t\geq U^{v,\tilde u}_t$.\qed\\

%%%%%%%%%%%%%%%%%%%%%%%%%%%%%%%%%%%%%%%%%%%%%%%%%%%%%%%%%%%%%%%%%%%%%%%%%%%%%%%%%%%%%%%%%%%%%%%%%%%%%%%%%%%%%%%%%%%%%%%%%%%%%%%%%%%%%%%%%%%%%%%%%%
%%%%%%%%%%%%%%%%%%%%%%%%%%%%%%%%%%%%%%%%%%%%%%%%%%%%%%%%%%%%%%%%%%%%%%%%%%%%%%%%%%%%%%%%%%%%%%%%%%%%%%%%%%%%%%%%%%%%%%%%%%%%%%%%%%%%%%%%%%%%%%%%%%
%%%%%%%%%%%%%%%%%%%%%%%%%%%%%%%%%%%%%%%%%%%%%%%%%%%%%%%%%%%%%%%%%%%%%%%%%%%%%%%%%%%%%%%%%%%%%%%%%%%%%%%%%%%%%%%%%%%%%%%%%%%%%%%%%%%%%%%%%%%%%%%%%%

\section{Application to robust impulse control\label{sec:robust-impulse}}
We now apply the above results to find weakly optimal solutions to robust impulse control problems. In particular, we are interested in finding a pair $(u^*,\alpha^*)\in\mcU^f\times\mcA$, a probability measure $\bbQ$ and a corresponding Brownian motion $W^\bbQ$ such that
\begin{align}
J(u^*,\alpha^*)=\sup_{u\in\mcU^f}\inf_{\alpha\in\mcA}J(u,\alpha)\label{ekv:maxmin}
\end{align}
when \eqref{ekv:forward-sde1}-\eqref{ekv:forward-sde2} admits its solution on $(\Omega,\mcF,\bbF,\bbQ,W^\bbQ)$.

Throughout, we assume the following forms on the drift and volatility terms in the forward SDE \eqref{ekv:forward-sde1}-\eqref{ekv:forward-sde2},
\begin{align*}
a(t,x,\alpha)=\left[\begin{array}{c}a_1(t,x)\\ a_2(t,x,\alpha)\end{array}\right] \quad{\rm and}\quad \sigma(t,x)=\left[\begin{array}{cc}\sigma_{1,1}(t,x)& 0\\ \sigma_{2,1}(t,x) & \sigma_{2,2}(t,x)\end{array}\right],
\end{align*}
where $a$ is of at most linear growth in the data $x$ and $\sigma$ is uniformly bounded. The drift is split into two terms $a_1:[0,T]\times \bbD\to\R^{d_1}$ (we let $\bbD$ denote the set of all \cadlag functions $x:[0,T]\to\R^d$) and $a_2:[0,T]\times \bbD\times A\to\R^{d_2}$, with $d=d_1+d_2$ the total dimension. The diffusion coefficient has a component $\sigma_{2,2}:[0,T]\times \bbD\to\R^{d_2\times d_2}$ that has an inverse, $\sigma_{2,2}^{-1}$, which is uniformly bounded on $[0,T]\times \bbD$.

For the purpose of solving \eqref{ekv:maxmin} we let $f^u$ be given by
\begin{align}
f^u(t,\omega,y,z):=\inf_{\alpha\in A}H^u(t,\omega,z,\alpha)=:H^{*,u}(t,\omega,z),\label{ekv:f-u-def}
\end{align}
where\footnote{We use the notation $(x_s)_{s\leq t}$ in arguments to emphasise that a function, for example, $\phi:[0,T]\times\bbD\times A\to \R$ at time $t$ only depend on the trajectory of $x$ on $[0,t]$.}
\begin{align*}
H^u(t,\omega,z,\alpha):=z\breve a(t,(X^{u}_s)_{ s\leq t},\alpha)+\phi(t,X^{u}_t,\alpha),
\end{align*}
with
\begin{align*}
\breve a(t,x,\alpha):=\left[\begin{array}{c}0 \\ \sigma_{2,2}^{-1}(t,x)a_2(t,x,\alpha)\end{array}\right]
\end{align*}
and $X^u$ is the unique solution to the impulsively controlled forward SDE
\begin{align}
X^{u}_t&=x_0+\int_0^t\tilde a(s,(X^{u}_r)_{ r\leq s})dt+\int_0^t\sigma(s,(X^{u}_r)_{ r\leq s})dW_s,\quad{\rm for}\, t\in [0,\tau_{1}),\label{ekv:driftless-sde1}
\\
X^{u}_{t}&=\Gamma(\tau_j,X^{[u]_{j-1}}_{\tau_j},\beta_j)+\int_{\tau_j}^t\tilde a(s,(X^{u}_r)_{ r\leq s})ds+\int_0^t\sigma(s,(X^{u}_r)_{ r\leq s})dW_s,\quad{\rm for}\, t\in [{\tau_j},\tau_{j+1}),\label{ekv:driftless-sde2}
\end{align}
with
\begin{align*}
\tilde a(t,x):=\left[\begin{array}{c}a_1(t,x)\\ 0\end{array}\right].
\end{align*}

Our approach to solving the above optimization problem is to define a measure $\bbQ^{u,\alpha}$ under which $W^{u,\alpha}_t=W_t-\int_0^t\breve a(s,(X^{u}_r)_{ r\leq s},\alpha^*_s)ds$ is a Brownian motion, where $\alpha^*$ is a measurable selection of a minimizer in \eqref{ekv:f-u-def}. In particular, we note that for any $(u,\alpha)\in\mcU^f\times \mcA$, the 6-tuple $(\Omega,\mcF,\bbF,\bbQ^{u,\alpha},X^u,W^{u,\alpha})$ is a weak solution to \eqref{ekv:forward-sde1}-\eqref{ekv:forward-sde2} with impulse control $u$ and continuous control $\alpha$.

Before we move on to show optimality of the above scheme, we give assumptions on $a,\sigma$ and $\Gamma$ and $\phi,\psi$ and $\ell$ under which the sequential system \eqref{ekv:seq-bsde-FH} with driver given by \eqref{ekv:f-u-def} attains a unique solution.

\begin{ass}\label{ass:onSFDE}
For any $t,t'\geq 0$, $b,b'\in U$, $\xi,\xi'\in\R^d$, $x,x'\in\bbD$ and $\alpha\in A$ and for some $\rho\geq 0$ we have:
\begin{enumerate}[i)]
  \item\label{ass:onSFDE-Gamma} The function $\Gamma:[0,T]\times\bbD\times U\to\R^d$ satisfies the Lipschitz condition
  \begin{align*}
    |\Gamma(t,(x_s)_{s\leq t},b)-\Gamma(t',(x'_s)_{s\leq t'},b')|&\leq C(\int_{0}^{t\wedge t'}|x'_s-x_s|ds+|x'_{t'}-x_t|
    \\
    &\quad+(|t'-t|+|b'-b|)(1+\sup_{s\leq t}|x_s|+\sup_{s\leq t'}|x'_s|))
  \end{align*}
  and the growth condition
  \begin{align*}
    |\Gamma(t,(x_s)_{s\leq t},b)|\leq K_\Gamma\vee |x_t|.
  \end{align*}
  for some constant $K_\Gamma>0$.
  \item\label{ass:onSFDE-a-sigma} The coefficients $a:[0,T]\times\bbD\times A\to\R^{d}$ and $\sigma:[0,T]\times\bbD\to\R^{d\times d}$ are continuous in $t$ (and $\alpha$ when applicable) and satisfy the growth conditions
  \begin{align*}
    |a(t,(x_s)_{s\leq t},\alpha)|&\leq C(1+\sup_{s\leq t}|x_s|),
    \\
    |\sigma(t,(x_s)_{s\leq t})|&\leq C
  \end{align*}
  and the Lipschitz continuity
  \begin{align*}
    |\tilde a(t,(x_s)_{s\leq t})-\tilde a(t,(x'_s)_{s\leq t})|+|\sigma(t,(x_s)_{s\leq t})-\sigma(t,(x'_s)_{s\leq t})|&\leq C\sup_{s\leq t}|x'_s-x_s|,
    \\
    \int_0^t |\tilde a(s,(x_r)_{r\leq s})-\tilde a(s,(x'_r)_{r\leq s})|ds&\leq C\int_{0}^t|x'_s-x_s|ds
    \\
    \int_0^t |\sigma(s,(x_r)_{r\leq s})-\sigma(s,(x'_r)_{r\leq s})|^2ds&\leq C\int_{0}^t|x'_s-x_s|^2ds.
  \end{align*}
  Moreover, for each $(t,x)\in[0,T]\times\bbD$, the matrix $\sigma_{2,2}(t,(x_s)_{s\leq t})$ has an inverse, $\sigma_{2,2}^{-1}(t,(x_s)_{s\leq t})$, that is uniformly bounded on $[0,T]\times \bbD$ and
  \begin{align*}
    \int_0^t |\breve a(s,(x_r)_{r\leq s},\hat\alpha(s))-\breve a(s,(x'_r)_{r\leq s},\hat\alpha(s))|^2 ds&\leq C\int_{0}^t|x'_s-x_s|^2ds,
  \end{align*}
for all measurable functions $\hat\alpha:[0,T]\to A$.
  \item\label{ass:onSFDE-phi} The running reward $\phi:[0,T]\times \R^d\times A\to\R$ is $\mcB([0,T]\times\R^d\times A)$-measurable, continuous in $\alpha$ and satisfies the growth condition
  \begin{align*}
    |\phi(t,\xi,\alpha)|\leq C^g_\phi(1+|\xi|^\rho)
  \end{align*}
  for some $C^g_\phi>0$ and all $\xi\in\R^d$, and locally Lipschitz in $\xi$, \ie there is a nondecreasing function $C^L_\phi:\R_+\to\R_+$ such that for each $K>0$,
  \begin{align*}
    |\phi(t,\xi,\alpha)-\phi(t,\xi',\alpha)|\leq C^L_\phi(K)|\xi-\xi'|,
  \end{align*}
  whenever $|\xi|\vee|\xi'|\leq K$.
  \item\label{ass:onSFDE-psi} The terminal reward $\psi:\R^d\to\R$ is $\mcB(\R^d)$-measurable, and satisfies the growth condition
  \begin{align*}
    |\psi(\xi)|\leq C^g_\psi(1+|\xi|^\rho)
  \end{align*}
  for some $C^g_\psi>0$ and all $\xi\in\R^d$, and locally Lipschitz, \ie there is a nondecreasing function $C^L_\psi:\R_+\to\R_+$ such that for each $K>0$,
  \begin{align*}
    |\psi(\xi)-\psi(\xi')|\leq C^L_\psi(K)|\xi-\xi'|,
  \end{align*}
  whenever $|\xi|\vee|\xi'|\leq K$.
  \item\label{ass:onSFDE-c} The intervention cost $\ell:[0,T]\times \R^d\times U\to \R_+$ is jointly continuous in $(t,\xi,b)$, bounded from below, \ie
  \begin{align*}
    \ell(t,\xi,b)\geq\delta >0,
  \end{align*}
  and locally Lipschitz in $\xi$ and H\"older continuous in $t$, \ie there is a nondecreasing function $C^L_\ell:\R_+\to\R_+$ such that for each $K>0$,
  \begin{align*}
    |\ell(t,\xi,b)-\ell(t',\xi',b)|\leq C^L_\ell(K)|\xi-\xi'|+C^H_\ell|t'-t|^\varsigma,
  \end{align*}
  whenever $|\xi|\vee|\xi'|\leq K$ for some $C^H_\ell,\varsigma>0$.
\end{enumerate}
\end{ass}

Under these assumptions, we note that $f^u$ as defined in \eqref{ekv:f-u-def} is stochastic Lipschitz with an admissible Lipschitz coefficient
\begin{align}
  L^{v}_t&:=\sup_{s\in[0,t]}\sup_{\alpha\in A}|\breve a(\cdot,(X^v_s)_{ s\leq \cdot},\alpha)|\vee C,\label{ekv:L-u-def}
\end{align}
where $C>0$ is chosen to eliminate jumps. Moreover, for some $k_L>0$, we have $L^v_t\leq k_L(1+\sup_{s\in[0,t]}|X^v_s|)$.

\subsection{Some preliminary estimates}
We now present some preliminary estimates of moments and stability of solutions to \eqref{ekv:driftless-sde1}-\eqref{ekv:driftless-sde2}. Towards the end of the section we will prove that any necessary changes of measure are admissible.% but for now we state our results in terms of solutions to \eqref{ekv:driftless-sde1}-\eqref{ekv:driftless-sde2} with an added coefficient in the drift term.

%In this regard we let $\zeta$ be a $\mcP_\bbF$-measurable process with trajectories in $\bbD$ and assume that the process $\tilde X^{v,\zeta}$ solves the FSDE \eqref{ekv:driftless-sde1}-\eqref{ekv:driftless-sde2} with $W$ replaced by $W+\int_\eta^\cdot \zeta_sds$, \ie
%\begin{align}
%d\tilde X^{u,\zeta}_t&=(\tilde a(t,(\tilde X^{u,\zeta}_s)_{ s\leq t})+\sigma(t,(\tilde X^{u,\zeta}_s)_{ s\leq t})\zeta_t)dt+\sigma(t,(\tilde X^{u,\zeta}_s)_{ s\leq t})dW_t,\mbox{ for }t\in [\tau_j,\tau_{j+1}),\label{ekv:driftadded-sde1}
%\\
%\tilde X^{[u]_j}_{\tau_j}&=\Gamma(\tau_j,\tilde X^{[u]_{j-1},\zeta}_{\tau_j},\beta_j),\label{ekv:driftadded-sde2}
%\\
%\tilde X^{\emptyset,\zeta}_0&=x_0,\label{ekv:driftadded-sde3}
%\end{align}
%Recall that $\mcK^v$ was defined to be the set of all $\mcP_\bbF$-measurable processes $(\zeta_s)_{t\leq s\leq T}$ such that $|\zeta_s|\leq L^{v}_s$.

\begin{prop}\label{prop:SFDEmoment}
Under Assumption~\ref{ass:onSFDE}, the FSDE \eqref{ekv:driftless-sde1}-\eqref{ekv:driftless-sde2} admits a unique solution for each $u\in\mcU$. Furthermore, the solution has moments of all orders, in particular we have for $p\geq 0$, that
\begin{align}\label{ekv:SFDEmoment}
\sup_{u\in\mcU^f}\E\Big[\sup_{t\in[0,T]}|X^{u}_t|^{p}\Big]\leq C,
\end{align}
where $C=C(p)$ and
\begin{align}\label{ekv:SFDEmoment2}
\sup_{v\in\mcU^f}\E\Big[\sup_{t\in [0,T]}\esssup_{u\in\mcU^f_t}|\E\big[\sup_{s\in[t,T]}|X^{v(t)\circ u}_s|^{\rho}\big|\mcF_t\big]|^p\Big]\leq C,
\end{align}
where $C=C(\rho,p)$.
\end{prop}

\noindent\emph{Proof.} By repeated use of Theorem 3.2 in~\cite{AgramOksen} existence and uniqueness of solutions to \eqref{ekv:driftless-sde1}-\eqref{ekv:driftless-sde2} follows as $N<\infty$, $\Prob$-a.s.

By Assumption~\ref{ass:onSFDE}.(\ref{ass:onSFDE-Gamma}) we get, for $t\in [\tau_{j},T]$, using integration by parts, that
\begin{align*}
|X^{[u]_j}_t|^2&= |X^{[u]_j}_{\tau_{j}}|^2+2\int_{\tau_{j}+}^t X^{[u]_j}_{s}dX^{[u]_j}_s+\int_{\tau_{j}+}^t d[X^{[u]_j},X^{[u]_j}]_s
\\
&\leq K^2_\Gamma\vee |X^{[u]_{j-1}}_{\tau_{j}}|^2+2\int_{\tau_{j}+}^t X^{[u]_j}_{s} dX^{[u]_j}_s+\int_{\tau_{j}+}^t d[X^{[u]_j},X^{[u]_j}]_s.
\end{align*}
We note that if $|X^{[u]_j}_t|> K_\Gamma$ and $|X^{[u]_j}_s|\leq K_\Gamma$ for some $s\in[0,t)$ then there is a largest time $\theta<t$ such that $|X^{[u]_j}_{\theta}|\leq K_\Gamma$. This means that during the interval $(\theta,t]$ interventions will not increase the magnitude $|X^{[u]_j}|$. By induction, since $|x_{0}|$ is finite, we find that
\begin{align}\label{ekv:X2-bound}
|X^{[u]_j}_t|^2&\leq C+\sum_{i=0}^{j} \Big\{2\int_{\theta\vee(\tilde\tau_{i}+)}^{t\wedge\tilde\tau_{i+1}} X^{[u]_i}_{s}dX^{[u]_i}_s+\int_{\theta\vee(\tilde\tau_{i}+)}^{t\wedge\tilde\tau_{i+1}} d[X^{[u]_i},X^{[u]_j}]_s\Big\}
\end{align}
for all $t\in[0,T]$, where $\theta:=\sup\{s\geq 0 : |X^u_s|\leq K_\Gamma\}\vee 0$, $\tilde\tau_0+=0$, $\tilde\tau_i=\tau_i$ for $i=1,\ldots,j$ and $\tilde\tau_{j+1}=\infty$.

Now, since $X^{[u]_i}$ and $X^{[u]_j}$ coincide on $[0,\tau_{i+1\wedge j+1})$ we have
\begin{align*}
\sum_{i=0}^{j}\int_{\theta\vee\tilde\tau_{i}+}^{t\wedge\tilde\tau_{i+1}} X^{[u]_i}_{s} dX^{[u]_i}_s
&=\int_{\theta}^t X^{[u]_j}_{s}\tilde a(s,(X^{[u]_j}_r)_{r\leq s})ds+\int_{\theta}^{t}X^{[u]_j}_{s}\sigma(s,(X^{[u]_j}_r)_{r\leq s})dW_s,
\end{align*}
and
\begin{align*}
\sum_{i=0}^{j} \int_{\theta\vee\tilde\tau_{i}+}^{t\wedge\tilde\tau_{i+1}} d[X^{[u]_i},X^{[u]_j}]_s&=\int_{\theta}^{t} \sigma^2(s,(X^{[u]_j}_r)_{r\leq s})ds\leq Ct.
\end{align*}
Inserted in \eqref{ekv:X2-bound} this gives
\begin{align*}
|X^{[u]_j}_t|^2&\leq C(1+t)+2\int_{\theta}^t X^{[u]_j}_{s}\tilde a(s,(X^{[u]_j}_r)_{r\leq s})ds+2\int_{\theta}^{t}X^{[u]_j}_{s}\sigma(s,(X^{[u]_j}_r)_{r\leq s})dW_s
\\
&\leq C\Big(1+t+\int_{0}^{t}\sup_{r\in [0,T]}|X^{[u]_j}_{r}|^2ds+\sup_{s\in[0,t]}\Big|\int_{0}^{s}X^{[u]_j}_r\sigma(r,(X^{[u]_j}_\varsigma)_{\varsigma\leq r})dW_r\Big|\Big).
\end{align*}
The Burkholder-Davis-Gundy inequality now gives that for $p\geq 4$
\begin{align*}
\E\Big[\sup_{r\in[0,t]}|X^{[u]_j}_r|^{p}\Big]\leq C\big(1+t+\int_{0}^{t}\E\Big[\sup_{r\in [0,s]}|X^{[u]_i}_{r}|^{p}\Big]ds+\E\Big[\big(\int_{0}^{t}|X^{[u]_j}_s|^{p/2}ds\big)\Big]\big)
\end{align*}
and Gr\"onwall's lemma gives that
\begin{align}\label{ekv:satisfies-Novikov}
\E\Big[\sup_{t\in [0,T]}|X^{[u]_j}_t|^{p}\Big]&\leq C.
\end{align}
Similarly, we find that
\begin{align}\label{ekv:moment_steg1}
\E\Big[\sup_{s\in[t,T]}|X^{[u]_j}_s|^{p}\big|\mcF_t\Big]&\leq C(1+ \sup_{s\in[0,t]}|X^{[u]_j}_s|^{p}),
\end{align}
$\Prob$-a.s., where the constant $C=C(T,p)$ does not depend on $u$, $\alpha$ or $j$ and \eqref{ekv:SFDEmoment} follows by letting $j\to\infty$ on both sides and using Fatou's lemma and dominated convergence. Applying \eqref{ekv:moment_steg1} to the left-hand side of \eqref{ekv:SFDEmoment2} we get
\begin{align}\nonumber
\E\Big[\sup_{t\in [0,T]}\esssup_{u\in\mcU^f_t}|\E\big[\sup_{s\in[t,T]}|X^{v(t)\circ u}_s|^{\rho}\big|\mcF_t\big]|^p\Big]&\leq C(1+\E\Big[\sup_{t\in [0,T]}|\E\big[\sup_{s\in[0,t]}|X^{v}_s|^{\rho}\big|\mcF_t\big]|^p
\\
&\qquad+\sup_{(t,b)\in [0,T]\times U}|\E\big[|\Gamma(t,X^{v}_t,b)|^{\rho}\big|\mcF_t\big]|^p\Big])\nonumber
%\\
%&\leq C(1+\E\Big[\sup_{t\in [0,T]}|\E\big[\sup_{s\in[0,t]}|X^{v}_s|^{\rho}\big|\mcF_t\big]|^p\Big])\nonumber
\\
&\leq C(1+\E\big[\sup_{t\in [0,T]}|X^{v}_t|^{qp}\big])\label{ekv:SFDE-moments-relation}
\end{align}
and the desired result follows by \eqref{ekv:SFDEmoment}.\qed\\

For any $\kappa\geq 1$ and all $\vecv:=(t_1,\ldots,t_\kappa;b_1,\ldots,b_\kappa)\in\mcD^\kappa$ and $\vecv':=(t'_1,\ldots,t'_\kappa;b'_1,\ldots,b'_\kappa)\in\mcD^\kappa$ we define the set $\Xi_{\vecv,\vecv'}:=\cup_{j=1}^\kappa [\underline t_j,\bar t_j)$, with $\underline t_j:=t_j\wedge t'_j$ and $\bar t_j:=t_j\vee t'_j$.

\begin{lem}\label{lem:SFDEflow}
For each $k,\kappa\geq 0$ and $p\geq 1$, there is a $C\geq 0$ such that
\begin{align*}
  \E\Big[\sup_{t\in [0,T]}\esssup_{u\in\mcU^k_t}\E\Big[\sup_{s\in[0,T]\setminus\Xi_{\vecv,\vecv'}}|X^{v\circ\vecv'\circ u}_{s}-X^{v\circ\vecv\circ u}_{s}|^{2p}\Big|\mcF_t\Big]\Big]\leq C\|\vecv'-\vecv\|^p_{\mcD^f},
\end{align*}
for all $v\in\mcU^f$ and all $\vecv,\vecv'\in\mcD^\kappa$.
\end{lem}

\noindent\emph{Proof.} To simplify notation we let $X^l:=X^{v\circ[\vecv\circ u]_l}$ and $\tilde X^l:=X^{v\circ[\vecv'\circ u]_l}$ for $l=0,\ldots,\kappa+k$. Moreover, we let $\delta X^l:=X^l-\tilde X^l$ and set $\delta X:=\delta  X^{\kappa+k}$. We have
\begin{align*}
|\delta  X^{j}_{\bar t_j}|&\leq C(\int_{0}^{\underline t_j}|\delta  X^{{j-1}}_{s}|ds+ |\delta X^{{j-1}}_{\bar t_j}|
+|X^{{j-1}}_{\bar t_j}-X^{{j-1}}_{t_j}| +|\tilde X^{{j-1}}_{\bar t_j}-\tilde X^{{j-1}}_{t'_j}|
    \\
    &\quad+(|t'_j-t_j|+|b'_j-b_j|)(1+\sup_{s\in [0,t_j]}|X^{{j-1}}_{s}|+\sup_{s\in [0,t'_j]}|\tilde X^{{j-1}}_{s}|))
\end{align*}
For $j=1$, this gives
\begin{align*}
|\delta X^{1}_{\bar t_1}|\leq C(|X^{v}_{t_1}-X^{v}_{t'_1}|+(|t'_1-t_1|+|b'_1-b_1|)(1+\sup_{s\in [0,\bar t_1]}|X^{v}_{s}|)).
\end{align*}
and by induction we find that
\begin{align*}
|\delta  X^{\kappa}_{\bar t_\kappa}|&\leq \sum_{j=1}^{\kappa} C^{\kappa+1-j}(\int_{0}^{\underline t_{j}}|\delta  X^{{j-1}}_{s}|ds+|\delta  X^{{j-1}}_{\bar t_{j}}|-|\delta  X^{{j-1}}_{\bar t_{j-1}}|
\\
&\quad+|X^{{j-1}}_{\bar t_{j}}-X^{{j-1}}_{t_{j}}| +|\tilde X^{{j-1}}_{\bar t_j}-\tilde X^{{j-1}}_{t'_j}|
\\
&\quad+(|t'_{j}-t_{j}|+|b'_{j}-b_{j}|)(1+\sup_{s\in [0,t_{j}]}|X^{{j-1}}_{s}|+\sup_{s\in [0,t'_{j}]}|\tilde X^{{j-1}}_{s}|)),
\end{align*}
with $\bar t_0=0$. In the above, we note that for $j=1,\ldots,\kappa$ and $r\geq \bar t_j$,
\begin{align*}
|\delta  X^{{j}}_{r}|&\leq |\delta  X^{{j}}_{\bar t_{j}}|+\int_{\bar t_{j}}^{r}|\tilde a(s,(X^{{j}}_\varrho)_{ \varrho\leq s})-\tilde a(s,(\tilde X^{{j}}_\varrho)_{ \varrho\leq s})|ds
\\
&\quad+\Big|\int_{\bar t_{j}}^{r}(\sigma(s,(X^{{j}}_\varrho)_{ \varrho\leq s})-\sigma(s,(\tilde X^{{j}}_\varrho)_{ \varrho\leq s}))dW_s\Big|
\end{align*}
The Burkholder-Davis-Gundy inequality now gives that
\begin{align*}
\E\Big[\sup_{r\in[\bar t_j,s]}|\delta X^j_{r}|^{2p}\Big] &\leq C\E\Big[|\delta X^j_{\bar t_j}|^{2p}+ \Big(\int_{\bar t_j}^s|\tilde a(r,(X^{{j}}_\varrho)_{ \varrho\leq r})-\tilde a(r,(\tilde X^{{j}}_\varrho)_{ \varrho\leq r})|dr\Big)^{2p}
\\
& + \Big(\int_{\bar t_j}^{s}|\sigma(r,(\tilde X^j_\varrho)_{\varrho\leq r})-\sigma(r,(X^j_\varrho)_{\varrho\leq r})|^2dr\Big)^{p}\Big]
\\
&\leq C\E\Big[|\delta X^j_{\bar t_j}|^{2p}+ \big(\int_{0}^s|\delta X^j_{r}|^{2}dr\big)^p\Big]
\end{align*}
by the integral Lipschitz conditions on the coefficients, which by Gr\"onwall's lemma implies that
\begin{align*}
\E\Big[\sup_{r\in[\bar t_j,T]}|\delta X^j_{r}|^{2p}\Big] &\leq C\E\Big[|\delta X^j_{\bar t_j}|^{2p}+ \big(\int_0^{\bar t_j}|\delta X^j_{r}|^{2}dr\big)^p\Big].
\end{align*}
Furthermore, we have
\begin{align*}
|X^{{j-1}}_{\bar t_{j}}-X^{{j-1}}_{t_{j}}|\leq \int_{t_{j}}^{\bar t_{j}}|\tilde a(s,(X^{{j-1}}_r)_{ r\leq s})|ds+\Big| \int_{t_{j}}^{\bar t_{j}}\sigma(s,(X^{{j-1}}_r)_{ r\leq s})dW_s\Big|.
\end{align*}
Combining the above we find that
\begin{align*}
\E\Big[\sup_{s\in [0,T]\setminus \Xi_{\vecv,\vecv'}}|X^{v}_{s}-X^{v'}_{s}|^{2p}\Big]&\leq C\E\Big[\sum_{j=1}^\kappa (1+\sup_{s\in [0,T]}(|X^{{j-1}}_{s}|^{2p}+|\tilde X^{{j-1}}_{s}|^{2p}))(|t_j-t'_j|^{p}+|b_j-b'_j|^{p})\Big]
\\
&\leq C\sum_{j=1}^\kappa (|t_j-t'_j|^{p}+|b_j-b'_j|^{p}).
\end{align*}
Now, for $u=(\tau_1,\ldots,\tau_N;\beta_1,\ldots,\beta_N)\in\mcU^k$ we let $\theta:=\max\{j:\tau_j\leq \bar t_\kappa\}\vee 0$ and get, by repeating the above argument, that
\begin{align}\label{ekv:diff-est-1}
\E\big[\sup_{s\in[0,T]\setminus\Xi_{\vecv,\vecv'}}|\delta X^{\kappa+\theta}_{s}|^{2p}\big]&\leq C\sum_{j=1}^{\kappa+k} (|t_j-t'_j|^{p}+|b_j-b'_j|^{p})
\end{align}
where we have used that fact that $|\tau_j\vee t_\kappa-\tau_j\vee t'_\kappa|\leq |t_\kappa- t'_\kappa|$.\\

Assumption~\ref{ass:onSFDE}.(\ref{ass:onSFDE-Gamma}) gives that for $l=\theta+1,\ldots,N$, we have
\begin{align*}
|\delta X^{\kappa+l}_{\tau_l}| &\leq C(\int_{0}^{\tau_l}|\delta X^{\kappa+l-1}_{s}|ds + |\delta X^{\kappa+l-1}_{\tau_l}|).
\end{align*}
Now, for $s\geq \bar t_\kappa$,
\begin{align*}
\delta X_{s}&=\delta X_{\bar t_\kappa}+\sum_{l=\theta}^{N}\int_{(\bar t_\kappa\vee\tau_l)+}^{s\wedge\tau_{l+1}} d(\delta X^{\kappa+l})_r+\sum_{l=\theta+1}^N\ett_{[s\geq\tau_l]}(\delta X^{\kappa+l}_{\tau_l}-\delta X^{\kappa+l-1}_{\tau_l}),
\end{align*}
with $\tau_0=0$ and $\tau_{N+1}:=\infty$. Combining these we find by taking the absolute value on both sides that
\begin{align*}
|\delta X_{s}|&\leq |\delta X_{\bar t_\kappa}|+ C(\int_{\bar t_\kappa}^s|\tilde a(r,(\tilde X_\zeta)_{\zeta\leq r})-\tilde a(r,(X_\zeta)_{\zeta\leq r})|dr
\\
&\quad +\sum_{l=\theta}^N |\int_{\bar t_\kappa\vee\tau_l}^{\tau_{l+1}\wedge s}\sigma(r,(\tilde X_\zeta)_{\zeta\leq r})-\sigma(r,(X_\zeta)_{\zeta\leq r})dW_r| +\int_t^s|\delta X_{r}|dr).
\end{align*}
As above, the Burkholder-Davis-Gundy inequality combined with the integral Lipschitz conditions on the coefficients and Gr\"onwall's lemma gives that
\begin{align*}%\setminus \Xi_{\vecv,\vecv'}
&\E\Big[\sup_{s\in[\bar t_\kappa,T]}|\delta X_s|^{2p}\Big] \leq C(\E\Big[|\delta X_{\bar t_\kappa}|^{2p}+\big(\int_0^{\bar t_j}|\delta X^j_{r}|^{2}dr\big)^p\Big],
\end{align*}
where $C$ does not depend on $u\in\mcU^k$. Now, the desired result follows by \eqref{ekv:diff-est-1}.\qed\\

A fundamental assumption in Section~\ref{sec:seq-rbsde-FH} is the existence of a $q'>1$ and a $C>0$ such that $\sup_{\zeta \in\mcK_0}\E^\bbQ\big[|\mcE(\zeta*W^\bbQ)_T|^{q'}\big]\leq C$ for all $\bbQ\in\PrM_0$. In the following two lemmas we show that since $L^v\leq k_L(1+\sup_{s\in[0,\cdot]}|X^v_s|)$, this statement is true.

\begin{lem}\label{lem:Lp-bound-if-martingale}
For $v\in\mcU^f$ and $\varsigma>1$, let $\Upsilon^{v,\varsigma}$ be the set of all $\mcP_\bbF$-measurable processes $\zeta$ with $|\zeta_t|\leq L^v_t$ for all $t\in[0,T]$ (outside of a $\Prob$-null set) such that $\E\big[\mcE(r\zeta*W)_T\big]=1$ for all $r\in[1,\varsigma]$. Then, there is a $q'>1$ such that $\sup_{u\in\mcU^f}\sup_{\zeta\in \Upsilon^{u,\varsigma}}\E\big[|\mcE(\zeta*W)_T|^{q'}\big]<\infty$.
\end{lem}

\noindent\emph{Proof.} For $q'\in[1,\varsigma]$ we have
\begin{align*}
\E\big[|\mcE(\zeta*W)_T|^{q'}\big]&=\E\big[\mcE(q'\zeta*W)_Te^{\frac{(q')^2-q'}{2}\int_0^T|\zeta_s|^2ds}\big]
\\
&=\E^{\bbQ^{q'\zeta}}\big[e^{\frac{(q')^2-q'}{2}\int_0^T|\zeta_s|^2ds}\big].
\end{align*}
Now, since $|\zeta|\leq k_L(1+\sup_{s\in[0,\cdot]}|X^v_s|)$ we have (see Lemma 1 in~\cite{Benes1971})
\begin{align*}
\E^{\bbQ^{q'\zeta}}\big[e^{\frac{(q')^2-q'}{2}\int_0^T|\zeta_s|^2ds}\big]\leq \E\big[e^{\frac{(q')^2-q'}{2}h(q')(1+\sup_{t\in [0,T]}|W_t|^2)}\big],
\end{align*}
where $h:\R_+\to\R_+$ is bounded on compacts and we conclude that the left hand side is finite for $q'>1$ sufficiently small.\qed\\

\begin{lem}\label{lem:is-martingale}
Let $\zeta$ be a $\mcP_\bbF$-measurable process with trajectories in $\bbD$ such that for some $v\in\mcU^f$, we have $|\zeta_t|\leq k_L(1+\sup_{s\in[0,t]}|X^v_s|)$ for all $t\in[0,T]$, then $\E\big[\mcE(\zeta*W)_T\big]=1$.
\end{lem}

\noindent\emph{Proof.} We will reach the result by adapting the proof of Lemma 7 in~\cite{Girsanov1960} to solutions of impulsively controlled FSDEs (see also Lemma 0 of Section 5 in \cite{Benes1971}). Since $(\mcE(\zeta*W)_t:0\leq t\leq T)$ is a $\Prob$-a.s.~non-negative local martingale it is a supermartingale and we only need to show that $\E\big[\mcE(\zeta*W)_T\big]\geq 1$. For $M\geq 0$ and $t\in[0,T]$, we define the sets
\begin{align*}
C_M(t):=\{x\in\bbD: k_L(1+\sup_{s\in[0,t]}|x_s|)<M\}
\end{align*}
Then for each $M\geq 0$, $(C_M(t))_{t\in [0,T]}$ is a non-increasing collection of open subsets of $\bbD$ and:
\begin{enumerate}[a)]
  \item If for some $x\in\bbD$ we have $x\in C_M(s)$ and $x\notin C_M(t)$ for some $0\leq s<t\leq T$, then there is a $t'\in (s,t]$ such that $x\in C_M(r)$ for all $r\in[0,t')$ and $x\notin C_M(t')$.
  \item For each $\epsilon>0$, there is a $M\geq 0$ such that $\Prob[X^v\in C_M(T)]>1-\epsilon$ for all $v\in\mcU^f$.
\end{enumerate}
Here, the second property follows from Proposition~\ref{prop:SFDEmoment}. Moreover, let
\begin{align*}
D_M(t):=\{\omega:X^{v,\zeta}\in C_M(t)\}\in\mcF_t,
\end{align*}
where $X^{u,\zeta}$ solves \eqref{ekv:forward-sde1}-\eqref{ekv:forward-sde2} with drift $a=\tilde{a}+\sigma\zeta$, \ie
\begin{align*}
X^{u,\zeta}_t&=x_0+\int_0^t(\tilde{a}(s,(X^{u,\zeta}_r)_{ r\leq s})+\sigma(s,(X^{u,\zeta}_r)_{ r\leq s})\zeta_s)ds+\int_0^t\sigma(s,(X^{u,\zeta}_r)_{ r\leq s})dW_s
\end{align*}
on $[0,\tau_1)$ and
\begin{align*}
X^{u,\zeta}_t&=\Gamma(\tau_j,X^{[u]_{j-1},\zeta}_{\tau_j},\beta_j)+\int_{\tau_j}^t(\tilde{a}(s,(X^{u,\zeta}_r)_{ r\leq s})+\sigma(s,(X^{u,\zeta}_r)_{ r\leq s})\zeta_s)ds+\int_{\tau_j}^t\sigma(s,(X^{u,\zeta}_r)_{ r\leq s})dW_s
\end{align*}
whenever $t\in [\tau_j,\tau_{j+1})$ for $j=1,\ldots,N$ with $\tau_{N+1}:=\infty$.

We first restrict our attention to the situation when $v\in\mcU^k$ for some $k\geq 0$, and note that (by arguing as in the proof of Lemma~\ref{lem:SFDEflow}) we have
\begin{align*}
|X^{v,\zeta}_t-X^{v}_t|&\leq C\big(\int_{0}^t\big(|\tilde{a}(s,(X^{v,\zeta}_r)_{ r\leq s})-\tilde{a}(s,(X^{v}_r)_{ r\leq s})|+|\sigma(s,(X^{v,\zeta}_r)_{ r\leq s})\zeta_s|+|X^{v,\zeta}_s-X^{v}_s|\big)ds
\\
&\quad+\Big|\int_0^{t}(\sigma(s,(X^{v,\zeta}_r)_{ r\leq s})-\sigma(s,(X^{v}_r)_{ r\leq s}))dW_s\Big|\big)
\\
&\leq C\big(\int_{0}^t\big(|X^{v,\zeta}_s-X^{v}_s|+|\zeta_s|\big)ds+\Xi^{v,\zeta}_t\big),
\end{align*}
where
\begin{align*}
\Xi^{v,\zeta}_{t}:=\sup_{t'\in[0,t]}\Big|\int_0^{t'}(\sigma(s,(X^{v,\zeta}_r)_{ r\leq s})-\sigma(s,(X^{v}_r)_{ r\leq s}))dW_s\Big|.
\end{align*}
Applying Gr\"onwall's inequality together with the fact that $|\zeta_t|\leq k_L(1+\sup_{s\in[0,t]}|X^{v}_s|)$ gives that for any $T'\in[0,T]$, we have
\begin{align*}
\sup_{t\in [0,T']}|X^{v,\zeta}_t-X^{v}_t|&\leq C\big(1+\int_{0}^{T'}\sup_{r\in [0,s]}|X^{v}_r|ds+\Xi^{v,\zeta}_{T'}\big),
\end{align*}
Now, for all $\omega\in D_M(T)$ we have $\sup_{t\in [0,T]}|X^{v,\zeta}_t|<M$ and we can apply Gr\"onwall's inequality once more to obtain
\begin{align*}
\sup_{t\in [0,T]}|X^{v}_t|&\leq C\big(1+M+\Xi^{v,\zeta}_T\big)
\end{align*}
Letting $E_M(t):=\{\omega\in D_M(t): \Xi^{v,\zeta}_t\leq M\}\in\mcF_t$ we note that
\begin{enumerate}[a)]\setcounter{enumi}{2}
  \item For $\omega\in E_M(t)$ we have $\zeta_t\leq C(1+M)$, where $C$ does not depend on $t$ or $M$.
\end{enumerate}
Now, set
\begin{align*}
\zeta^{M}_t:=\ett_{E_M(t)}\zeta_t
\end{align*}
and let $X^{v,\zeta,M}:=X^{v,\zeta^{M}}$.

Since $|\zeta^{M}_t|\leq C(1+M)$, the Novikov condition holds for any constant multiple of $\zeta^{M}$. In particular, we conclude that the $\bbQ_M$ defined by $d\bbQ_M:=\mcE(\zeta^{M}*W)_Td\Prob$ is a probability measure. Moreover, for some $q'>1$ we have
\begin{align*}
\E^{\bbQ_M}\Big[(\Xi^{v,\zeta}_T)^2\Big]&\leq C\E\big[|\mcE(\zeta^{M}*W)_T|^{q'}\big]^{1/q'}\E\Big[\big(\int_0^T(|\sigma(s,(X^{v,\zeta}_r)_{ r\leq s})|^2+|\sigma(s,(X^{v}_r)_{ r\leq s})|^2)ds\big)^{q}\Big]^{1/q}
\\
&\leq C
\end{align*}
by Assumption~\ref{ass:onSFDE}.\ref{ass:onSFDE-a-sigma}, where $\frac{1}{q'}+\frac{1}{q}=1$ and, by Lemma~\ref{lem:Lp-bound-if-martingale}, $C>0$ can be chosen independently from $M$. This gives that
\begin{enumerate}[a)]\setcounter{enumi}{3}
  \item For each $\epsilon>0$, there is a $M\geq 0$ such that $\bbQ_{M'}(\{\omega:\Xi^{v,\zeta}_T\leq M'\})>1-\epsilon$ for all $M'\geq M$.
\end{enumerate}

Making use of b) and d) we find that for each $\epsilon>0$, there is a $M\geq 0$ such that
\begin{align*}
1-\epsilon<\Prob(\{\omega:X^v\in C_M(T)\})&=\bbQ_M(\{\omega:X^{v,\zeta,M}\in C_M(T)\})% =
\end{align*}
and
\begin{align*}
\bbQ_M( \{\omega:\Xi^{v,\zeta}_T> M\})<\epsilon.
\end{align*}
%Moreover, we note that $X^{v,\zeta,M}=\tilde X^{v,\zeta^M}$ where $\tilde X^{v,\zeta^M}$ solves \eqref{ekv:driftless-sde1}-\eqref{ekv:driftless-sde2} with $W$ replaced by $W^{\zeta^M}:=W+\int_0^\cdot \zeta^M_sds$ which gives that $\tilde X^{v,\zeta^M}=X^{v,\zeta}$ on $ \{\omega:\Xi^{v,\zeta}_T\leq M\}$.
Combining these gives that
\begin{align*}
1-\epsilon&<\bbQ_M(\{\omega:X^{v,\zeta,M}\in C_M(T)\}\cap \{\omega:\Xi^{v,\zeta}_T\leq M\})+\bbQ_M(\{\omega:\Xi^{v,\zeta}_T> M\})
\\
&\leq \bbQ_M(\{\omega:X^{v,\zeta,M}\in C_M(T)\}\cap \{\omega:\Xi^{v,\zeta}_T\leq M\})+\epsilon
\end{align*}
Moreover, by property a) above and right-continuity we have that
\begin{align*}
\{\omega: X^{v,\zeta,M}\in C_M(t)\}\cap \{\omega:\Xi^{v,\zeta}_t\leq M\}=\{\omega: X^{v,\zeta}\in C_M(t)\}\cap \{\omega:\Xi^{v,\zeta}_t\leq M\}=E_M(t)
\end{align*}
so that
\begin{align*}
1-2\epsilon<\bbQ_M\big(E_M(T))=\E\big[\mcE(\zeta^{M}*W)_T \ett_{E_M(T)}\big]\leq \E\big[\mcE(\zeta*W)_T\big].
\end{align*}
Since $\epsilon>0$ was arbitrary, this proves the assertion whenever $v\in\mcU^k$ for some finite $k\geq 0$. To get the result for arbitrary $v:=(\tau_1,\ldots,\tau_N;\beta_1,\ldots,\beta_N)\in\mcU^f$, we define the sets
\begin{align*}
F_k(t):=\{\omega:|\zeta_t|\leq k_L(1+\sup_{s\in[0,t]}|X^{[v]_k}_s|)\}\supset \{\omega: N\leq k\}
\end{align*}
for all $k\geq 0$ and let
\begin{align*}
\zeta^{[v]_k}_t=\ett_{F_k(t)}\zeta_t.
\end{align*}
Now, for any $v\in\mcU^f$ we have by definition that $\Prob(\{\omega:N>k,\,\forall k\geq 0\})=0$ and so we can by again appealing to Lemma~\ref{lem:Lp-bound-if-martingale} find a $k\geq 0$ such that
\begin{align*}
\bbQ^{[v]_k}(N> k):=\E\big[\mcE(\zeta^{[v]_k}*W)_T \ett_{[N> k]}\big]\leq \E\big[|\mcE(\zeta^{[v]_k}*W)_T|^{q'}\big]^{1/q'}(\Prob(N> k))^{1/q}<\epsilon
\end{align*}
implying that
\begin{align*}
1-\epsilon<\bbQ^{[v]_k}(N \leq  k)=\E\big[\mcE(\zeta^{[v]_k}*W)_T\ett_{[N\leq k]}\big]=\E\big[\mcE(\zeta^{v}*W)_T\ett_{[N\leq k]}\big]\leq \E\big[\mcE(\zeta^{v}*W)_T\big]
\end{align*}
and the assertion follows as $\epsilon>0$ was arbitrary.\qed\\

\begin{cor}
There is a $q'>1$ such that $\sup_{\zeta\in\mcK_0}\E\big[|\mcE(\zeta*W)_T|^{q'}\big]<\infty$.
\end{cor}

\noindent\emph{Proof.} Lemma~\ref{lem:is-martingale} shows that for each $\varsigma>1$ the $\Upsilon^{v,\varsigma}$ in Lemma~\ref{lem:Lp-bound-if-martingale} is in fact all $\mcP_\bbF$-measurable processes $\zeta$ with $|\zeta_t|\leq k_L(1+\sup_{s\in[0,t]}|X^v_s|)$ for all $t\in[0,T]$ (outside of a $\Prob$-null set).\qed\\

The above corollary gives the following:

\begin{prop}\label{prop:SFDEmomentQ}
Under Assumption~\ref{ass:onSFDE}, the FSDE \eqref{ekv:forward-sde1}-\eqref{ekv:forward-sde2} admits a weak solution for each $(u,\alpha)\in\mcU\times \mcA$. Furthermore, the solution has moments of all orders on compacts, in particular we have for $p\geq 0$, that
\begin{align}\label{ekv:SFDEmomentQ}
\sup_{u\in\mcU^f,\bbQ\in\PrM^u}\E^{\bbQ}\Big[\sup_{t\in[0,T]}|X^{u}_t|^{p}\Big]\leq C,
\end{align}
where $C=C(p)$ and
\begin{align}\label{ekv:SFDEmomentQ2}
\sup_{v\in\mcU^f}\E\Big[\sup_{t\in [0,T]}\esssup_{u\in\mcU^f_t,\bbQ\in\PrM^{v(t)\circ u}}|\E^{\bbQ}\big[\sup_{s\in[t,T]}|X^{v(t)\circ u}_s|^{\rho}\big|\mcF_t\big]|^p\Big]\leq C,
\end{align}
where $C=C(\rho,p)$.

Moreover, there is a $q'>1$ and a $C>0$ such that for each $v\in\mcU^f$ and all $\zeta\in\mcK^v$ and $\bbQ\in\PrM^v$ we have $\E^{\bbQ}[|\mcE(\zeta*W^{\bbQ})_T|^{q'}]\leq C$ (where $W^\bbQ$ is a Brownian motion under $\bbQ$).
\end{prop}

\noindent\emph{Proof.} Existence of a weak solution to \eqref{ekv:forward-sde1}-\eqref{ekv:forward-sde2} follows by taking $\zeta_t=\breve a(t,(X^u_s)_{ s\leq t},\alpha_t)$ and using Lemma~\ref{lem:is-martingale} to conclude that under $\bbQ^{u,\alpha}$, the process $W^{u,\alpha}:=W-\int_0^\cdot \breve a(t,(X_s)_{ s\leq t},\alpha_t)dt$ is a Brownian motion.% Moreover, weak uniqueness follows by strong uniqueness of solutions to \eqref{ekv:driftless-sde1}-\eqref{ekv:driftless-sde1}

The moment estimates \eqref{ekv:SFDEmomentQ} and \eqref{ekv:SFDEmomentQ2} now follow by repeating the steps in the proof of Proposition~\ref{prop:SFDEmoment} and the last assertion follows by repeating the steps in Lemma~\ref{lem:Lp-bound-if-martingale} and Lemma~\ref{lem:is-martingale} while referring to the bound \eqref{ekv:SFDEmomentQ} rather than \eqref{ekv:SFDEmoment}.\qed\\

\subsection{The sequential system of reflected BSDEs\label{sec:robust-rbsde-FH}}
In the present section we show that there is a unique family $(Y^{v},Z^{v},K^{v})$ that solves the sequential system of reflected BSDEs
\begin{align}\label{ekv:seq-bsde-robust-FH}
\begin{cases}
  Y^{v}_t=\psi(X^v_T)+\int_t^T H^{*,v}(s,Z^{v}_s)ds-\int_t^T Z^{v}_sdW_s+ K^{v}_T-K^{v}_t,\quad\forall t\in[0,T], \\
  Y^{v}_t\geq\sup_{b\in U}\{Y^{v\circ(t,b)}_t-\ell^v(t,X^{v}_t,b))\},\quad\forall t\in[0,T],\\
  \int_0^T(Y^{v}_t-\sup_{b\in U}\{Y^{v\circ(t,b)}_t-\ell^v(t,X^{v}_t,b))\})dK^{v}_t=0,
\end{cases}
\end{align}
where $\ell^u(t,X^{u}_t,b)):=\infty\ett_{[0,\tau_N)}(t)+\ell(t,X^{u}_t,b))$ making \eqref{ekv:seq-bsde-robust-FH} a non-reflected BSDE on $[0,\tau_N)$. In the remainder of the article we will drop the subscript $v$ in $\ell^v$ but remind ourselves that no reflection can occur before the time of the last intervention in $v$.

Then, we will leverage the result in Theorem~\ref{thm-verification-FH} to find a weak solution to the robust impulse control problem in finite horizon.

Letting
\begin{align*}
\bar K^{v,p}:=2^{p-1}\big((T(C^g_\phi)^p +(C^g_\psi)^p+(k_L)^p)(1+\esssup_{u\in\mcU^f_t}\E\Big[\sup_{s\in[0,T]}|X^{v(t)\circ u}_s|^{\rho p}\Big|\mcF_t\Big]):0\leq t\leq T\big),
\end{align*}
we note that
\begin{align*}
\esssup_{u\in\mcU_f}\E\Big[|L^{v\circ u}_T|^p+|\psi(X^{v\circ u}_T)|^p+\int_t^T|H^{*,v\circ u}(s,0)|^pds\Big|\mcF_t\Big]\leq \bar K^{v,p}_t.
\end{align*}
%for all $\bbQ\in\PrM^L$.
Moreover, by Proposition~\ref{prop:SFDEmoment} we have that $\sup_{v\in\mcU^f}\|\bar K^{v,p}\|_{\mcS^2}<\infty$, for all $p\geq 1$ and by \eqref{ekv:moment_steg1} we have that $(\bar K^{v,p}:v\in\mcU^f,p\geq 1)$ satisfies the relation in \eqref{ekv:barK-Doob-equivalent}.

On the other hand $H^{*,v'}(t,z')-H^{*,v}(t,z)$ contains the term $\phi(t,X^{v'}_t)-\phi(t,X^{v}_t)$ and so $H^{*,v}$ generally fails to satisfy the conditions in Assumption~\ref{ass:on-coeff-FH} since $\phi$ is only locally Lipschitz in $x$. The same thing applies to $\psi$ and $\ell$ and we will rely on a localization argument leading us to introduce
\begin{align*}
H^{*,u,m,n}(t,\omega,z):=\inf_{\alpha\in A}H^{u,m,n}(t,\omega,z,\alpha),
\end{align*}
where
\begin{align*}
H^{u,m,n}(t,\omega,z,\alpha):=z\breve a(t,(X^{u}_s)_{ s\leq t},\alpha)+\phi^{m,n}(t,X^{u}_t,\alpha),%\phi^+(t,X^{u}_t,\alpha)\ett_{[|X^{u}_t|\leq m]}-\phi^-(t,X^{u}_t,\alpha)\ett_{[|X^{u}_t|\leq n]}%(H^{u}(t,\omega,z,\alpha))^+\ett_{[|X^{u}_t|\leq m]}-(H^{u}(t,\omega,z,\alpha))^-\ett_{[|X^{u}_t|\leq n]}.
\end{align*}
with $\phi^{m,n}:[0,T]\times \R^d\times A\to\R$ given by $\phi^{m,n}(t,x,\alpha):=\phi^{+,m}(t,x,\alpha)-\phi^{-,n}(t,x,\alpha)$, where $(\phi^{+,m})_{m\geq 0}$ and $(\phi^{-,n})_{n\geq 0}$ are both non-decreasing sequences of Borel-measurable, non-negative functions that are Lipschitz continuous in $x$ and continuous in $\alpha$ such that $\phi^{+,m}=\phi^+(x)$ on $|x|\leq m$ and $\phi^{-,n}(x)= \phi^-(x)$ on $|x|\leq n$.

Similarly, for $m,n\geq 0$, we let $\psi^{m,n}:\R^d\to\R$ be given by $\psi^{m,n}(x):=\psi^{+,m}(x)-\psi^{-,n}(x)$, where $(\psi^{+,m})_{m\geq 0}$ and $(\psi^{-,n})_{n\geq 0}$ are both non-decreasing sequences of non-negative, Lipschitz continuous functions such that $\psi^{+,m}=\psi^+(x)$ on $|x|\leq m$ and $\psi^{-,n}(x)= \psi^-(x)$ on $|x|\leq n$
% x\rfloor_m:=\frac{m}{|x|\vee m}x$ and set
%\begin{align*}
%\psi^{m,n}(x)&:=\psi^+(\lfloor x\rfloor_m)-\frac{|x|-m}{\inf\{r>0:\psi(\lfloor x\rfloor_{m+r})=0\}}\psi^+(\lfloor x\rfloor_m)\ett_{[|x|\leq m]}
%\\
%&\qquad-(\psi^-(x))^-\ett_{[|x|\leq n]}.
%\end{align*}
and let $\ell^{n}:[0,T]\times \R^d\times U\to [\delta,\infty)$ be a non-decreasing sequence of jointly continuous functions that are Lipschitz continuous in $x$ and H\"older continuous in $t$ (uniformly in the other variables) and satisfy $\ell^{n}(t,x)=\ell(t,x)$ on $|x|\leq n$ and $\ell^{n}(x)\geq \delta$ for all $(t,x)\in [0,T]\times\R^d$.

We now consider the following localized form of \eqref{ekv:seq-bsde-robust-FH}
\begin{align}\label{ekv:seq-bsde-robust-FH-local}
\begin{cases}
  Y^{v,m,n}_t=\psi^{m,n}(X^v_T)+\int_t^T H^{*,v,m,n}(s,Z^{v,m,n}_s)ds-\int_t^T Z^{v,m,n}_sdW_s+ K^{v,m,n}_T-K^{v,m,n}_t,\quad\forall t\in[0,T], \\
  Y^{v,m,n}_t\geq\sup_{b\in U}\{Y^{v\circ(t,b),m,n}_t-\ell^{n}(t,X^{v}_{t},b))\},\quad\forall t\in[0,T],\\
  \int_0^T(Y^{v,m,n}_t-\sup_{b\in U}\{Y^{v\circ(t,b),m,n}_t-\ell^{n}(t,X^{v}_{t},b))\})dK^{v,m,n}_t=0.
\end{cases}
\end{align}

%\begin{align}
%  \esssup_{u\in\mcU^f_t}\E\big[|L^{v\circ u}_T|^p\big|\mcF_t\big]&\leq \bar K^{v,p}_t,\label{ekv:L-bound}
%  \\
%  \esssup_{u\in\mcU^f_t}\E\big[|\phi^{\vecv,\vecv', u}_T|^p\big|\mcF_t\big]&\leq \bar K^{\vecv,\vecv',p}_t\label{ekv:phi-bound}
%  \\
%  \esssup_{u\in\mcU^f_t}\E\big[\sup_{s\in [t,T]}|\bar K^{v\circ u(s),p}_s|^{r}\big|\mcF_t\big]&\leq C_r\bar K^{v,pr}_t,\label{ekv:barK-Doob-equivalent}
%\end{align}

Since
\begin{align*}
|H^{*,v',m,n}(t,z')-H^{*,v,m,n}(t,z)|&\leq C((1+|L^{v'}_t|\vee|L^{v}_t|)|z'-z|
\\ &+(|z|+|z'|)\sup_{\alpha\in A}|\breve a(t,(X^{v'}_s)_{ s\leq t},\alpha)-\breve a(t,(X^{v}_s)_{ s\leq t},\alpha)|)
\\
&\quad+|\phi^{m,n}(t,X^{v'}_t)-\phi^{m,n}(t,X^{v}_t)|,
\end{align*}
Lemma~\ref{lem:SFDEflow} and Assumption~\ref{ass:onSFDE} implies the existence of a family $(\Lambda^{\vecv,\vecv',u}:(\vecv,\vecv')\in \cup_{\kappa\geq 1}\mcD^\kappa\times\mcD^\kappa,\,v\in\mcU^f)$ and a family $(\bar K^{\vecv,\vecv',k,p}:(\vecv,\vecv')\in \cup_{\kappa\geq 1}\mcD^\kappa\times\mcD^\kappa,k,p\geq 0)$ satisfying the conditions in Definition~\ref{defn:bounding-fam} and Assumption~\ref{ass:on-coeff-FH} and it follows by Proposition~\ref{prop:seq-existence} and Theorem~\ref{thm-verification-FH} that there is a unique family $(Y^{v,m,n},Z^{v,m,n},K^{v,m,n})$ that solves the sequential system of reflected BSDEs in \eqref{ekv:seq-bsde-robust-FH-local}.

We have,
\begin{lem}\label{lem:U-V-mn-bounds}
For $v,u\in\mcU^f$, let $(U^{v,u,m,n},V^{v,u,m,n})\in\mcS^2_l\times\mcH^2$ solve
\begin{align}\nonumber
U_t^{v,u,m,n}&=\psi^{m,n}(X^{v\circ u}_T)+\int_t^{T}H^{*,v\circ u,m,n}(s,V^{v,u,m,n}_s)ds
-\int_t^{T} V^{v,u,m,n}_sdW_s
\\
&\quad-\sum_{j=1}^N \ett_{[t\leq\tau_j]}\ell^{n}(\tau_j,X^{v\circ [u]_{j-1}}_{\tau_j},\beta_j),\label{ekv:U-V-mn-def}
\end{align}
whenever it has a unique solution and set $U^{v,u,m,n}\equiv -\infty$, otherwise. Then, there is a $C>0$ (that does not depend on $m,n$) such that, whenever $u^*\in\mcU^f$ is such that $\esssup_{u\in\mcU^f}U_t^{v,u,m,n}=U_t^{v,u^*,m,n}$, we have
\begin{align}\label{ekv:UV-mn-bound}
|U^{v,u^*,m,n}_t|^2+\E^{\bbQ}\Big[\int_t^T |V^{v,u^*,m,n}_s|^2ds+(N^*)^2\Big|\mcF_t\Big]\leq C(1+\esssup_{u\in\mcU^f_t,\zeta\in\mcK^{v\circ u}}\E^{\bbQ^\zeta}\Big[\sup_{s\in[t,T]}|X^{v\circ u}_s|^{2\rho}\Big|\mcF_t\Big]).
\end{align}
for all $\bbQ\in\PrM^{v\circ u^*}$.
\end{lem}

\noindent\emph{Proof.} First note that whenever $u^*\in\mcU^f$ is a maximizer then \eqref{ekv:U-V-mn-def} admits a unique solution. The bounds on $|U^{v,u^*,m,n}_t|^2$ and $\E^\bbQ\big[\int_t^T |V^{v,u^*,m,n}_s|^2ds\big|\mcF_t\big]$ now follow by repeating the argument in the proof of Lemma~\ref{lem:seq-U-V-bound} while noting that
\begin{align*}
\E^{\bbQ}\Big[|\psi^{m,n}(X^{v\circ u}_T)|^p+\int_t^{T}|H^{*,v\circ u,m,n}(s,0)|^pds\Big|\mcF_t\Big]\leq C\big(1+\E^\bbQ\Big[\sup_{s\in[0,T]}|X^{v\circ u}_s|^{p \rho}\Big|\mcF_t\Big]\big).
\end{align*}
From this, the bound on $N^*$ is immediate from \eqref{ekv:interv-cost-bound}.\qed\\

The statement of Lemma~\ref{lem:U-V-mn-bounds} holds for all $\bbQ\in\PrM^{v\circ u^*}$. Here it is notable that, since for any $m,n,m',n>0$ the drivers $H^{*,m,n}$ and $H^{*,m',n'}$ have the same stochastic Lipschitz coefficients, the set $\PrM^{v\circ u^*}$ is not parameterized by $m,n$. This is a key property when deriving the following stability result:

\begin{lem}\label{lem:Y-mn-unif-conv}
For each $m\geq 0$ and $p\geq 1$ we have
\begin{align}\label{ekv:mn-diff-to-0}
\|\sup_{b\in U}|Y^{v\circ(\cdot,b),m,n'}_\cdot-Y^{v\circ(\cdot,b),m,n}_\cdot|\|_{\mcS^p}\to 0
\end{align}
as $n,n'\to\infty$.
%Moreover, the limit $\bar Y^{v,m}:=\lim_{n\to\infty}Y^{v,m,n}$ and
%\begin{align*}
%\lim_{m,m'\to\infty}\lim_{n\to\infty}\{\|Y^{v,m',n}-Y^{v,m,n}\|_{\mcS^2}+\|Z^{v,m',n}-Z^{v,m,n}\|_{\mcH^2}\}=0.
%\end{align*}
\end{lem}

\noindent\emph{Proof.} %Let
%\begin{align*}
%U_t^{v,u,m,n}&=\psi^{m,n}(X^{v\circ u}_T)+\int_t^{T}H^{*,v\circ u,m,n}(s,V^{v,u,m,n}_s)ds
%-\int_t^{T} V^{v,u,m,n}_sdW_s
%\\
%&\quad-\sum_{j=1}^N \ett_{[t\leq\tau_j]}\ell^{n}(\tau_j,X^{v\circ [u]_{j-1}}_{\tau_j},\beta_j).
%\end{align*}
We have,
\begin{align*}
|Y^{v\circ(t,b),m,n'}_t-Y^{v\circ(t,b),m,n}_t|&\leq |U_t^{v\circ (t,b),u^*,m,n}-U_t^{v\circ (t,b),u^*,m,n'}|+|U_t^{v\circ (t,b),u',m,n}-U_t^{v\circ (t,b),u',m,n'}|,
\end{align*}
where $u^*$ and $u'$ are elements of $\mcU^f_t$ such that $\sup_{u\in\mcU^f}U_t^{v\circ (t,b),u,m,n}=U_t^{v\circ (t,b),u^*,m,n}$ and\\ $\sup_{u\in\mcU^f}U_t^{v\circ (t,b),u,m,n'}=U_t^{v\circ (t,b),u',m,n'}$. We now consider the first term and suppress the references to the control strategies (\ie $v\circ (t,b)$ and $u^*$) in the superscript, we have
\begin{align*}
U_t^{m,n}-U_t^{m,n'}&=(\psi^{m,n}(X_T)-\psi^{m,n'}(X_T))+\int_t^T (H^{*,m,n}(s,V^{m,n}_s)-H^{*,m,n'}(s,V^{m,n'}_s)ds
\\
&\quad -\int_t^{T} (V^{m,n}_s-V^{m,n'}_s)dW_s-\sum_{j=1}^{N^*}(\ell^{n}(\tau_j,X_{\tau_j},\beta_j)-\ell^{n'}(\tau_j,X_{\tau_j},\beta_j)).
\end{align*}
Taking the conditional expectation under the measure $\bbQ^{n,n'}$ where $W_t-\int_0^t\zeta^{n,n'}_sds$ is a martingale, with
\begin{align*}
\zeta^{n,n'}_s=\frac{H^{*,m,n'}(s,V^{m,n}_s)-H^{*,m,n'}(s,V^{m,n'}_s)}{|V^{m,n}_s -V^{m,n'}_s|^2} (V^{m,n}_s-V^{m,n'}_s)^\top\ett_{[V^{m,n}_s\neq V^{m,n'}_s]},
\end{align*}
gives
\begin{align*}
|U_t^{m,n}-U_t^{m,n'}|&\leq C\E^{\bbQ^{n,n'}}\Big[|\psi^{m,n}(X_T)-\psi^{m,n'}(X_T)|+\int_t^T |H^{*,m,n}(s,V^{m,n}_s)-H^{*,m,n'}(s,V^{m,n}_s)|ds
\\
&\quad +\sum_{j=1}^{N^*} |\ell^{n}(\tau_j,X_{\tau_j},\beta_j)-\ell^{n'}(\tau_j,X_{\tau_j},\beta_j)|\Big|\mcF_t\Big]
\\
&\leq C\E^{\bbQ^{n,n'}}\Big[\ett_{[|X_{T}|>\underline n]}(1+|X_T|^\rho)+\int_t^T \ett_{[|X_{s}|>\underline n]}(1+|X_s|^\rho)ds
\\
&\quad +\sum_{j=1}^{N^*} (1+|X_{\tau_j}|^\rho)\ett_{[|X_{\tau_j}|>\underline n]}\Big|\mcF_t\Big],
\end{align*}
where $\underline n:=n\wedge n'$. As $H^{*,m,n}$ and $H^{*,m,n'}$ have the same $z$-coefficient and thus also the same stochastic Lipschitz coefficient we find that $\bbQ^{n,n'}\in \PrM^{v\circ u^*}$. In particular, we have
\begin{align*}
|U_t^{m,n}-U_t^{m,n'}|&\leq C\E^{\bbQ^{n,n'}}\Big[\ett_{[\sup_{s\in[t,T]}|X_{s}|>\underline n]}(1+\sup_{s\in[t,T]}|X_{s}|^\rho)(1+N^*)\Big|\mcF_t\Big]
\\
&\leq C\E^{\bbQ^{n,n'}}\Big[\ett_{[\sup_{s\in[t,T]}|X_{s}|>\underline n]}(1+\sup_{s\in[t,T]}|X_{s}|^{2\rho})\Big|\mcF_t\Big]^{1/2}\E^{\bbQ^{n,n'}}\Big[1+(N^*)^2\Big|\mcF_t\Big]^{1/2}
\\
&\leq C\esssup_{u\in\mcU^f_t,\bbQ\in\PrM^{v\circ u}}\E^\bbQ\Big[\ett_{[\sup_{s\in[t,T]}|X^{v\circ u}_{s}|>\underline n]}\Big|\mcF_t\Big]^{1/4}(1
\\
&\quad+\esssup_{u\in\mcU^f_t,\bbQ\in\PrM^{v\circ u}}\E^\bbQ\Big[\sup_{s\in[t,T]}|X^{v\circ u}_{s}|^{4\rho})\Big|\mcF_t\Big]).
%\\
%&\leq C\esssup_{u\in\mcU^f_t,\alpha\in\mcA_t}\E\Big[\ett_{[\sup_{s\in[t,T]}|X^{v\circ u,\alpha}_{s}|>\underline n]}\Big|\mcF_t\Big]^{1/4}(1 +\E\Big[|\!X^{v}_{t}|^{4\rho})\Big|\mcF_t\Big]).
\end{align*}
Combining the above and applying the relation in \eqref{ekv:SFDE-moments-relation}, we find that for $p\geq 2$, we have
\begin{align*}
\|\sup_{b\in U}|Y^{v\circ(\cdot,b),m,n'}_\cdot-Y^{v\circ(\cdot,b),m,n}_\cdot|\|_{\mcS^p}^p&\leq C\E\Big[\sup_{t\in[0,T]}\esssup_{u\in\mcU^f_t,\bbQ\in\PrM^{v\circ u}}\E^\bbQ\Big[\ett_{[\sup_{s\in[t,T]}|X^{v\circ u}_{s}|>\underline n]}\Big|\mcF_t\Big]\Big]^{1/2}
\\
&\quad\times(1 +\esssup_{u\in\mcU^f,\bbQ\in\PrM^{u}}\E^\bbQ\Big[\sup_{t\in[0,T]}|X^{u}_{t}|^{4pq}\Big]^{1/2})
\end{align*}
where the first term tends to 0 as $\underline n\to\infty$ and the second term is bounded by Proposition~\ref{prop:SFDEmomentQ}.\qed\\

Now, as clearly $\|\sup_{b\in U}|\ell^{n'}(\cdot,X^v_\cdot,b)-\ell^{n}(\cdot,X^v_\cdot,b)|\|_{\mcS^p}\to 0$ as $n,n'\to\infty$ for all $p\geq 1$, Lemma~\ref{lem:Y-mn-unif-conv} and Proposition~\ref{prop:rbsde-solu-FH} implies that
\begin{align*}
\lim_{n,n'\to\infty}\{\|Y^{v,m,n'}-Y^{v,m,n}\|_{\mcS^2}+\|Z^{v,m,n'}-Z^{v,m,n}\|_{\mcH^2}+\|K^{v,m,n'}-K^{v,m,n}\|_{\mcS^2}\}=0.
\end{align*}

For each $m\geq 0$ we note that $(Y^{v,m,n})_{n\geq 0}$ is a non-increasing sequence of continuous processes that is $\Prob$-a.s.~bounded and we have that $Y^{v,m,n}$ converges pointwisely to a progressively measurable process $Y^{v,m}:=\lim_{n\to\infty}Y^{v,m,n}$. Furthermore, by Lemma~\ref{lem:Y-mn-unif-conv} we find that $Y^{v,m}$ is continuous and thus belongs to $\mcS^2$.

\begin{prop}\label{prop:seq-bsde-Ym}
For each $m\geq 0$, there is a family of pairs $(Z^{v,m},K^{v,m}:v\in\mcU^f)$ such that $(Y^{v,m},Z^{v,m},K^{v,m}:v\in\mcU^f)$ is the unique solution to
\begin{align}\label{ekv:seq-bsde-Ym}
\begin{cases}
  Y^{v,m}_t=\psi^{m}(X^v_T)+\int_t^T H^{*,v,m}(s,Z^{v,m}_s)ds-\int_t^T Z^{v,m}_sdW_s+ K^{v,m}_T-K^{v,m}_t,\quad\forall t\in[0,T], \\
  Y^{v,m}_t\geq\sup_{b\in U}\{Y^{v\circ(t,b),m}_t-\ell(t,X^{v}_{t},b))\},\quad\forall t\in[0,T],\\
  \int_0^T(Y^{v,m}_t-\sup_{b\in U}\{Y^{v\circ(t,b),m}_t-\ell(t,X^{v}_{t},b))\})dK^{v,m}_t=0,
\end{cases}
\end{align}
where $H^{*,u,m}(t,\omega,z):=\inf_{\alpha\in A}H^{u,m}(t,\omega,z,\alpha)$, with
\begin{align*}
H^{u,m}(t,\omega,z,\alpha):=z\breve a(t,(X^{u}_s)_{ s\leq t},\alpha)+\phi^{+,m}(t,X^{u}_t,\alpha)-\phi^{-}(t,X^{u}_t,\alpha),
\end{align*}
and $\psi^{m}(x):=\psi^{+,m}(x)-\psi^{-}(x)$.
\end{prop}

\noindent\emph{Proof.} Let $\eta_k:=\inf\{s\geq 0:|X^{v}_s|\geq k\}\wedge T$. Then, for all $n\geq k$ it follows that $(Y^{v,m,n},Z^{v,m,n},K^{v,m,n})$ solves
\begin{align*}
\begin{cases}
  Y^{v,m,n}_t=Y^{v,m,n}_{\eta_k}+\int_t^{\eta_k} H^{*,v,m}(s,Z^{v,m,n}_s)ds-\int_t^{\eta_k} Z^{v,m,n}_sdW_s+ K^{v,m,n}_{\eta_k}-K^{v,m,n}_t,\quad\forall t\in[0,{\eta_k}], \\
  Y^{v,m,n}_t\geq\sup_{b\in U}\{Y^{v\circ(t,b),m,n}_t-\ell(t,X^{v}_{t},b))\},\quad\forall t\in[0,{\eta_k}],\\
  \int_0^{\eta_k}(Y^{v,m,n}_t-\sup_{b\in U}\{Y^{v\circ(t,b),m,n}_t-\ell(t,X^{v}_{t},b))\})dK^{v,m,n}_t=0.
\end{cases}
\end{align*}
Moreover, by Proposition~\ref{prop:seq-existence} and Theorem~\ref{thm-verification-FH} there is a unique family of triples $(\hat Y^{v,m},\hat Z^{v,m},\hat K^{v,m})$ that solves
\begin{align*}
\begin{cases}
   \hat Y^{v,m}_t=Y^{v,m}_{\eta_k}+\int_t^{\eta_k} H^{*,v,m}(s,\hat Z^{v,m}_s)ds-\int_t^{\eta_k} \hat Z^{v,m}_sdW_s+ \hat K^{v,m}_{\eta_k}-\hat K^{v,m}_t,\quad\forall t\in[0,{\eta_k}], \\
  \hat Y^{v,m}_t\geq\sup_{b\in U}\{\hat Y^{v\circ(t,b),m}_t-\ell(t,X^{v}_{t},b))\},\quad\forall t\in[0,{\eta_k}],\\
  \int_0^{\eta_k}(\hat Y^{v,m}_t-\sup_{b\in U}\{\hat Y^{v\circ(t,b),m}_t-\ell(t,X^{v}_{t},b))\})d\hat K^{v,m}_t=0.
\end{cases}
\end{align*}
Letting, $n\to\infty$ we have by Lemma~\ref{lem:Y-mn-unif-conv} and Proposition~\ref{prop:rbsde-solu-FH} that
\begin{align*}
\|(\hat Y^{v,m} - Y^{v,m,n})\ett_{[0,\eta_k]}\|_{\mcS^2} + \|(\hat Z^{v,m} - Z^{v,m,n})\ett_{[0,\eta_k]}\|_{\mcH^2} + \|(\hat K^{v,m} - K^{v,m,n})\ett_{[0,\eta_k]}\|_{\mcS^2} \to 0
\end{align*}
and we find that there is a family of pairs $(Z^{v,m},K^{v,m}:v\in\mcU^f)$ such that for each $k\geq 0$,
\begin{align*}
\begin{cases}
   Y^{v,m}_t=Y^{v,m}_{\eta_k}+\int_t^{\eta_k} H^{*,v,m}(s,Z^{v,m}_s)ds-\int_t^{\eta_k} Z^{v,m}_sdW_s+ K^{v,m}_{\eta_k}-K^{v,m}_t,\quad\forall t\in[0,{\eta_k}], \\
  Y^{v,m}_t\geq\sup_{b\in U}\{Y^{v\circ(t,b),m}_t-\ell(t,X^{v}_{t},b))\},\quad\forall t\in[0,{\eta_k}],\\
  \int_0^{\eta_k}(Y^{v,m}_t-\sup_{b\in U}\{Y^{v\circ(t,b),m}_t-\ell(t,X^{v}_{t},b))\})dK^{v,m}_t=0.
\end{cases}
\end{align*}
Now, since there is a $\Prob$-a.s.~finite $k_0(\omega)$ such that $\eta_k=T$ for all $k\geq k_0$, existence of a solution to \eqref{ekv:seq-bsde-Ym} follows.

Uniqueness is established by repeating steps in the proof of Theorem~\ref{thm-verification-FH}.\qed\\

By an identical argument to that used in Lemma~\ref{lem:Y-mn-unif-conv} we find that
\begin{align*}%\label{ekv:m-diff-to-0}
\lim_{m,m'\to\infty}\|\sup_{b\in U}|Y^{v\circ(\cdot,b),m'}_\cdot-Y^{v\circ(\cdot,b),m}_\cdot|\|_{\mcS^p}= 0
\end{align*}
for all $p\geq 1$ and we conclude that:

\begin{prop}
The sequential system of reflected BSDEs \eqref{ekv:seq-bsde-robust-FH} has a unique solution.
\end{prop}

\noindent\emph{Proof.} The result follows by letting $m\to\infty$ and using an identical argument to that of Proposition~\ref{prop:seq-bsde-Ym}.\qed\\

%%%%%%%%%%%%%%%%%%%%%%%%%%%%%%%%%%%%%%%%%%%%%%%%%%%%%%%%%%%%%%%%%%%%%%%%%%%%%%%%%%%%%%%%%%%%%%%%%%%%%%%%%%%%%%%%%%%%%%%%%%%%%%%%%%%%%%%%%%%%%%%%%%

\subsection{Robust impulse control in finite horizon}
We are now ready to solve the robust impulse control problem by relating optimal controls to solutions of the sequential system of reflected BSDEs \eqref{ekv:seq-bsde-robust-FH}. However, before doing this we need to narrow down the set of admissible impulse controls that we search over in order to guarantee that \eqref{ekv:U-V-mn-def} admits a unique solution.

\begin{defn}
We let $\mcU^m$ be the set subset of $\mcU^f$ with all $u=(\tau_1,\ldots,\tau_N;\beta_1,\ldots,\beta_N)$ such that $N$ has moments of all orders, \ie for each $k\geq 0$ we have $\E\big[(N)^k\big]<\infty$.
\end{defn}

\begin{lem}\label{lem-optimal-in-mcU-m}
We have
\begin{align*}%\label{ekv:mcU-m-contain-optimal}
\sup_{u\in\mcU^f}\inf_{\alpha\in\mcA}J(u,\alpha)=\sup_{u\in\mcU^m}\inf_{\alpha\in\mcA}J(u,\alpha).
\end{align*}
\end{lem}

\noindent\emph{Proof.} For $(u,\alpha)\in \mcU^f\times\mcA$ we let $(B^{u,\alpha},E^{u,\alpha})$ solve
\begin{align}\label{ekv:B-E-def}
B^{u,\alpha}_t=\psi(X^{u,\alpha}_T)+\int_t^T\phi(s,X^{u,\alpha}_s,\alpha_s)ds-\int_t^TE^{u,\alpha}_sdW_s - \sum_{j=1}^n\ett_{[t\leq \tau_j]}\ell(\tau_j,X^{[u]_{j-1},\alpha}_{\tau_j},\beta_j),
\end{align}
whenever a solution exists in $\mcS^2_l\times\mcH^2$. We define the set of sensible impulse controls, $\mcU^s$, as the subset of $u\in\mcU^f$ such that for each $t\in[0,T]$ and $\alpha\in\mcA$,
\begin{align}\nonumber
&\E\Big[\psi(X^{u,\alpha}_T)+\int_t^T\phi(s,X^{u,\alpha}_s,\alpha_s)ds - \sum_{t\leq \tau_j}\ell(\tau_j,X^{[u]_{j-1},\alpha}_{\tau_j},\beta_j)\Big|\mcF_t\Big]
\\
&\geq -(C_\psi^g+C_\phi^g T)(1+\esssup_{u'\in\mcU^f_t,\alpha'\in\mcA_t}\E\Big[\sup_{s\in [t,T]}|X^{u(t-)\circ u',\alpha(t-)\circ \alpha'}_s|^\rho\Big|\mcF_t\Big]),\label{ekv:sensible-control}
\end{align}
where $u(t-)=(\tau_1,\ldots,\tau_{N(t-)};\beta_1,\ldots,\beta_{N(t-)})$ with $N(t-):=\max\{j\geq 0:\tau_j<t\}$ and $\alpha(t-)\circ \alpha':=\ett_{[0,t)}\alpha+\ett_{[t,T]}\alpha'$. Then for each $u\in\mcU^f$ we obtain a $u'\in\mcU^s$ by removing all future interventions whenever \eqref{ekv:sensible-control} does not hold. Moreover, $u'$ dominates $u$ in the sense that $J(u',\alpha)>J(u,\alpha)$ for all $\alpha\in\mcA$. In particular, we note that
\begin{align*}%\label{ekv:mcU-m-contain-optimal}
\sup_{u\in\mcU^f}\inf_{\alpha\in\mcA}J(u,\alpha)=\sup_{u\in\mcU^s}\inf_{\alpha\in\mcA}J(u,\alpha).
\end{align*}
Now, whenever $u\in\mcU^s$ there is a $(B^{u,\alpha},E^{u,\alpha})\in\mcS^2_l\times\mcH^2$ that solves \eqref{ekv:B-E-def}. We will build on the argument in Lemma~\ref{lem:seq-U-V-bound} to show that $\mcU^s\subset\mcU^m$. We thus assume that $u\in\mcU^s$. Rearranging the terms in \eqref{ekv:B-E-def} gives
\begin{align}\label{ekv:B-interv-cost}
\sum_{j=1}^{N} \ell(\tau_j,X^{[u]_{j-1},\alpha}_{\tau_j},\beta_j)=\psi(X^{u,\alpha}_T)-B^{u,\alpha}_0 +\int_0^{T}\phi(s,X^{u,\alpha}_s,\alpha_s)ds-\int_0^{T}E^{u,\alpha}_sdW_s,
\end{align}
where we know that all terms on the right hand side, except for the last (martingale) term, have moments of all orders. By Ito's formula we have
\begin{align*}
|B^{u,\alpha}_0|^2+\int_0^{T}|E^{u,\alpha}_s|^2ds&=|\psi(X^{u,\alpha}_T)|^2+2\int_0^{T}B^{u,\alpha}_s\phi(s,X^{u,\alpha}_s,\alpha_s)ds -2\int_0^TB^{u,\alpha}_sE^{u,\alpha}_sdW_s
\\
&\quad+\sum_{j=1}^{N} (-2B^{u,\alpha}_{\tau_j}\ell(\tau_j,X^{[u]_{j-1},\alpha}_{\tau_j},\beta_j) +|\ell(\tau_j,X^{[u]_{j-1},\alpha}_{\tau_j},\beta_j)|^2)
\\
&\leq |\psi(X^{u,\alpha}_T)|^2+(1+2\kappa)\sup_{t\in[0,T]}|B^{u,\alpha}_t|^2+\int_0^{T}|\phi(s,X^{u,\alpha}_s,\alpha_s)|^2ds
\\
&\quad-2\int_0^TB^{u,\alpha}_sE^{u,\alpha}_sdW_s+\frac{2}{\kappa}\big(\sum_{j=1}^{N} \ell(\tau_j,X^{[u]_{j-1},\alpha}_{\tau_j},\beta_j)\big)^2,
\end{align*}
for $\kappa>0$ (where we have used that $\ell(\tau_j,X^{[u]_{j-1},\alpha}_{\tau_j},\beta_j)\leq 2\sup_{t\in[0,T]}|B^{u,\alpha}_t|$). Using \eqref{ekv:B-interv-cost}, the growth conditions on $\phi$ and $\psi$ and the fact that $u\in\mcU^s$ together with \eqref{ekv:SFDE-moments-relation} gives that
\begin{align*}
\int_0^{T}|E^{u,\alpha}_s|^2ds&\leq C(1+\kappa+\frac{1}{k})(1+\sup_{t\in [0,T]}|X^{u,\alpha}|^{2\rho})-2\int_0^TB^{u,\alpha}_sE^{u,\alpha}_sdW_s+\frac{4}{\kappa}\Big|\int_0^{T}E^{u,\alpha}_sdW_s\Big|^2.
\end{align*}
Raising both sides to $p\geq 1$ followed by taking the expectation and applying BDG gives
\begin{align*}
\E\Big[\big(\int_0^{T}|E^{u,\alpha}_s|^2ds\big)^p\Big]&\leq C(1+\kappa+\frac{1}{k})+C\E\Big[\big(\int_0^T|B^{u,\alpha}_sE^{u,\alpha}_s|^2ds\big)^{p/2}+\frac{1}{\kappa}\big(\int_0^{T}|E^{u,\alpha}_s|^2\big)^p\Big]
\\
&\leq C(1+\kappa+\frac{1}{k})+\frac{C}{\kappa}\E\Big[\big(\int_0^{T}|E^{u,\alpha}_s|^2\big)^p\Big]
\end{align*}
where we have used the relation $2ab\leq \kappa a^2+\frac{1}{\kappa}b^2$ together with the bound on $\E\big[\sup_{t\in [0,T]}|B^{u,\alpha}_t|^{2p}\big]$ resulting from the fact that $u\in\mcU^s$ to reach the last inequality. Now, choosing $\kappa$ sufficiently large it follows that the left hand side is finite. Finally, as the left hand side of \eqref{ekv:B-interv-cost} is greater that $\delta N$ we conclude that $u\in\mcU^m$.\qed\\

By Benes’ selection Theorem (\cite{Benes1971}, Lemma 5, pp. 460), there exists, for each $v\in\mcU^f$, a
$\mcP\otimes\mcB(\R^d)/\mcB(A)$-measurable function $\alpha^v(t,\omega,z)$ such that for any given $(t,\omega, z) \in [0,T]\times\Omega\times \R^{d}$, we have
\begin{align*}
H^v(t,\omega,z,\alpha^v(t,\omega,z))=\inf_{\alpha\in A}H^v(t,\omega,z,\alpha),
\end{align*}
$\Prob$-a.s.

The following theorem shows that we can extract the optimal pair $(u^*,\alpha^*)$ from the family of maps $(\alpha^v(t,\omega,z):v\in\mcU^f)$ and the solution to \eqref{ekv:seq-bsde-robust-FH}.
\begin{thm}\label{thm:is-robust-FH}
Let the family $(Y^v,Z^v,K^v:v\in\mcU^f)$ be a solution to \eqref{ekv:seq-bsde-robust-FH}. Then the pair $(u^*,\alpha^*)\in\mcU^f\times \mcA$, with $u^*=(\tau^*_1,\ldots,\tau^*_{N^*};\beta^*_1,\ldots,\beta^*_{N^*})\in\mcU^f_0$ defined as:
\begin{itemize}
  \item $\tau^*_{j}:=\inf \Big\{s \geq \tau^*_{j-1}:\:Y_s^{v\circ[u^*]_{j-1}}=\sup_{b\in U} \{Y^{v\circ[u^*]_{j-1}\circ(s,b)}_s-\ell(s,X^{v\circ[u^*]_{j-1}}_s,b)\}\Big\}\wedge T$
  \item $\beta^*_j\in\mathop{\arg\max}_{b\in U}\{Y^{v\circ [u^*]_{j-1}\circ(\tau^*_j,b)}_{\tau^*_j}-\ell(\tau^*_j,X^{v\circ[u^*]_{j-1}}_{\tau^*_j},b)\}$
\end{itemize}
and $N^*=\sup\{j:\tau^*_j<T\}$, with $\tau_0^*:=0$ and
\begin{align*}
\alpha^*_t:=\sum_{j=0}^{N^*}\ett_{[\tau^*_j,\tau^*_{j+1})}(t)\alpha^{[u^*]_j}(t,Z^{[u^*]_j}_t),
\end{align*}
with $\tau^*_{N+1}:=\infty$ is an optimal pair in the sense that
\begin{align}\label{ekv:is-robust-FH}
Y^\emptyset_0=J(u^*,\alpha^*)=\sup_{u\in\mcU^f}\inf_{\alpha\in\mcA}J(u,\alpha).
\end{align}
\end{thm}

\noindent\emph{Proof.} For $u\in\mcU^m$ we let $(U^u,V^u)\in \mcS^2_l\times\mcH^2$ be the unique solution to
\begin{align}\label{ekv:U-V-def-robust}
U_t^{u}&=\psi(X^u_T)+\int_t^T H^{*,u}(s,V^u_s)ds
-\int_t^{T} V^{u}_sdW_s-\sum_{j=1}^N \ett_{[\tau_j\geq t]}\ell(\tau_j,X^{[u]_{j-1}}_{\tau_j},\beta_j).
\end{align}
Then, by Theorem~\ref{thm:bsde-no-reflection}, $\|V^u\|_{\mcH^p}<\infty$ for all $p\geq 1$ implying that $V^u\in\mcH^2_{\bbQ}$ for all $\bbQ\in\PrM_0$ and, since $U^u_0$ is $\mcF_0$-measurable and $\mcF_0$ is trivial, we have
\begin{align*}
U^u_0&=\E^{\bbQ^{u,\alpha^*}}\Big[\psi(X^u_T)+\int_0^T H^{*,u}(s,V^u_s)ds
-\int_0^{T} V^{u}_sdW_s-\sum_{j=1}^N \ell(\tau_j,X^{[u]_{j-1}}_{\tau_j},\beta_j)\Big]
\\
&=\E^{\bbQ^{u,\alpha^*}}\Big[\psi(X^u_T)+\int_0^T \phi(s,X^u_s,\alpha^*_s)ds
-\int_0^{T} V^{u}_sdW^{u,\alpha^*}_s-\sum_{j=1}^N \ell(\tau_j,X^{[u]_{j-1}}_{\tau_j},\beta_j)\Big]
\\
&= J(u,\alpha^*)
\end{align*}
where now $\bbQ^{u,\alpha}$ is the measure, equivalent to $\Prob$, under which $W^{u,\alpha}:=W-\int_0^\cdot \breve a(s,(X^v_r)_{ r\leq s},\alpha_s)ds$ is a martingale. Moreover, by Theorem~\ref{thm-verification-FH} and Lemma~\ref{lem-optimal-in-mcU-m} we have
\begin{align*}
Y^\emptyset_0=\sup_{u\in\mcU^f}U^u_0=\sup_{u\in\mcU^m}J(u,\alpha^*)=\sup_{u\in\mcU^f}J(u,\alpha^*).
\end{align*}
To show that $\alpha^*$ is an optimal response it is, in light of Lemma~\ref{lem-optimal-in-mcU-m}, enough to show that $\alpha^*$ is a minimizer for all $u\in\mcU^m$. However, for any $(u,\alpha)\in\mcU^m\times\mcA$, we have
\begin{align*}
U^u_0&=\E^{\bbQ^{u,\alpha}}\Big[\psi(X^u_T)+\int_0^T H^{*,u}(s,V^u_s)ds
-\int_0^{T} V^{u}_sdW_s-\sum_{j=1}^N \ell(\tau_j,X^{[u]_{j-1}}_{\tau_j},\beta_j)\Big]
\\
&=\E^{\bbQ^{u,\alpha}}\Big[\psi(X^u_T)+\int_0^T \phi(s,X^u_s,\alpha_s)ds
-\int_0^{T} V^{u}_sdW^{u,\alpha}_s-\sum_{j=1}^N \ell(\tau_j,X^{[u]_{j-1}}_{\tau_j},\beta_j)\Big]
\\
&\quad +\E^{\bbQ^{u,\alpha}}\Big[\int_0^T (H^{*,u}(s,V^{u}_s)-H^{u}(s,V^{u}_s,\alpha_s))ds\Big]
\\
&\leq J(u,\alpha)
\end{align*}
and \eqref{ekv:is-robust-FH} follows.\qed\\

%%%%%%%%%%%%%%%%%%%%%%%%%%%%%%%%%%%%%%%%%%%%%%%%%%%%%%%%%%%%%%%%%%%%%%%%%%%%%%%%%%%%%%%%%%%%%%%%%%%%%%%%%%%%%%%%%%%%%%%%%%%%%%%%%%%%%%%%%%%%%%%%%%

\bibliographystyle{plain}
\bibliography{RBSDEimp-FH_ref}

\begin{thebibliography}{10}

\bibitem{AgramOksen}
N.~Agram and B.~{\O}ksendal.
\newblock Stochastic control of memory mean-field processes.
\newblock {\em Appl. Math. Optim.}, 79:181--204, 2019.

\bibitem{BaseiImpulse}
M.~Basei.
\newblock Optimal price management in retail energy markets: an impulse control
  problem with asymptotic estimates.
\newblock {\em Math Meth Oper Res}, 89:355--383, 2019.

\bibitem{BayraktarRobust}
E.~Bayraktar, A.~Cosso, and H.~Pham.
\newblock Robust feedback switching control: dynamic programming and viscosity
  solutions.
\newblock {\em SIAM J. Control Optim.}, 54(5):2594--2628, 2016.

\bibitem{Bender00}
C.~Bender and M.~Kohlmann.
\newblock Bsdes with stochastic lipschitz condition.
\newblock {\em CoFE Discussion Paper, No. 00/08, University of Konstanz, Center
  of Finance and Econometrics (CoFE), Konstanz}, 2000.

\bibitem{Benes1971}
V.~E. Benes.
\newblock Existence of optimal stochastic control laws.
\newblock {\em SIAM Journal on Control}, 9(3):446--472, 1971.

\bibitem{BensLionsImpulse}
A.~Bensoussan and J.L. Lions.
\newblock {\em Impulse Control and Quasivariational inequalities}.
\newblock Gauthier-Villars, Montrouge, France, 1984.

\bibitem{BertsekasShreve}
D.~P. Bertsekas and S.~E. Shreve.
\newblock {\em Stochastic optimal control: The discrete-time case}.
\newblock Academic Press, 1978.

\bibitem{Briand08}
P.~Briand and F.~Confortola.
\newblock Bsdes with stochastic lipschitz condition and quadratic pdes in
  hilbert spaces.
\newblock {\em Stochastic Process. Appl.}, 118:818--838, 2008.

\bibitem{CarmLud}
R.~Carmona and M.~Ludkovski.
\newblock Pricing asset scheduling flexibility using optimal switching.
\newblock {\em Appl. Math. Finance}, 15:405--447, 2008.

\bibitem{CohenElliottBook}
S.~N. Cohen and R.~J. Elliott.
\newblock {\em Stochastic Calculus and Applications}.
\newblock Birkh\"auser, New York, NY, 2 edition, 2015.

\bibitem{DjehiceImpulse}
B.~Djehiche, S.~Hamad\'ene, and I.~Hdhiri.
\newblock Stochastic impulse control of non-markovian processes.
\newblock {\em Appl Math Optim}, 61(1):1--26, 2010.

\bibitem{DjehicheInfHorImp}
B.~Djehiche, S.~Hamad\'ene, I.~Hdhiri, and H.~Zaatra.
\newblock Infinite horizon stochastic impulse control with delay and random
  coefficients.
\newblock {\em arXiv:1904.11924}, 2019.

\bibitem{ElAsri2020}
B.~{El Asri}, S.~Hamad{\'e}ne, and K.~Oufdil.
\newblock On the stochastic control-stopping problem.
\newblock {\em arXiv:2005.06789}, 2020.

\bibitem{ElKaroui1}
N.~El-Karoui, C.~Kapoudjian, E.~Pardoux, S.~Peng, and M.~C. Quenez.
\newblock Reflected solutions of backward {SDEs} and related obstacle problems
  for {PDEs}.
\newblock {\em Ann. Probab.}, 25(2):702--737, 1997.

\bibitem{Girsanov1960}
I.~V. Girsanov.
\newblock On transforming a certain class of stochastic processes by absolutely
  continuous substitution of measures.
\newblock {\em Theory of Probability and its Applications}, 5(3):285--301,
  1960.

\bibitem{HamJean}
S.~Hamad{\'e}ne and M.~Jeanblanc.
\newblock On the starting and stopping problem: application in reversible
  investments.
\newblock {\em Math. Oper. Res.}, 32(1):182--192, 2007.

\bibitem{HamZhang}
S.~Hamad\'ene and J.~Zhang.
\newblock Switching problem and related system of reflected backward {SDEs}.
\newblock {\em Stochastic Processes and their Applications}, 120(4):403--426,
  2010.

\bibitem{Hdhiri}
I.~Hdhiri and M.~Karouf.
\newblock Optimal stochastic impulse control with random coefficients and
  execution delay.
\newblock {\em Stochastics}, 90(2):151--164, 2018.

\bibitem{HuTang}
Y.~Hu and S.~Tang.
\newblock Multi-dimensional {BSDE} with oblique reflection and optimal
  switching.
\newblock {\em Prob. Theory and Related Fields}, 147(1-2):89--121, 2008.

\bibitem{Imkeller}
P.~Imkeller and G.~{Dos Reis}.
\newblock Path regularity and explicit convergence rate for bsdewith truncated
  quadratic growth.
\newblock {\em Stochastic Processes and their Applications}, 120(3):348--379,
  2010.

\bibitem{JonteSFDE}
J.~J{\"o}nsson and M.~Perninge.
\newblock Finite horizon impulse control of stochastic functional differential
  equations.
\newblock {\em arXiv:2006.09768}, 2020.

\bibitem{ElKarouiTan}
N.~El Karoui and X.~Tan.
\newblock Capacities, measurable selection and dynamic programming part i:
  Abstract framework.
\newblock {\em arXiv:1310.3363}, 2013.

\bibitem{KolmogorovFomin}
A.~N. Kolmogorov and S.~V. Fomin.
\newblock {\em Introductory Real Analysis}.
\newblock Dover Publications Inc., 2000.

\bibitem{Korn}
R.~Korn.
\newblock Some applications of impulse control in mathematical finance.
\newblock {\em Math Meth Oper Res}, 50:493--518, 1999.

\bibitem{OksenSulemBok}
B.~{\O}ksendal and A.~Sulem.
\newblock {\em Applied Stochastic Control of Jump Diffusions}.
\newblock Springer, 2007.

\bibitem{PSImpulsive}
J.~Palczewski and L.~Stettner.
\newblock Impulsive control of portfolios.
\newblock {\em Appl Math Optim}, 56:67--103, 2007.

\bibitem{SwitchElephant}
M.~Perninge.
\newblock A finite horizon optimal switching problem with memory and
  application to controlled sddes.
\newblock {\em Math Meth Oper Res}, 2019.

\bibitem{ImpulsIH}
M.~Perninge.
\newblock Infinite horizon impulse control of stochastic functional
  differential equations.
\newblock {\em arXiv:2003.08833}, 2020.

\bibitem{Protter}
P.~Protter.
\newblock {\em Stochastic Integration and Differential Equations}.
\newblock Springer, Berlin, 2nd edition, 2004.

\end{thebibliography}
\end{document}